\documentclass[reqno]{amsart}
\usepackage[english]{babel}
\usepackage{amssymb,amsmath,mathtools,cite,booktabs}
\usepackage{enumerate,xcolor}
\usepackage[T1]{fontenc}
\usepackage[a4paper,margin=1in]{geometry}
\usepackage[colorlinks=true,linkcolor=blue,citecolor=blue,urlcolor=blue]{hyperref}
\usepackage{comment}

\newtheorem{theorem}{Theorem}[section]
\newtheorem{proposition}[theorem]{Proposition}
\newtheorem{lemma}[theorem]{Lemma}
\newtheorem{corollary}[theorem]{Corollary}
\theoremstyle{definition}
\newtheorem{definition}[theorem]{Definition}
\theoremstyle{remark}
\newtheorem{remark}[theorem]{Remark}
\numberwithin{equation}{section}

\newcommand{\C}{\mathbb C}
\newcommand{\D}{\mathbb D}
\newcommand{\T}{\mathbb T}
\newcommand{\B}{\mathbb B}
\newcommand{\Acal}{\mathcal A}
\newcommand{\Bcal}{\mathcal B}
\newcommand{\norm}[1]{\left\lVert #1\right\rVert}
\newcommand{\abs}[1]{\left\lvert #1\right\rvert}
\newcommand{\ip}[2]{\left\langle #1,#2\right\rangle}
\newcommand{\dd}{\,d}
\DeclareMathOperator{\Ree}{Re}
\DeclareMathOperator{\diag}{diag}
\DeclareMathOperator{\cond}{cond}

\DeclareMathOperator{\alg}{alg}

\title[The Complete Crouzeix Conjecture in Dimension Three and the COR conjecture]{The Complete Crouzeix Conjecture in Dimension Three and the Clou\^atre--Ostermann--Ransford conjecture}

\author{Per \AA hag}
\address{Department of Mathematics and Mathematical Statistics\\
Ume\aa\ University\\SE-901~87 Ume\aa\\ Sweden}
\email{per.ahag@math.umu.se}

\author{Rafa\l\ Czy{\.z}}
\address{Faculty of Mathematics and Computer Science, Jagiellonian
University, \L ojasiewicza~6, 30-348 Krak\'ow, Poland}
\email{rafal.czyz@im.uj.edu.pl}

\author{Jani Virtanen}
\address{University of Eastern Finland, Department of Physics and
Mathematics, 80100 Joensuu, Finland. University of Reading, Reading,
UK\@. University of Helsinki, Helsinki, Finland}
\email{jani.virtanen@uef.fi, j.a.virtanen@reading.ac.uk, and \newline
jani.virtanen@helsinki.fi}

\begin{document}

\begin{abstract}
We settle the complete Crouzeix conjecture for matrices of order at most three and prove the Q-algebra case of the completely bounded Clou\^atre--Ostermann--Ransford (COR) conjecture for homomorphisms into matrices of order at most three, via sharp abstract column and row estimates. Our approach also establishes the scalar COR conjecture in a stronger form, for homomorphisms with commutative range on Banach algebras with unity satisfying von Neumann's inequality. Under contractivity of the symmetrized map, this result holds on arbitrary Hilbert spaces without initial boundedness assumptions on the homomorphism or the antilinear map. For arbitrary operator algebras, we disprove the complete COR conjecture by an exact three-dimensional example with target matrix order two. We also prove the sharp complete bound for every matrix subalgebra containing the diagonal, in arbitrary matrix order and on arbitrary target Hilbert spaces. We also obtain column and row square-function inequalities with sharp norm bounds, strict scalar bounds for operators similar to normal operators, sharp complete spectral constants for scaled $q$-numerical ranges in dimensions two and three, and rigidity, stability, and representing-measure results.

\end{abstract}

\maketitle

\section{Introduction}\label{sec:introduction}

In 2004, Crouzeix \cite{Crouzeix2004} conjectured that for $A\in M_n(\C)$
\begin{equation}\label{eq:intro-crouzeix}
        \norm{p(A)}
        \leq 2\max_{z\in W(A)}\abs{p(z)}
\end{equation}
for every complex polynomial $p$, where
$$
    W(A)=\{\ip{Ax}{x}:\norm x=1\}
$$
is the numerical range of $A$. As evidence for his conjecture, the paper included a proof of the $2\times 2$ case. It is easy to see that the constant $2$ is sharp already in dimension two. 

In 2007, Crouzeix~\cite{Crouzeix2007} obtained the first dimension-free bound, replacing $2$ in \eqref{eq:intro-crouzeix} by $11.08$. In the same paper, he introduced the completely bounded version: For matrix-valued polynomials $F=(f_{ij})$, with $F(A)=(f_{ij}(A))$, he conjectured that
\begin{equation}\label{eq:intro-crouzeix-cb}
        \norm{F(A)}
        \leq 2\max_{z\in W(A)}\norm{F(z)},
\end{equation}
and proved the corresponding inequality with $2$ replaced by $11.08$. We refer to~\eqref{eq:intro-crouzeix-cb} as the \emph{complete Crouzeix conjecture}, which is equivalent to asking whether $W(A)$ is a complete $2$-spectral set for $A$.

The spectral-set approach to numerical ranges was initiated by Delyon and Delyon~\cite{DelyonDelyon1999}, and subsequently used by Crouzeix and Palencia~\cite{CrouzeixPalencia2017} to lower the bound in \eqref{eq:intro-crouzeix-cb} to $1+\sqrt2$. The constant $1+\sqrt2$ is optimal for the abstract lemma used in the Crouzeix--Palencia argument~\cite{RansfordSchwenninger2018}, although for every fixed numerical-range domain a strictly smaller scalar bound is available~\cite[Theorems~2 and~3]{MalmanEtAl2025}. Recently, Jin~\cite{Jin2026} and Lorist and Schwenninger~\cite{LoristSchwenninger2026} independently proved the scalar conjecture. Both papers explicitly note that their methods do not directly yield the completely bounded inequality. Indeed, in the Lorist--Schwenninger argument, the obstruction is a commutativity hypothesis that is not available after matrix amplification. For \eqref{eq:intro-crouzeix-cb}, $1+\sqrt2$ therefore remains the best universal bound.

Our first main result settles the complete conjecture for matrices of order at most three.

\begin{theorem}\label{thm:intro-complete-three}
Let $A\in M_n(\C)$ with $n\leq3$. Then
\begin{equation}\label{eq:intro-complete-three}
    \norm{F(A)}\leq2\max_{z\in W(A)}\norm{F(z)}
\end{equation}
for every matrix-valued polynomial $F$ of arbitrary size. Thus $W(A)$ is a complete $2$-spectral set for $A$.

More generally, if \eqref{eq:intro-crouzeix-cb} fails for some $A\in M_n(\C)$, then it already fails for an $s\times r$ matrix-valued polynomial with
\begin{equation}\label{eq:intro-rs}
    r,s\geq2,\qquad r+s\leq n;
\end{equation}
in particular $n\geq4$, and in dimension four the only possible matrix level is $2\times2$.
\end{theorem}

The case $n=2$ of the complete inequality was already known \cite[Section~5.2]{BadeaCrouzeixDelyon2006} while the new result is the case $n=3$.

The proof of Theorem~\ref{thm:intro-complete-three} combines new square-function inequalities, which control scalar, column, and row matrix levels, with a finite-dimensional reduction that restricts the size of a possible counterexample.

The same finite-dimensional reduction, combined with abstract analogs of the column and row estimates, also yields a complete abstract result. Clou\^atre, Ostermann, and Ransford formulated abstract scalar and completely bounded versions of Crouzeix's conjecture for symmetrized homomorphisms \cite{ClouatreOstermannRansford2023}. We prove their complete conjecture through dimension
three for Q-algebras, that is, quotients of uniform algebras equipped with the quotient matrix norms; see Definition~\ref{v3:def-Q}. In the scalar setting, Theorem~\ref{v3:scalar-automatic}
extends the uniform-algebra theorem of Badea, O'Loughlin,
and Virtanen~\cite[Theorem~2.1]{BadeaOloughlinVirtanen2026}
to homomorphisms with commutative range on Banach algebras
with unity satisfying von Neumann's inequality.
Under contractivity of the symmetrized map, it removes
the initial boundedness assumptions on the homomorphism
and the antilinear map.

\begin{theorem}\label{thm:intro-abstract-complete-three}
Let $\Bcal$ be a Q-algebra with unity, let $\mathcal H$ be a Hilbert space,
let $\theta:\Bcal\to B(\mathcal H)$ be a completely bounded homomorphism preserving the identity, and let $\alpha:\Bcal\to\Bcal$ be an antilinear map preserving the identity. Define
\[
 \Phi_{\theta,\alpha}(b)
 =\frac12\bigl(\theta(b)+\theta(\alpha(b))^*\bigr),
 \qquad b\in\Bcal.
\]
Suppose that $\Phi_{\theta,\alpha}$ is completely contractive. If
$b_1,\ldots,b_q\in\Bcal$ satisfy
\[
        \norm{(b_1,\ldots,b_q)^{\mathsf T}}\leq1,
\]
then
\begin{equation}\label{eq:intro-abstract-column-row}
 \norm{\begin{pmatrix}
 \theta(b_1)\\ \vdots\\ \theta(b_q)
 \end{pmatrix}}
 \leq2,
 \qquad
 \norm{\begin{pmatrix}
 \theta(b_1)&\cdots&\theta(b_q)
 \end{pmatrix}}
 \leq2.
\end{equation}
Consequently, if $\mathcal H=\C^d$ and $d\leq3$, then
\[
        \norm{\theta}_{\mathrm{cb}}\leq2.
\]
In this finite-dimensional assertion, $\theta$ need not initially be
assumed bounded. The constant $2$ is best possible.
\end{theorem}

Notice that no boundedness assumption on $\alpha$ is needed in
Theorem~\ref{thm:intro-abstract-complete-three}.
The complete COR conjecture was formulated for arbitrary operator
algebras, with $\alpha$ completely contractive
\cite[Conjecture~6.1]{ClouatreOstermannRansford2023}.
In this generality it is false. We use the same notation
$\Phi_{\theta,\alpha}$ for the symmetrized map on an arbitrary operator
algebra, and write $T_d$ for the upper triangular subalgebra of
$M_d(\C)$. Complete contractivity of an antilinear map refers to its
entrywise amplifications.

\begin{theorem}\label{v5:thm:cor-counterexample}
There are a three-dimensional noncommutative operator algebra
$\Acal\subset M_3(\C)$, an injective completely bounded homomorphism
$\theta:\Acal\to M_2(\C)$ preserving the identity, and an antilinear map
$\alpha:\Acal\to\Acal$ preserving the identity such that
\[
 \norm{\alpha}_{\mathrm{cb}}
 =\norm{\Phi_{\theta,\alpha}}_{\mathrm{cb}}=1,
\]
but
\begin{equation}\label{v5:eq:cor-bound}
 2<\frac{2(48+\sqrt{2665})}{19(3+\sqrt5)}
 \leq\norm{\theta}=\norm{\theta}_{\mathrm{cb}}
 \leq\frac{8\sqrt{29}}{19}.
\end{equation}
Moreover, $\theta(\Acal)=T_2$ and $\Acal\cap\Acal^*=\C I_3$.
\end{theorem}

 The target matrix order and the dimension of the domain algebra are the smallest
possible, by Corollary~\ref{v5:cor:minimal}.
The homomorphism in Theorem~\ref{v5:thm:cor-counterexample} has
noncommutative range $T_2$, whereas the numerical-range functional
calculus of a single matrix has commutative range.

The complete bound nevertheless holds for every matrix subalgebra
containing the diagonal. These are the digraph algebras
\cite[p.~78]{Donsig1998}.

\begin{theorem}\label{v5:thm:digraph}
Let $\Bcal\subset M_n(\C)$ be a subalgebra containing every diagonal
matrix. Let $\theta:\Bcal\to B(\mathcal H)$ be a homomorphism preserving
the identity, and let $\alpha:\Bcal\to\Bcal$ be an antilinear contraction
preserving the identity. If $\Phi_{\theta,\alpha}$ is completely
contractive, then
\[
 \norm{\theta}_{\mathrm{cb}}\leq2.
\]
The constant $2$ is best possible.
\end{theorem}

There is no restriction on $n$ or on the target Hilbert space, and only
ordinary contractivity of $\alpha$ is assumed. The theorem includes
upper triangular matrix algebras, block triangular matrix algebras,
and all finite block incidence algebras. In particular, the bound holds
for $T_2$ with its usual matrix norms, but fails for the algebra in
Theorem~\ref{v5:thm:cor-counterexample}, which is algebraically
isomorphic to $T_2$. Thus the bound depends on the concrete
operator-algebra structure, not only on the algebraic isomorphism class.
\par

We next state the square-function inequalities underlying Theorem~\ref{thm:intro-complete-three}. 
Let $K_A=\overline{W(A)}$, and denote by $\Acal(K_A)$ the functions continuous on $K_A$ and holomorphic in its interior. For $k\ge 1$, let $\mathcal W_k=\{1, \ldots, q\}^k$ and $\mathcal W_0=\{\emptyset\}$. Given $f_1,\ldots,f_q\in\Acal(K_A)$ and $\alpha=(\alpha_1,\ldots,\alpha_k)\in \mathcal W_k$, we set $f_\alpha=f_{\alpha_1}\cdots f_{\alpha_k}$ and $f_\emptyset=1$. For self-adjoint operators $X$ and $Y$, we write $X\preceq Y$ when $Y-X$ is positive.

\begin{theorem}\label{thm:tuple-square}
Let $A\in\Bcal(H)$ and let $f_1,\ldots,f_q\in\Acal(K_A)$ satisfy
\begin{equation}\label{eq:intro-pointwise-column}
        \sum_{j=1}^q\abs{f_j(z)}^2\leq1
        \qquad(z\in K_A).
\end{equation}
For every $0<\tau<1$, the series
\begin{align}
 P_\tau(f,A)
 &=\sum_{k=0}^{\infty}\tau^{2k}
   \sum_{\alpha\in\mathcal W_k}
   f_\alpha(A)^*f_\alpha(A),
   \label{eq:intro-column-square}\\
 Q_\tau(f,A)
 &=\sum_{k=0}^{\infty}\tau^{2k}
   \sum_{\alpha\in\mathcal W_k}
   f_\alpha(A)f_\alpha(A)^*
   \label{eq:intro-row-square}
\end{align}
converge in operator norm and satisfy
\begin{equation}\label{eq:intro-square-bounds}
        I\preceq P_\tau(f,A)\preceq\frac1{1-\tau}I,
        \qquad
        I\preceq Q_\tau(f,A)\preceq\frac1{1-\tau}I.
\end{equation}
In particular,
\begin{equation}\label{eq:intro-column-row-two}
 \norm{(f_1(A),\ldots,f_q(A))^{\mathsf T}}\leq2,
 \qquad
 \norm{(f_1(A),\ldots,f_q(A))}\leq2.
\end{equation}
The constant $2$ is best possible.
\end{theorem}

Choosing $q=1$ and $\tau=1/2$, the preceding theorem recovers and extends \eqref{eq:intro-crouzeix} to every $f\in\Acal(K_A)$. The terms with $k\geq2$ and the other values of $\tau$ will be used later to obtain rigidity and stability results.

The proof of Theorem~\ref{thm:tuple-square} combines the positive operator-valued double-layer calculus of Delyon--Delyon and Crouzeix--Palencia \cite{DelyonDelyon1999, CrouzeixPalencia2017, SchwenningerDeVries2025} with a multivariable version of Jin's sampling argument~\cite{Jin2026}. Jin's proof produces the scalar Gramian estimate by sampling a matrix-valued Herglotz kernel on the disk. We carry out the sampling on the unit ball of $\C^q$, which produces the two tuple square functions above.

The passage from Theorem~\ref{thm:tuple-square} to Theorem~\ref{thm:intro-complete-three} is finite-dimensional. A rectangular form of the finite-level argument associated with Smith's lemma \cite{Smith1983}, together with a variational lemma of Hartz and McCarthy \cite{HartzMcCarthy2026}, shows that if the complete functional-calculus norm for $A\in M_d(\C)$ exceeds $2$, then this is already detected at an $s\times r$ matrix level with $r,s\geq2$ and $r+s\leq d$. The inequalities $r,s\geq2$ follow from \eqref{eq:intro-column-row-two}. Thus no such pair exists when $d\leq3$, while for $d=4$ only $r=s=2$ remains. This dimension count is the reason that the tuple square functions settle the complete conjecture through dimension three.

\medskip

For the remaining statements, define
$$
    \psi(A)=\sup\{\norm{f(A)}: f\in\Acal(K_A),\ \max_{z\in K_A}\abs{f(z)}\leq1\}
$$
and let $\psi_{\mathrm{cb}}(A)$ be the corresponding supremum over matrix-valued functions of arbitrary size. Thus, $\psi(A)\leq\psi_{\mathrm{cb}}(A)$, and Theorem~\ref{thm:intro-complete-three} gives $\psi_{\mathrm{cb}}(A)\leq2$ when $d\leq3$.

For $d\geq1$, put
\[
 Q_{\mathrm{cb}}(d)
 =\sup\{\psi_{\mathrm{cb}}(A):A\in M_d(\C)\}.
\]
Thus $Q_{\mathrm{cb}}(d)$ is the least constant in the complete
numerical-range inequality for all matrices of order $d$.
\par

The higher-order terms in the square-function estimates yield the following sharpening of Theorem~\ref{thm:intro-complete-three}.

\begin{theorem}\label{thm:intro-complete-equality}
Let $A\in M_d(\C)$ with $d\leq3$. Then
$$
    \psi_{\mathrm{cb}}(A)=2 \iff \psi(A)=2.
$$
\end{theorem}

The square-function estimates also improve the elementary bound obtained from a diagonalization. For an invertible matrix $V$, write
$\cond(V)=\norm V\norm{V^{-1}}$, and set
\begin{equation}\label{eq:intro-Phi}
        \Phi(\kappa)
        =\left(
        \frac{8\kappa^2}
        {\kappa^2+1+\sqrt{(\kappa^2-1)^2+4}}
        \right)^{1/2},
        \qquad \kappa\geq1.
\end{equation}

\begin{theorem}\label{thm:intro-complete-diagonalizable}
Let $d\leq3$ and let $A\in M_d(\C)$ be diagonalizable, say $A=V\Lambda V^{-1}$ with $\Lambda$ diagonal. Then
\begin{equation}\label{eq:intro-complete-diagonalizable}
        \psi_{\mathrm{cb}}(A)
        \leq
        \min\{\cond(V),\Phi(\cond(V))\}
        <2.
\end{equation}
\end{theorem}

The scalar counterpart holds in every dimension whenever $A$ is similar to a normal operator. Indeed, every diagonalizable matrix satisfies the scalar Crouzeix inequality strictly; conversely, if $A=V\Lambda V^{-1}$ and $\psi(A)=s<2$, then
$$
        \cond(V)^4>
        \frac{s^4}{4(4-s^2)}.
$$
Thus, if $A_j=V_j\Lambda_jV_j^{-1}$ and $\psi(A_j)\to2$, then $\cond(V_j)\to\infty$ for every choice of eigenvector matrices $V_j$.

This condition-number viewpoint is closely related to recent work of Aouichaoui and O'Loughlin~\cite{AouichaouiOloughlin2026}, which characterizes scaled $q$-numerical ranges as unions of numerical ranges under similarities with bounded condition number and determines their sharp scalar spectral constants. We obtain a completely bounded analog of their similarity-amplification result; see Theorem~\ref{thm:cb-similarity-amplification}. For $A\in M_d(\C)$ with $d\geq2$ and $0<\abs q\leq1$, let
$$
    W_q(A)=\{\ip{Ax}{y}:\norm x=\norm y=1,\ \ip{x}{y}=q\},
    \qquad\Omega_q(A)=q^{-1}W_q(A).
$$
Theorem~\ref{thm:q-complete-constants} shows that, for $d=2,3$, the sharp universal complete spectral constant for $\Omega_q(A)$ is $\max\{1,2\abs q/(1+(1-\abs q^2)^{1/2})\}$. Thus, $\Omega_q(A)$ is a complete spectral set for every matrix of order two or three if and only if $\abs q\leq4/5$; in every dimension at least two, this holds whenever $|q|\leq1/\sqrt2$.

The freedom in choosing $\tau$ and the higher-order terms in Theorem~\ref{thm:tuple-square} also give information at and near the sharp constant. If an operator $T$ satisfies the two square-function estimates at $\tau=1/2$ and $\norm T=2$ is attained, then $T$ has a two-dimensional reducing subspace on which it is unitarily equivalent to $\left(\begin{smallmatrix}0&0\\2&0\end{smallmatrix}\right)$. There is a quantitative version: Let $A\in M_n(\C)$. If $f\in\Acal(K_A)$ satisfies
$$
    \max_{z\in K_A}\abs{f(z)}=1,\qquad\norm{f(A)}\geq2-\varepsilon,
    \qquad 0\leq\varepsilon\leq2-\frac{\sqrt{15}}2,
$$
then $f(A)$ is within $8\sqrt\varepsilon$ of an operator having this block as a reducing summand. Thus near equality forces a concrete approximate nilpotent structure.

There is also a function-theoretic form of this rigidity. If $f$ maximizes $\norm{f(A)}$ among functions bounded by one on a convex domain $\Omega$ and $x$ is a unit vector with $\norm{f(A)x}=\norm{f(A)}$, then the associated representing probability measure $\mu$ on $\partial\Omega$ is characterized by
$$
    \ip{h(A)x}{x}=\int_{\partial\Omega}h\dd\mu \qquad(h\in\Acal(\overline\Omega));
$$
for its existence and uniqueness, see \cite{BickelEtAl2020}. We show that near equality, all nonzero Fourier coefficients of the pushforward measure $f_{\ast}\mu$ under $f$ are quantitatively small. At equality, $f_{\ast}\mu$ is normalized arclength measure on the unit circle, and along sequences approaching equality these measures converge weakly to it.

Jin's scalar rigidity theorem shows that $\psi(A)=2$ forces $W(A)$
to be a nondegenerate disk \cite[Theorem~1]{JinRigidity2026}.
Together with Theorem~\ref{thm:intro-complete-equality} and Crouzeix's
disk classification \cite[Lemma~10.2]{Crouzeix2016}, this gives the full
complete equality classification in dimension three;
see Corollary~\ref{v3:cor:three-classification}. In particular, a
nondisk numerical range has complete constant strictly less than two
in dimensions at most three. We also classify scalar polynomial
equality in every finite dimension in
Theorem~\ref{v3:thm:polynomial-classification}. After translating the
disk centre to zero, every such polynomial is a monomial, and the
matrix has the corresponding reducing Crabb--Choi--Crouzeix block.
Consequently, a polynomial of degree at least $d$ with $p(A)\ne0$
satisfies the scalar inequality strictly for $A\in M_d(\C)$.

Finally, we apply the general theory to the family
\begin{equation}\label{eq:intro-family}
    A_{\mu,r}=\begin{pmatrix}
    i\mu&0&ir\\
    0&-i\mu&ir\\
    ir&ir&-1
\end{pmatrix},
    \qquad\mu>0,\qquad r\geq0.
\end{equation}
Up to multiplication by $-1$ and unitary equivalence, this is a real family studied numerically by Crouzeix \cite{Crouzeix2016} while related numerical work appears in
\cite{GreenbaumOverton2018, CrouzeixGreenbaumLi2024}. Since
$$
    \Ree A_{\mu,r}=\diag(0,0,-1),
$$
the boundary of its numerical range contains $[-i\mu,i\mu]$. For $r>0$, the matrix is unitarily irreducible. Since $\mu>0$, the nontrivial boundary segment $[-i\mu,i\mu]$ also shows that its numerical range is not a disk. 

Corollary~\ref{v3:cor:nondisk-strict} therefore gives
$\psi_{\mathrm{cb}}(A_{\mu,r})<2$ for every $\mu>0$ and $r\geq0$.
We retain the quantitative bound obtained from an eigenvector matrix
at diagonalizable parameters. The nondiagonalizable parameters form
the curve
\[
    \mu=t\sqrt{1-2t},\qquad r=(1-t)\sqrt t,
    \qquad 0<t<\frac12.
\]
At $t=1/3$, the matrix has Jordan form $J_3(-1/3)$ and is affinely
equivalent to the $3\times3$ Kac--Murdock--Szeg\H{o} matrix. We retain
an independent proof of strictness there using the rational-domain
construction of Crouzeix, Greenbaum, and Li
\cite{CrouzeixGreenbaumLi2024}, with containment certified in exact
rational arithmetic by Sturm sequences \cite{BasuPollackRoy2006}.
Although every point of the defective curve satisfies strictness,
$\psi_{\mathrm{cb}}(A_{\mu(t),r(t)})\to2$ as $t\to1/2^-$.
Thus the strict bounds cannot be uniform along this curve.

The paper is organized as follows. Section~\ref{sec:quotient} sets up
the quotient operator algebra of the numerical-range functional calculus.
Section~\ref{sec:square-functions} proves the tuple square-function
inequalities. Sections~\ref{sec:stein} and~\ref{sec:rectangular} develop
the quantitative similarity estimates and the rectangular extremal
argument. Section~\ref{sec:complete-three} proves
Theorems~\ref{thm:intro-complete-three},
\ref{thm:intro-complete-equality}, and
\ref{thm:intro-complete-diagonalizable}.
Section~\ref{sec:abstract-cor} proves the scalar abstract estimates and
the Q-algebra extension. Section~\ref{v5:sec:noncommutative} gives the
counterexample to the complete COR conjecture and proves the bound for
digraph algebras. Section~\ref{sec:q-range} contains the scaled
$q$-numerical-range results. Sections~\ref{sec:equality}
and~\ref{sec:measures} treat equality, stability, and representing
measures, and Section~\ref{sec:dimension-three} classifies scalar
polynomial equality and complete equality in dimension three. Finally,
Section~\ref{sec:family} studies the family \eqref{eq:intro-family}.
\par

\section{The numerical-range functional calculus}\label{sec:quotient}

For $K\subset\C$ compact and convex, give $M_{s,r}(\Acal(K))$ the
matrix norm
\[
        \norm{F}_{K}=\max_{z\in K}\norm{F(z)}.
\]
These matrix norms make $\Acal(K)$ an operator algebra with unity. Let
$A\in M_d(\C)$ and put $K=W(A)$. 
Polynomials are dense in $\Acal(K)$ by Mergelyan's theorem, including
when $K$ has empty interior. The Crouzeix--Palencia inequality
\cite{CrouzeixPalencia2017}, applied to differences of polynomials,
therefore defines $f(A)$ by continuous extension of the polynomial
calculus. This agrees with the holomorphic functional calculus whenever
the latter applies.
\par
The functional calculus is the homomorphism preserving the identity
\[
        \rho_A:\Acal(K)\longrightarrow M_d(\C),
        \qquad \rho_A(f)=f(A).
\]
Let
\[
        \mathcal I_A=\ker\rho_A,
        \qquad
        \mathfrak A_A=\Acal(K)/\mathcal I_A.
\]
We equip $\mathfrak A_A$ with the quotient matrix norms. Thus, for
$F\in M_{s,r}(\Acal(K))$,
\begin{equation}\label{eq:quotient-matrix-norm}
 \norm{F+M_{s,r}(\mathcal I_A)}
 =\inf\{\norm{F+G}_{K}:G\in M_{s,r}(\mathcal I_A)\}.
\end{equation}
The Crouzeix--Palencia inequality \cite{CrouzeixPalencia2017} shows that
$\rho_A$ is bounded. Since $\rho_A(1)=I_d$, its kernel $\mathcal I_A$
is a proper closed two-sided ideal. By
\cite[Corollary~3.2]{BlecherRuanSinclair1990}, the quotient
$\mathfrak A_A$ is an operator algebra with unity. When a concrete
representation is needed, we use an identity-preserving completely
isometric representation given by
\cite[Theorem~3.1]{BlecherRuanSinclair1990}.
\par

The homomorphism $\rho_A$ induces an injective homomorphism preserving the identity
\[
        \theta_A:\mathfrak A_A\longrightarrow M_d(\C),
        \qquad \theta_A(f+\mathcal I_A)=f(A).
\]

\begin{proposition}\label{prop:quotient-calculus}
Let $A\in M_d(\C)$. Then
\begin{equation}\label{eq:quotient-image}
        \theta_A(\mathfrak A_A)=\alg(A),
        \qquad
        \dim\mathfrak A_A=\deg m_A\leq d,
\end{equation}
where $m_A$ is the minimal polynomial of $A$. Moreover,
\begin{equation}\label{eq:quotient-norms}
        \norm{\theta_A}=\psi(A),
        \qquad
        \norm{\theta_A}_{\mathrm{cb}}=\psi_{\mathrm{cb}}(A).
\end{equation}
The same equality holds at every fixed rectangular level.
\end{proposition}

\begin{proof}
The image of the polynomial calculus is $\alg(A)$. Since this algebra
is finite-dimensional, it is closed. Polynomial approximation on $K$
therefore gives
$\rho_A(\Acal(K))=\alg(A)$. The first isomorphism theorem gives the
first equality in \eqref{eq:quotient-image}. A monogenic matrix algebra
has basis $I,A,\ldots,A^{m-1}$, where $m=\deg m_A$. This proves the
dimension formula.

Fix $r,s$ and let $c_{s,r}(A)$ be the norm of the rectangular
functional calculus on $M_{s,r}(\Acal(K))$. For
$X=F+M_{s,r}(\mathcal I_A)$, the definition of the quotient norm gives
\[
 \norm{\theta_A^{(s,r)}(X)}
 =\norm{F(A)}
 \leq c_{s,r}(A)\norm{F+G}_{K}
 \qquad(G\in M_{s,r}(\mathcal I_A)).
\]
Taking the infimum over $G$ shows
$\norm{\theta_A^{(s,r)}}\leq c_{s,r}(A)$. Conversely,
\[
 \norm{F+M_{s,r}(\mathcal I_A)}\leq\norm{F}_{K},
\]
so
$c_{s,r}(A)\leq\norm{\theta_A^{(s,r)}}$. Thus, equality holds at each
rectangular level. Taking $r=s=1$ and then the supremum over all
levels proves \eqref{eq:quotient-norms}. It is immaterial whether the
suprema are taken over functions in $\Acal(K)$ or over polynomials,
because entrywise polynomial approximation is uniform at every fixed
matrix level.
\end{proof}

We shall repeatedly use approximate lifts from the quotient.

\begin{lemma}\label{lem:quotient-lifts}
Let $X\in M_{s,r}(\mathfrak A_A)$ and suppose that $\norm{X}\leq1$.
For every $\varepsilon>0$, there is
$F\in M_{s,r}(\Acal(W(A)))$ such that
\begin{equation}\label{eq:quotient-lift}
        F+M_{s,r}(\mathcal I_A)=X,
        \qquad
        \norm{F}_{W(A)}\leq1+\varepsilon.
\end{equation}
\end{lemma}

\begin{proof}
This is the definition of the quotient matrix norm in
\eqref{eq:quotient-matrix-norm}.
\end{proof}

The scalar and complete ratios are unchanged by taking adjoints.

\begin{lemma}\label{lem:adjoint-symmetry}
For $A\in M_d(\C)$,
\[
        \psi(A^*)=\psi(A),
        \qquad
        \psi_{\mathrm{cb}}(A^*)=\psi_{\mathrm{cb}}(A).
\]
\end{lemma}

\begin{proof}
One has $W(A^*)=\{\overline z:z\in W(A)\}$. For a scalar function $f$, set
$f^\#(z)=\overline{f(\overline z)}$. For a matrix-valued function
$F=(f_{ij})$, set
\[
        F^\#(z)=\bigl(\overline{f_{ji}(\overline z)}\bigr).
\]
Then $F^\#(A^*)=F(A)^*$ and
$\max_{W(A^*)}\norm{F^\#}=\max_{W(A)}\norm{F}$. The assertions follow
in both directions.
\end{proof}

\section{Positive kernels and iterated square functions}\label{sec:square-functions}

For $q\geq1$, let
\[
        \B_q=\{w\in\C^q:\norm{w}<1\}.
\]
We write $w^{\mathsf T}a=\sum_{j=1}^q w_ja_j$. This expression is
bilinear.

\begin{definition}
A scalar kernel $K$ on $\B_q$ is positive if
\[
        [K(w_i,w_j)]_{i,j=1}^N\succeq0
\]
for every finite family $w_1,\ldots,w_N\in\B_q$. An
$M_d(\C)$-valued kernel is positive when the corresponding block
matrices are positive.
\end{definition}

\begin{lemma}\label{lem:ball-kernel}
Let $a\in\overline{\B_q}$ and put
\[
        c_a(w)=\frac{1+w^{\mathsf T}a}{1-w^{\mathsf T}a}.
\]
Then
\begin{equation}\label{eq:ball-kernel}
        K_a(w,z)=\frac{c_a(w)+\overline{c_a(z)}}{1-z^*w}
\end{equation}
is a positive kernel on $\B_q$.
\end{lemma}

\begin{proof}
A direct calculation gives
\begin{equation}\label{eq:ball-factor}
 K_a(w,z)=
 \frac{2}{(1-w^{\mathsf T}a)(1-\overline{z^{\mathsf T}a})}
 \left(
 1+\frac{z^*\{I-\overline a(\overline a)^*\}w}{1-z^*w}
 \right).
\end{equation}
The matrix $I-\overline a(\overline a)^*$ is positive. If $M\succeq0$,
then
\[
        \frac{z^*Mw}{1-z^*w}
        =\sum_{k=0}^{\infty}(z^*Mw)(z^*w)^k
\]
is a positive kernel. Hence the expression in parentheses in
\eqref{eq:ball-factor} is positive. Multiplication by a scalar function
of $w$ and its conjugate at $z$ preserves positivity.
\end{proof}

We next give the double-layer formula in the form used below. This is
the positive operator-valued double-layer calculus of Delyon--Delyon
and Crouzeix--Palencia
\cite{DelyonDelyon1999,CrouzeixPalencia2017}. A detailed recent account
is given in \cite{SchwenningerDeVries2025}.

\begin{lemma}\label{lem:double-layer}
Let $B\in M_d(\C)$, and let $\Omega$ be a bounded convex domain with
smooth boundary such that $W(B)\subset\Omega$. There is a positive map preserving the identity
\[
        \Phi:C(\partial\Omega)\longrightarrow M_d(\C)
\]
with the following property. If $h$ is holomorphic on a neighborhood
of $\overline\Omega$, then
\begin{equation}\label{eq:double-layer-split}
        2\Phi(h)=h(B)+C_h^*,
        \qquad C_h\in\alg(B).
\end{equation}
\end{lemma}

\begin{proof}
Let $\nu(\zeta)$ be the outward unit normal to $\partial\Omega$, viewed
as a complex number, and let $s$ be arclength. Set
\begin{equation}\label{eq:double-layer-density}
 \mu(\zeta,B)=\frac{1}{2\pi}\left(
 \nu(\zeta)(\zeta I-B)^{-1}
 +\overline{\nu(\zeta)}(\overline\zeta I-B^*)^{-1}\right).
\end{equation}
Put $Z=\zeta I-B$ and $R=Z^{-1}$. Convexity gives
\[
        \nu Z^*+\overline\nu Z\succeq0.
\]
Indeed, its quadratic form at a unit vector $x$ is
\[
 2\operatorname{Re}\left\{
 \overline\nu(\zeta-x^*Bx)\right\}\geq0.
\]
Congruence by $R$ gives
$\nu R+\overline\nu R^*\succeq0$. Thus
$\mu(\zeta,B)\succeq0$.

For the positive orientation of $\partial\Omega$, one has
$\nu\dd s=-i\dd\zeta$. The Cauchy formula and its adjoint give
\[
        \int_{\partial\Omega}\mu(\zeta,B)\dd s=2I.
\]
Define
\[
        \Phi(g)=\frac12\int_{\partial\Omega}
        g(\zeta)\mu(\zeta,B)\dd s.
\]
Then $\Phi$ is positive and $\Phi(1)=I$. Splitting the two terms in
\eqref{eq:double-layer-density} gives \eqref{eq:double-layer-split},
where
\[
 C_h=\frac{1}{2\pi i}\int_{\partial\Omega}
 \overline{h(\zeta)}(\zeta I-B)^{-1}\dd\zeta.
\]
Every resolvent of $B$ belongs to $\alg(B)$. Hence $C_h\in\alg(B)$.
\end{proof}

Let $F=(f_1,\ldots,f_q)^{\mathsf T}$, where the $f_j$ are holomorphic
near $\overline\Omega$ and
\begin{equation}\label{eq:F-bound-domain}
        \sum_{j=1}^q\abs{f_j(z)}^2\leq1
        \qquad(z\in\overline\Omega).
\end{equation}
For $w\in\B_q$, define
\begin{equation}\label{eq:h-w}
        h_w(\zeta)=
        \frac{1+w^{\mathsf T}F(\zeta)}
             {1-w^{\mathsf T}F(\zeta)},
        \qquad H(w)=\Phi(h_w).
\end{equation}
The denominator does not vanish. The function $H$ is holomorphic and
$H(0)=I$.

\begin{lemma}\label{lem:H-kernel}
The kernel
\begin{equation}\label{eq:H-kernel}
        \mathcal L_H(w,z)
        =\frac{H(w)+H(z)^*}{1-z^*w}
\end{equation}
is positive on $\B_q$.
\end{lemma}

\begin{proof}
For $\zeta\in\partial\Omega$, put $a=F(\zeta)$. By
\eqref{eq:F-bound-domain}, $a\in\overline{\B_q}$. For finite points
$w_1,\ldots,w_N$, Lemma~\ref{lem:ball-kernel} and
\eqref{eq:double-layer-density} show that
\[
 [K_{F(\zeta)}(w_i,w_j)\mu(\zeta,B)]_{i,j=1}^N\succeq0.
\]
Integration gives
\[
 [\mathcal L_H(w_i,w_j)]_{i,j=1}^N
 =\frac12\int_{\partial\Omega}
 [K_{F(\zeta)}(w_i,w_j)\mu(\zeta,B)]_{i,j=1}^N\dd s.
\]
\end{proof}

Put
\[
        T_j=f_j(B),
        \qquad T(w)=\sum_{j=1}^q w_jT_j.
\]
The double-layer formula applied to $h_w$ gives
\begin{equation}\label{eq:completion}
        H(w)-(I-T(w))^{-1}\in\alg(B^*).
\end{equation}
Indeed,
\[
 h_w(B)=(I+T(w))(I-T(w))^{-1}
       =2(I-T(w))^{-1}-I,
\]
and \eqref{eq:completion} follows from
\eqref{eq:double-layer-split}. The sample points and the additional
point at the origin used below are the multivariable form of the
one-variable sampling argument in \cite{Jin2026}.

\begin{proposition}\label{prop:simple-square}
Let $B\in M_d(\C)$ have simple spectrum. Let $\Omega$ be a bounded
smooth convex domain with $W(B)\subset\Omega$. Suppose that
$f_1,\ldots,f_q$ are holomorphic near $\overline\Omega$ and satisfy
\eqref{eq:F-bound-domain}. Then, for every $0<\tau<1$,
\begin{equation}\label{eq:simple-column-square}
 \sum_{k=0}^{\infty}\tau^{2k}
 \sum_{\alpha\in\mathcal W_k}
 f_\alpha(B)^*f_\alpha(B)
 \preceq\frac1{1-\tau}I.
\end{equation}
\end{proposition}

\begin{proof}
Choose a diagonalization
\[
        B=S\diag(\beta_1,\ldots,\beta_d)S^{-1}
\]
and put $G=S^*S$. For $1\leq i\leq d$, set
\[
        \lambda_i=
        \begin{pmatrix}f_1(\beta_i)\\ \vdots\\ f_q(\beta_i)\end{pmatrix}
        \in\overline{\B_q}.
\]
For $1\leq j\leq q$, let
\[
        \Lambda_j=\diag(f_j(\beta_1),\ldots,f_j(\beta_d)).
\]
Thus $T_j=S\Lambda_jS^{-1}$. Put
\[
 D(w)=\sum_{j=1}^q w_j\Lambda_j
 =\diag(w^{\mathsf T}\lambda_1,\ldots,w^{\mathsf T}\lambda_d).
\]
Since the eigenvalues $\beta_i$ are distinct, polynomial interpolation
gives
\[
 \alg(B^*)=
 \{(S^{-1})^*\Delta S^*: \Delta\text{ is diagonal}\}.
\]
It follows from \eqref{eq:completion} that there is a diagonal
holomorphic matrix $\Psi(w)$, with $\Psi(0)=0$, such that
\begin{equation}\label{eq:balanced-H}
        S^*H(w)S=G(I-D(w))^{-1}+\Psi(w)G.
\end{equation}

Fix $0<\tau<1$. Set
\[
        c_{ij}=\lambda_i^*\lambda_j,
        \qquad w_i=\tau\overline{\lambda_i}.
\]
Some of the points $w_i$ may coincide. Repeated sampling points are
allowed in the positivity condition. We have
\[
        w_j^*w_i=\tau^2c_{ij},
        \qquad w_i^{\mathsf T}\lambda_j=\tau c_{ij}.
\]
Define Hermitian matrices $P,Q,Y$ by
\begin{equation}\label{eq:PQY}
        P_{ij}=\frac{G_{ij}}{1-\tau^2c_{ij}},
        \qquad
        Q_{ij}=\frac{G_{ij}}{1-\tau c_{ij}},
        \qquad Y=Q-P.
\end{equation}

Let $e_1,\ldots,e_d$ be the standard basis of $\C^d$. Fix
$u=(u_1,\ldots,u_d)^{\mathsf T}\in\C^d$. In the positive-kernel
inequality for \eqref{eq:H-kernel}, first take the congruence by $S$ in
each block. Use the points $w_i$ and the vectors $u_ie_i$. Add the
point $0$ with the vector
\begin{equation}\label{eq:origin-vector}
        v=-G^{-1}Pu.
\end{equation}
Write $\psi_i$ for the $i$th diagonal entry of $\Psi$. The contribution
of $\Psi$ to the resulting quadratic form is
\[
 2\operatorname{Re}\sum_{i=1}^d
 \overline{u_i}\,\psi_i(w_i)(Gv+Pu)_i,
\]
which is zero by \eqref{eq:origin-vector}.

The contribution of $G(I-D(w))^{-1}$ from the nonzero sample points is
the matrix $M$ with entries
\[
        M_{ij}=
        \frac{2G_{ij}}
        {(1-\tau c_{ij})(1-\tau^2c_{ij})}.
\]
The scalar identity
\[
 \frac{2}{(1-\tau x)(1-\tau^2x)}
 =\frac{2}{1-\tau^2x}
  +\frac{2}{1-\tau}
   \left(\frac{1}{1-\tau x}-\frac{1}{1-\tau^2x}\right)
\]
gives
\begin{equation}\label{eq:M-PY}
        M=2P+\frac{2}{1-\tau}Y.
\end{equation}
The block between the nonzero sample points and the origin is $G+Q$,
and the origin block is $2G$. Positivity and
\eqref{eq:origin-vector} give
\[
 0\preceq
 M-(G+Q)G^{-1}P-PG^{-1}(G+Q)+2PG^{-1}P.
\]
Using \eqref{eq:M-PY} and $Q=P+Y$, we obtain
\begin{equation}\label{eq:key-PY}
        \frac{2}{1-\tau}Y-YG^{-1}P-PG^{-1}Y\succeq0.
\end{equation}

Define
\[
 \widetilde T_j=G^{1/2}\Lambda_jG^{-1/2},\qquad
 \widehat P=G^{-1/2}PG^{-1/2},\qquad
 \widehat Q=G^{-1/2}QG^{-1/2},\qquad
 \widehat Y=\widehat Q-\widehat P.
\]
For $\alpha=(\alpha_1,\ldots,\alpha_k)$, put
\[
 \Lambda_\alpha=\Lambda_{\alpha_1}\cdots\Lambda_{\alpha_k},
 \qquad
 \widetilde T_\alpha=
 \widetilde T_{\alpha_1}\cdots\widetilde T_{\alpha_k},
\]
with $\Lambda_\varnothing=\widetilde T_\varnothing=I$. For every
$k\geq0$,
\[
 \sum_{\alpha\in\mathcal W_k}
 (\Lambda_\alpha^*G\Lambda_\alpha)_{ij}
 =G_{ij}c_{ij}^k.
\]
Expanding the denominators in \eqref{eq:PQY} gives norm-convergent
series
\begin{align}
 \widehat P
 &=\sum_{k=0}^{\infty}\tau^{2k}
   \sum_{\alpha\in\mathcal W_k}
   \widetilde T_\alpha^*\widetilde T_\alpha,
   \label{eq:P-series}\\
 \widehat Q
 &=\sum_{k=0}^{\infty}\tau^k
   \sum_{\alpha\in\mathcal W_k}
   \widetilde T_\alpha^*\widetilde T_\alpha.
   \label{eq:Q-series}
\end{align}
Since
\[
        \tau^k-\tau^{2k}
        \geq\frac{1-\tau}{\tau}\tau^{2k}
        \qquad(k\geq1),
\]
we have $\widehat P\succeq I$ and
\begin{equation}\label{eq:Y-lower}
        \widehat Y\succeq
        \frac{1-\tau}{\tau}(\widehat P-I).
\end{equation}
Congruence of \eqref{eq:key-PY} by $G^{-1/2}$ gives
\begin{equation}\label{eq:key-balanced}
        \frac{2}{1-\tau}\widehat Y
        -\widehat Y\widehat P
        -\widehat P\widehat Y\succeq0.
\end{equation}
Let $e$ be a unit eigenvector of $\widehat P$, with
$\widehat Pe=\gamma e$, and put $y=e^*\widehat Ye$. Then $\gamma\geq1$
and
\[
        0\leq2\left(\frac{1}{1-\tau}-\gamma\right)y.
\]
If $\gamma>1/(1-\tau)$, then \eqref{eq:Y-lower} gives $y>0$, which is
impossible. Hence
\begin{equation}\label{eq:P-upper}
        \widehat P\preceq\frac{1}{1-\tau}I.
\end{equation}

The polar decomposition $S=UG^{1/2}$ gives
$T_j=U\widetilde T_jU^*$. Thus \eqref{eq:P-series} and
\eqref{eq:P-upper} are unitarily equivalent to
\eqref{eq:simple-column-square}.
\end{proof}

We use two elementary approximation facts.

\begin{lemma}\label{lem:outer-approximation}
Let $K\subset\C$ be compact and convex, and let $V$ be an open set
containing $K$. There is a bounded convex domain $\Omega$ with smooth
boundary such that
\[
        K\subset\Omega,
        \qquad \overline\Omega\subset V.
\]
Moreover, matrices with simple spectrum are dense in $M_d(\C)$, and
\begin{equation}\label{eq:numerical-range-continuity}
        d_{\mathrm H}(W(A),W(B))\leq\norm{A-B}.
\end{equation}
\end{lemma}

\begin{proof}
Choose $\delta>0$ such that
$K+\delta\overline\D\subset V$. Let $d_K$ be the distance to $K$.
This is a convex $1$-Lipschitz function. Let $v_\eta$ be a smooth
convolution of $d_K$ at scale $\eta$, and put
\[
        u_\eta(z)=v_\eta(z)+\eta\abs{z}^2.
\]
The Lipschitz estimate gives $\abs{v_\eta-d_K}\leq\eta$. Choose $R$ with $K\subset R\overline\D$, and take $\eta>0$ so small
that
\[
        \eta(1+R^2)<\delta/4,
        \qquad \delta/2-\eta>\delta/4.
\]
Then the set
\[
        \Omega=\{z:u_\eta(z)<\delta/4\}
\]
contains $K$. If $d_K(z)\geq\delta/2$, then
$v_\eta(z)\geq d_K(z)-\eta>\delta/4$, so
$\Omega\subset K+(\delta/2)\D$. The function
$u_\eta$ is smooth, proper, and strictly convex. It has at most one
critical point, which is its global minimum. The minimum is less than
$\delta/4$, so the gradient does not vanish on the displayed level
set. Thus $\Omega$ has smooth boundary.

Simple-spectrum matrices are dense because the discriminant of the
characteristic polynomial is a nonzero polynomial in the matrix
entries. Finally, for every unit vector $x$,
\[
        \abs{x^*(A-B)x}\leq\norm{A-B}.
\]
Applying this with $A$ and $B$ interchanged proves
\eqref{eq:numerical-range-continuity}.
\end{proof}

\begin{proof}[Proof of Theorem~\ref{thm:tuple-square} for matrices]
First suppose that every $f_j$ is holomorphic on an open neighborhood
of $W(A)$. Fix $\varepsilon>0$. Lemma~\ref{lem:outer-approximation}
gives a bounded smooth convex domain $\Omega$ such that
$W(A)\subset\Omega$, every $f_j$ is holomorphic near
$\overline\Omega$, and
\[
        \sup_{z\in\overline\Omega}
        \sum_{j=1}^q\abs{f_j(z)}^2\leq1+\varepsilon.
\]
Put $\widetilde f_j=(1+\varepsilon)^{-1/2}f_j$. Choose
simple-spectrum matrices $B_n\to A$. By
\eqref{eq:numerical-range-continuity}, $W(B_n)\subset\Omega$ for all
large $n$. Proposition~\ref{prop:simple-square} applies to $B_n$ and
the functions $\widetilde f_j$.

Fix $N$. Retain the terms with $k\leq N$ and let $n\to\infty$.
Continuity of the holomorphic functional calculus gives
\[
 \sum_{k=0}^{N}\tau^{2k}
 \sum_{\alpha\in\mathcal W_k}
 \widetilde f_\alpha(A)^*\widetilde f_\alpha(A)
 \preceq\frac{1}{1-\tau}I.
\]
Letting $\varepsilon\to0^+$ gives the same inequality with $f_\alpha$
in place of $\widetilde f_\alpha$. The partial sums are increasing and
bounded. In finite dimensions their positive increments have summable
traces, and therefore summable norms. This proves norm convergence and
the column estimate.

For general $f_j\in\Acal(W(A))$, choose polynomial tuples converging
uniformly to $(f_1,\ldots,f_q)$. After multiplication by a scalar
converging to one, the polynomial tuples satisfy the pointwise column
bound. Applying the preceding result to each finite partial sum and
passing to the limit proves the column estimate. The row estimate
follows by applying the column estimate to $A^*$ and to
\[
        f_j^\#(z)=\overline{f_j(\overline z)}.
\]
The functions of $A$ commute, so reversal of a word does not change
its product.
\end{proof}

We now pass to bounded operators. The following argument is the
finite-dimensional Krylov compression used in numerical-range
functional-calculus estimates \cite{Crouzeix2007}.

\begin{proof}[Completion of the proof of Theorem~\ref{thm:tuple-square}]
Suppose first that $f_1,\ldots,f_q$ are polynomials, each of degree at
most $m$. Fix a unit vector $x\in H$ and an integer $N\geq0$. Let
\[
 \mathcal K_{N,x}
 =\operatorname{span}\{x,Ax,\ldots,A^{mN}x\},
 \qquad
 B=P_{\mathcal K_{N,x}}A|_{\mathcal K_{N,x}}.
\]
The numerical range of the compression $B$ is contained in $K_A$.
For every word $\alpha$ with $|\alpha|\leq N$,
\begin{equation}\label{eq:krylov-word}
        f_\alpha(B)x=f_\alpha(A)x.
\end{equation}
Indeed, the degree of $f_\alpha$ is at most $mN$, and
$B^jx=A^jx$ for $0\leq j\leq mN$. Applying the matrix case to $B$ and
using \eqref{eq:krylov-word} gives
\[
 \sum_{k=0}^{N}\tau^{2k}
 \sum_{\alpha\in\mathcal W_k}
 \norm{f_\alpha(A)x}^2
 \leq\frac1{1-\tau}.
\]
Since this holds for every unit vector $x$, the corresponding operator
partial sum is bounded by $(1-\tau)^{-1}I$.

Now let $f_1,\ldots,f_q\in\Acal(K_A)$ satisfy
\eqref{eq:intro-pointwise-column}. For each $n$, choose polynomials
$p_{j,n}$ converging uniformly to $f_j$ on $K_A$. There are numbers
$\varepsilon_n\to0^+$ such that
\[
 \sum_{j=1}^q\abs{p_{j,n}(z)}^2
 \leq(1+\varepsilon_n)^2
 \qquad(z\in K_A).
\]
Apply the polynomial result to $p_{j,n}/(1+\varepsilon_n)$. For each
fixed $N$, uniform convergence and continuity of the functional
calculus permit passage to the limit. Thus every column partial sum is
bounded by $(1-\tau)^{-1}I$. 
The increasing sequence of partial sums therefore converges
strongly to an operator satisfying the same bound.

In fact, the convergence is in operator norm. Fix
$\tau<\sigma<1$ and put
\[
    S_k=\sum_{\alpha\in\mathcal W_k}
    f_\alpha(A)^*f_\alpha(A).
\]
For $M>N$, the bound already proved at $\sigma$ gives
\[
\begin{aligned}
    0
    &\preceq \sum_{k=N+1}^{M}\tau^{2k}S_k\\
    &\preceq
    \left(\frac{\tau}{\sigma}\right)^{2(N+1)}
    \sum_{k=N+1}^{M}\sigma^{2k}S_k\\
    &\preceq
    \frac{(\tau/\sigma)^{2(N+1)}}{1-\sigma}I.
\end{aligned}
\]
Thus the partial sums are Cauchy in operator norm.

Apply the column estimate to $A^*$ and the reflected functions to
obtain the row estimate. Keeping only the terms of degrees zero and one
gives
\[
 \sum_{j=1}^q f_j(A)^*f_j(A)
 \preceq\frac1{\tau(1-\tau)}I,
\]
and the corresponding row inequality. The minimum of the scalar
factor occurs at $\tau=1/2$, which proves
\eqref{eq:intro-column-row-two}. For sharpness, take $q=1$,
$f_1(z)=z$, and
\[
        A=\begin{pmatrix}0&2\\0&0\end{pmatrix}.
\]
Then $W(A)=\overline\D$ and $\norm{A}=2$.
\end{proof}

\begin{corollary}\label{cor:small-tau-square}
Under the hypotheses of Theorem~\ref{thm:tuple-square}, let $P_\tau$
and $Q_\tau$ denote its column and row square functions. For
$0<\tau\leq1/2$,
\[
       I\preceq P_\tau,Q_\tau\preceq(1+4\tau^2)I.
\]
The upper bound is sharp for every $\tau$ in this interval.
\end{corollary}
\begin{proof}
For $k\geq1$ and $0<\tau\leq1/2$, one has
$\tau^{2k}\leq4\tau^2 4^{-k}$. Positivity of each term gives
\[
 P_\tau-I\preceq4\tau^2(P_{1/2}-I)\preceq4\tau^2I,
\]
and the same argument applies to $Q_\tau$. For sharpness, take
$A=\left(\begin{smallmatrix}0&0\\2&0\end{smallmatrix}\right)$ and
$f(z)=z$. Then $W(A)=\overline\D$, $A^2=0$, and both square
functions have norm $1+4\tau^2$.
\end{proof}
\par

\section{Stein equations and normal similarities}\label{sec:stein}

Let $T=(T_1,\ldots,T_q)$ be a commuting tuple of bounded operators. For
$X\in\Bcal(H)$, put
\[
 \mathcal L_T(X)=\sum_{j=1}^qT_j^*XT_j,
 \qquad
 \mathcal R_T(X)=\sum_{j=1}^qT_jXT_j^*.
\]
The two square functions in Theorem~\ref{thm:tuple-square} are the
resolvent series
\[
        P_\tau=\sum_{k=0}^{\infty}\tau^{2k}\mathcal L_T^k(I),
        \qquad
        Q_\tau=\sum_{k=0}^{\infty}\tau^{2k}\mathcal R_T^k(I).
\]

\begin{proposition}\label{prop:tuple-stein}
Under the assumptions of Theorem~\ref{thm:tuple-square}, with
$T_j=f_j(A)$, one has
\begin{align}
 P_\tau&=I+\tau^2\sum_{j=1}^qT_j^*P_\tau T_j,
 \label{eq:column-stein}\\
 Q_\tau&=I+\tau^2\sum_{j=1}^qT_jQ_\tau T_j^*.
 \label{eq:row-stein}
\end{align}
Moreover,
\[
        I\preceq P_\tau,Q_\tau\preceq\frac1{1-\tau}I.
\]
\end{proposition}

\begin{proof}
Let $P_{\tau,N}$ and $Q_{\tau,N}$ be the partial sums. They converge
strongly. Multiplication on the left and right by a fixed bounded
operator preserves strong convergence on bounded sets. Shifting the
two series gives \eqref{eq:column-stein} and
\eqref{eq:row-stein}. The bounds are contained in
Theorem~\ref{thm:tuple-square}.
\end{proof}

The two Stein equations give, in general, two different similarities.
One transforms the tuple into a column contraction and the other into a
row contraction. Both similarities depend on the tuple.

\begin{corollary}\label{cor:tuple-similarity}
Under the assumptions of Proposition~\ref{prop:tuple-stein}, put
\[
 C_j=\sqrt\tau\,P_\tau^{1/2}T_jP_\tau^{-1/2},
 \qquad
 D_j=\sqrt\tau\,Q_\tau^{-1/2}T_jQ_\tau^{1/2}.
\]
Then
\begin{equation}\label{eq:column-row-contractions}
        \sum_{j=1}^qC_j^*C_j\preceq I,
        \qquad
        \sum_{j=1}^qD_jD_j^*\preceq I.
\end{equation}
Furthermore,
\[
 \cond(P_\tau^{1/2}),\ \cond(Q_\tau^{1/2})
 \leq(1-\tau)^{-1/2}.
\]
For $0<\tau\leq1/2$, the condition numbers are at most
$\sqrt{1+4\tau^2}$, and the same similarities satisfy
\[
 \sum_j C_j^*C_j\preceq\frac{4\tau}{1+4\tau^2}I,
 \qquad
 \sum_j D_jD_j^*\preceq\frac{4\tau}{1+4\tau^2}I.
\]
\par
\end{corollary}

\begin{proof}
Equation \eqref{eq:column-stein} gives
\[
 \sum_{j=1}^qC_j^*C_j
 =\tau P_\tau^{-1/2}
   \left(\sum_{j=1}^qT_j^*P_\tau T_j\right)
   P_\tau^{-1/2}
 =\frac1\tau(I-P_\tau^{-1})\preceq I.
\]
The last inequality follows from
$P_\tau^{-1}\succeq(1-\tau)I$. The row calculation is
identical. The condition-number estimates follow from
$I\preceq P_\tau,Q_\tau\preceq(1-\tau)^{-1}I$.
For the additional estimates, use Corollary~\ref{cor:small-tau-square}
and the identities
$\sum_j C_j^*C_j=\tau^{-1}(I-P_\tau^{-1})$ and
$\sum_j D_jD_j^*=\tau^{-1}(I-Q_\tau^{-1})$.
\par
\end{proof}

We next prove a quantitative estimate for tuples similar to commuting
normal tuples. The following moment inequality is the required form of
H\"older's inequality.

\begin{lemma}\label{lem:joint-normal-moments}
Let $N_1,\ldots,N_q$ be commuting normal operators, and assume
\[
        \sum_{j=1}^qN_j^*N_j\preceq I.
\]
For a word $\alpha$, put
$N_\alpha=N_{\alpha_1}\cdots N_{\alpha_k}$. Then, for every vector $z$
and every $k\geq1$,
\begin{equation}\label{eq:joint-moment}
 \left(\sum_{|\alpha|=k}\norm{N_\alpha z}^2\right)
 \norm{z}^{2(k-1)}
 \geq
 \left(\sum_{j=1}^q\norm{N_jz}^2\right)^k.
\end{equation}
\end{lemma}

\begin{proof}
By the joint spectral theorem, there is a projection-valued measure
$E$ and bounded scalar functions $n_j$ such that
$N_j=\int n_j\dd E$. Put
\[
        a=\sum_{j=1}^q|n_j|^2,
        \qquad
        \omega(\Gamma)=\ip{E(\Gamma)z}{z}.
\]
The hypothesis gives $0\leq a\leq1$, and
\[
 \sum_{|\alpha|=k}\norm{N_\alpha z}^2
 =\int a^k\dd\omega,
 \qquad
 \sum_{j=1}^q\norm{N_jz}^2
 =\int a\dd\omega.
\]
H\"older's inequality gives \eqref{eq:joint-moment}.
\end{proof}

For $0<s<2$ and $0<\tau<1$, put
\begin{equation}\label{eq:L-tau-definition}
 L_\tau(s)
 =\frac{\tau^3(1-\tau)s^4}
 {1-\tau(1-\tau)s^2},
 \qquad
 L(s)=\sup_{0<\tau<1}L_\tau(s).
\end{equation}
The denominator is positive because $s<2$.

\begin{theorem}\label{thm:quantitative-tuple-normal}
Let $T_1,\ldots,T_q\in\Bcal(H)$ satisfy the column square-function
estimate in Theorem~\ref{thm:tuple-square}. Suppose that
\[
        T_j=SN_jS^{-1},
        \qquad 1\leq j\leq q,
\]
where $S$ is invertible, the operators $N_j$ are commuting and normal,
and
$\sum_jN_j^*N_j\preceq I$. Put
\[
        s=\norm{(T_1,\ldots,T_q)^{\mathsf T}},
        \qquad \kappa=\cond(S).
\]
Then, for every $0<\tau<1$,
\begin{equation}\label{eq:tuple-all-tau}
 \tau^2s^2+
 \frac{\tau^4s^4}
 {\kappa^4(1-\tau^2s^2/\kappa^2)}
 \leq\frac{\tau}{1-\tau}.
\end{equation}
Consequently,
\begin{equation}\label{eq:tuple-Phi-bound}
        s\leq\min\{\kappa,\Phi(\kappa)\}<2,
\end{equation}
where $\Phi$ is given by \eqref{eq:intro-Phi}. If $0<s<2$, then
\begin{equation}\label{eq:tuple-condition-lower}
        \kappa^4\geq L(s)
        >\frac{s^4}{4(4-s^2)}.
\end{equation}
The same assertions hold for a row tuple if the row square-function
estimate is assumed.
\end{theorem}

\begin{proof}
The column operator associated with $(N_1,\ldots,N_q)$ is a
contraction, so $s\leq\kappa$. Choose unit vectors $x_m$ such that
\[
        \sum_{j=1}^q\norm{T_jx_m}^2\longrightarrow s^2,
\]
and put $z_m=S^{-1}x_m$. For $k\geq1$, Lemma~\ref{lem:joint-normal-moments}
gives
\begin{align*}
 \sum_{|\alpha|=k}\norm{T_\alpha x_m}^2
 &\geq\frac1{\norm{S^{-1}}^2}
       \sum_{|\alpha|=k}\norm{N_\alpha z_m}^2\\
 &\geq\frac1{\norm{S^{-1}}^2}
 \frac{\left(\sum_j\norm{N_jz_m}^2\right)^k}
      {\norm{z_m}^{2(k-1)}}\\
 &\geq\frac1{\kappa^{2k}}
       \left(\sum_j\norm{T_jx_m}^2\right)^k.
\end{align*}
In the column square function, keep the exact term of degree one and
use the last lower bound from degree two onward. First let $m\to\infty$
in each finite partial sum and then let the number of terms tend to
infinity. We obtain
\[
 1+\tau^2s^2+
 \sum_{k=2}^{\infty}\frac{\tau^{2k}s^{2k}}{\kappa^{2k}}
 \leq\frac1{1-\tau}.
\]
The geometric series converges because $s\leq\kappa$ and $\tau<1$.
This is \eqref{eq:tuple-all-tau}.

If $\kappa=1$, then $s\leq1$. Suppose that $\kappa>1$. At
$\tau=1/2$, put $u=s^2/4$. Inequality
\eqref{eq:tuple-all-tau} becomes
\[
        u+\frac{u^2}{\kappa^4(1-u/\kappa^2)}\leq1.
\]
Equivalently,
\[
 (\kappa^2-1)u^2
 -\kappa^2(\kappa^2+1)u
 +\kappa^4\geq0.
\]
The smaller zero is
\[
        u_\kappa
        =\frac{2\kappa^2}
        {\kappa^2+1+\sqrt{(\kappa^2-1)^2+4}}.
\]
The first two terms of the square function give $0\leq u\leq1$. The
quadratic is positive at $0$ and negative at $1$, so $u\leq u_\kappa$.
This proves \eqref{eq:tuple-Phi-bound}.

If $0<s<2$, retaining only the terms of degrees zero, one, and two in
the square function gives
\[
        1+\tau^2s^2+\frac{\tau^4s^4}{\kappa^4}
        \leq\frac1{1-\tau}.
\]
Rearrangement gives $\kappa^4\geq L_\tau(s)$ for every
$0<\tau<1$, and hence $\kappa^4\geq L(s)$. At $\tau=1/2$,
\[
        L_{1/2}(s)=\frac{s^4}{4(4-s^2)}.
\]
A direct differentiation gives
\[
 \left.\frac{\partial}{\partial\tau}L_\tau(s)
 \right|_{\tau=1/2}
 =\frac{s^4}{4-s^2}>0,
\]
which proves the strict inequality in
\eqref{eq:tuple-condition-lower}.

For a row tuple, apply the column assertion to
$(T_1^*,\ldots,T_q^*)$. The similarity operator is
$(S^{-1})^*$ and has the same condition number. The operators $T_j$
commute, so the word sums agree after reversal.
\end{proof}

In the next corollary, the conclusion that $\psi(A)<2$
for every diagonalizable matrix $A$ also follows from
Theorem~\ref{thm:equality-block};
see Remark~\ref{rem:similar-normal-proof2}.

\begin{corollary}\label{cor:similar-normal}
Let $A\in\Bcal(H)$ and suppose that
\[
        A=VNV^{-1},
\]
where $V$ is invertible and $N$ is normal. Put
$\kappa=\cond(V)$. For every nonzero $f\in\Acal(K_A)$, set
\[
        s_f=\frac{\norm{f(A)}}{\norm{f}_{K_A}}.
\]
Then
\begin{equation}\label{eq:similar-normal-scalar}
        s_f\leq\min\{\kappa,\Phi(\kappa)\}<2.
\end{equation}
If $0<s_f<2$, then
\begin{equation}\label{eq:similar-normal-condition}
        \kappa^4\geq L(s_f)
        >\frac{s_f^4}{4(4-s_f^2)}.
\end{equation}
Consequently, every diagonalizable matrix, of arbitrary order, has
scalar Crouzeix ratio strictly smaller than $2$.
\end{corollary}

\begin{proof}
Normalize $f$ so that $\norm{f}_{K_A}=1$. Then
\[
        f(A)=Vf(N)V^{-1}.
\]
The operator $f(N)$ is normal, and
$\norm{f(N)}\leq1$ because $\sigma(N)=\sigma(A)\subset K_A$.
Apply Theorem~\ref{thm:quantitative-tuple-normal} with $q=1$.
\end{proof}

\section{Rectangular completely bounded extremals}\label{sec:rectangular}

Let $\mathfrak A$ be a finite-dimensional operator algebra with unity and
let
\[
        \theta:\mathfrak A\longrightarrow M_d(\C)
\]
be a homomorphism preserving the identity. For
$F\in M_{s,r}(\mathfrak A)$, the amplification
$\theta^{(s,r)}(F)$ acts from $\C^r\otimes\C^d$ to
$\C^s\otimes\C^d$.

For $x\in\C^r\otimes\C^d$, let $E_x$ be the smallest subspace of
$\C^r$ such that $x\in E_x\otimes\C^d$. Its dimension is the Schmidt
rank of $x$. The compression below is the rectangular form of the
finite-level argument associated with Smith's lemma \cite{Smith1983}.

\begin{lemma}\label{lem:abstract-levels}
One has
\begin{equation}\label{eq:abstract-finite-levels}
 \norm{\theta}_{\mathrm{cb}}
 =\max_{1\leq r,s\leq d}\norm{\theta^{(s,r)}}.
\end{equation}
At every fixed rectangular level, the norm is attained on the closed
unit ball of $M_{s,r}(\mathfrak A)$.
\end{lemma}

\begin{proof}
Let $F\in M_n(\mathfrak A)$ with $\norm{F}\leq1$, and choose unit
vectors $x,y\in\C^n\otimes\C^d$ such that
\[
        \abs{\ip{\theta^{(n)}(F)x}{y}}
        =\norm{\theta^{(n)}(F)}.
\]
Let $E_x,E_y\subset\C^n$ be the coefficient subspaces of $x$ and $y$,
and put $r=\dim E_x$, $s=\dim E_y$. Then $r,s\leq d$. The compression
\[
        P_{E_y}F|_{E_x}:E_x\longrightarrow E_y
\]
does not increase the matrix norm and preserves the displayed matrix
coefficient. This proves that the left side of
\eqref{eq:abstract-finite-levels} is bounded by the right side. The
reverse inequality follows by adjoining zero rows and columns.

The space $M_{s,r}(\mathfrak A)$ is finite-dimensional. Its closed
unit ball is compact, and the map
$F\mapsto\norm{\theta^{(s,r)}(F)}$ is continuous. Hence the norm is
attained.
\end{proof}

The next lemma uses the variational argument of Hartz and McCarthy
\cite[Lemma~4.1]{HartzMcCarthy2026}. We allow perturbations in a
$C^*$-subalgebra on which the homomorphism preserves adjoints.
The scalar-matrix case gives the rectangular reduction below, and the
more general form will be used for digraph algebras.

\begin{lemma}\label{lem:potapov}
Let $\mathcal A$ be an operator algebra with unity, and let
$\theta:\mathcal A\to M_d(\C)$ be a homomorphism preserving the identity.
Suppose that $F\in M_n(\mathcal A)$ satisfies $\norm{F}\leq1$ and
\[
 C=\norm{\theta^{(n)}}>1,\qquad
 R=\theta^{(n)}(F),\qquad \norm{Rx}=C
\]
for some unit vector $x\in\C^n\otimes\C^d$. Then
\begin{equation}\label{eq:potapov-orthogonality}
 \ip{Rx}{(D\otimes I)x}=0\qquad(D\in M_n(\C)).
\end{equation}
More generally, if $\mathcal D\subset\mathcal A$ is a $C^*$-subalgebra
with the same identity and $\theta|_{\mathcal D}$ is a
$*$-homomorphism, then
\begin{equation}\label{v5:eq:diagonal-orthogonality}
 \ip{Rx}{\theta^{(n)}(B)x}=0\qquad(B\in M_n(\mathcal D)).
\end{equation}
\end{lemma}

\begin{proof}
We prove the more general assertion. By
\cite[Theorem~3.1]{BlecherRuanSinclair1990}, we may represent
$\mathcal A$ completely isometrically as a subalgebra of $B(\mathcal L)$
containing the identity. This representation preserves adjoints on
$\mathcal D$ by \cite[Proposition~1.2.8]{Arveson1969}.
For $\norm{E}<1$, the transformation
\[
 m_E(X)=(I-EE^*)^{-1/2}(X+E)
        (I+E^*X)^{-1}(I-E^*E)^{1/2}
\]
maps the closed operator ball into itself. Indeed, put
$\Delta_E=(I-E^*E)^{1/2}$. For $\norm{X}\leq1$, direct calculation gives
\[
 I-m_E(X)^*m_E(X)
 =\Delta_E(I+X^*E)^{-1}(I-X^*X)(I+E^*X)^{-1}\Delta_E\geq0.
\]
Fix $B\in M_n(\mathcal D)$ and put $K=\theta^{(n)}(B)$.
For all sufficiently small $t\in\C$, the two square-root factors in
$m_{tB}(F)$ belong to $M_n(\mathcal D)$, and the remaining inverse is
given by a norm-convergent Neumann series. Hence
$m_{tB}(F)\in M_n(\mathcal A)$ and $\norm{m_{tB}(F)}\leq1$.
Since $C=\norm{Rx}<\infty$, the map $\theta^{(n)}$ is bounded.
Applying it to the Neumann and binomial series, and using
$\theta^{(n)}(B^*)=K^*$, gives
\[
 \theta^{(n)}(m_{tB}(F))=m_{tK}(R),\qquad
 \norm{m_{tK}(R)}\leq C.
\]
The target inverse exists for all sufficiently small $t$, since
$|t|\norm{K}\norm{R}<1$.

Since $\norm{R}\leq C$ and $\norm{Rx}=C$, one has
$R^*Rx=C^2x$. Expansion at $t=0$ gives
\[
 m_{tK}(R)x=Rx+tKx-\overline t\,RK^*Rx+O(|t|^2).
\]
With inner products linear in the first variable, it follows that
\[
 \norm{m_{tK}(R)x}^2
 =C^2+2(1-C^2)\operatorname{Re}
       \bigl(\overline t\,\ip{Rx}{Kx}\bigr)+O(|t|^2).
\]
The left side is at most $C^2$. Varying the argument of $t$ and then
letting $t\to0$ proves \eqref{v5:eq:diagonal-orthogonality}.
Taking $\mathcal D=\operatorname{span}\{1\}$ proves \eqref{eq:potapov-orthogonality}.
\end{proof}
\par

\begin{theorem}\label{thm:abstract-rectangular-extremal}
Suppose that $\mathfrak A$ is finite-dimensional and
$\norm{\theta}_{\mathrm{cb}}>1$. There are positive integers $r,s$ and
an element $F\in M_{s,r}(\mathfrak A)$ such that
\begin{equation}\label{eq:abstract-extremal}
        r+s\leq d,
        \qquad
        \norm{F}\leq1,
        \qquad
        \norm{\theta^{(s,r)}(F)}
        =\norm{\theta}_{\mathrm{cb}}.
\end{equation}
If $r+s$ is minimal, norming input and output vectors have Schmidt
ranks $r$ and $s$, and their Schmidt-vector families in $\C^d$ are
mutually orthogonal.
\end{theorem}

\begin{proof}
Put $C=\norm{\theta}_{\mathrm{cb}}>1$. By
Lemma~\ref{lem:abstract-levels}, choose a pair $(r,s)$ with $r+s$
minimal among the pairs at which $C$ is attained, and choose an
extremal $F\in M_{s,r}(\mathfrak A)$. There are unit vectors
\[
        x\in\C^r\otimes\C^d,
        \qquad
        y\in\C^s\otimes\C^d
\]
such that
\begin{equation}\label{eq:abstract-singular-vectors}
        \theta^{(s,r)}(F)x=Cy,
        \qquad
        \theta^{(s,r)}(F)^*y=Cx.
\end{equation}
Minimality implies that the Schmidt ranks of $x$ and $y$ are $r$ and
$s$. Otherwise, the compression in Lemma~\ref{lem:abstract-levels}
would give an extremal with a smaller sum.

On $\C^s\oplus\C^r$, define
\[
        \widehat F=
        \begin{pmatrix}0&F\\0&0\end{pmatrix},
        \qquad
        \widehat x=\begin{pmatrix}0\\x\end{pmatrix},
        \qquad
        \widehat y=\begin{pmatrix}y\\0\end{pmatrix}.
\]
Then $\norm{\widehat F}\leq1$ and
$\theta^{(r+s)}(\widehat F)\widehat x=C\widehat y$. Hence
$\norm{\theta^{(r+s)}}=C$. Apply
Lemma~\ref{lem:potapov} to $\widehat F$ and $\widehat x$. Using
\eqref{eq:abstract-singular-vectors}, we obtain
\[
        \ip{\widehat y}{(D\otimes I)\widehat x}=0
        \qquad(D\in M_{r+s}(\C)).
\]
Taking
\[
        D=\begin{pmatrix}0&D_0\\0&0\end{pmatrix},
        \qquad D_0\in M_{s,r}(\C),
\]
gives
\begin{equation}\label{eq:rectangular-orthogonality}
        \ip{y}{(D_0\otimes I)x}=0
        \qquad(D_0\in M_{s,r}(\C)).
\end{equation}

Take Schmidt decompositions
\[
        x=\sum_{j=1}^r a_je_j\otimes\xi_j,
        \qquad
        y=\sum_{i=1}^s b_if_i\otimes\eta_i,
\]
where $a_j,b_i>0$, and where $(\xi_j)$ and $(\eta_i)$ are orthonormal
families in $\C^d$. Taking $D_0$ to be a matrix unit in
\eqref{eq:rectangular-orthogonality} gives
\[
        \ip{\eta_i}{\xi_j}=0
        \qquad(1\leq i\leq s,\ 1\leq j\leq r).
\]
Thus the two Schmidt-vector families are mutually orthogonal, and
$r+s\leq d$.
\end{proof}

The first variation also gives simplicity of the largest singular
value at a minimal rectangular extremal in the following cases.

\begin{proposition}\label{v3:prop:extremal-simple}
Let $F\in M_{s,r}(\mathfrak A)$ satisfy
\eqref{eq:abstract-extremal}, with $r+s$ minimal, and put
$C=\norm{\theta}_{\mathrm{cb}}>1$.
If $r+s\geq d-1$, then the largest singular value of
$R=\theta^{(s,r)}(F)$ is simple.
\end{proposition}

\begin{proof}
Let $\mathcal K=\ker(R^*R-C^2I)$. Every nonzero vector in
$\mathcal K$ norms $R$ after normalization. Applying
Lemma~\ref{lem:potapov} to the square matrix obtained by placing $F$
in an upper-right corner, as in the preceding proof, gives
\[
        \ip{Rx}{(D\otimes I)x}=0
        \qquad(x\in\mathcal K,\ D\in M_{s,r}(\C)).
\]
For each fixed $D$, polarization on the complex vector space
$\mathcal K$ therefore gives
\[
        \ip{Ru}{(D\otimes I)v}=0
        \qquad(u,v\in\mathcal K).
\]
Let $E,G\subset\C^d$ be the smallest subspaces for which
\[
        \mathcal K\subset\C^r\otimes E,
        \qquad
        R\mathcal K\subset\C^s\otimes G.
\]
Taking $D$ to be a matrix unit shows that $E\perp G$.
Minimality of $r+s$ implies that every nonzero $x\in\mathcal K$
and $Rx$ have Schmidt ranks $r$ and $s$, respectively. Hence
\[
        \dim E\geq r,
        \qquad
        \dim G\geq s.
\]

Since $E\perp G$,
\[
 0\leq(\dim E-r)+(\dim G-s)\leq d-r-s\leq1.
\]
Consequently, either $\dim E=r$ or $\dim G=s$.
In the first case, identify $\C^r\otimes E$ with $M_r(\C)$.
Every nonzero element of $\mathcal K$ is invertible. If $X,Y$ were
linearly independent elements of $\mathcal K$, then $\det(X+tY)$
would have degree $r$ with nonzero leading coefficient. A complex
root would give a nonzero singular element of $\mathcal K$.
Thus $\dim\mathcal K=1$. In the second case, apply the same argument
to $R\mathcal K\subset\C^s\otimes G\simeq M_s(\C)$.
Since $R$ is injective on $\mathcal K$, again $\dim\mathcal K=1$.
\par
\end{proof}

The corresponding scalar functional-calculus statement is the exact
similarity amplification theorem of Aouichaoui and O'Loughlin
\cite[Theorem~3.2]{AouichaouiOloughlin2026}.

\begin{theorem}\label{thm:cb-similarity-amplification}
Let $\mathfrak A$ be an operator algebra with unity and let
$\theta:\mathfrak A\to M_d(\C)$ be a bounded homomorphism preserving
the identity with
\[
        C=\norm{\theta}_{\mathrm{cb}}>1.
\]
For every $\gamma\geq1$, there is a positive invertible matrix
$S_\gamma\in M_d(\C)$ such that
\begin{equation}\label{eq:cb-similarity-amplification}
 \cond(S_\gamma)=\gamma,
 \qquad
 \norm{\operatorname{Ad}_{S_\gamma^{-1}}\circ\theta}_{\mathrm{cb}}
 =\gamma C,
\end{equation}
where
$\operatorname{Ad}_{S^{-1}}(X)=S^{-1}XS$.
\end{theorem}

\begin{proof}
Every bounded map into $M_d(\C)$ is completely bounded by the
coefficient-space compression argument in Lemma~\ref{lem:abstract-levels}.
The ideal $J=\ker\theta$ is closed, and $\mathfrak A/J$, with its
quotient matrix norms, is finite-dimensional, since it is linearly
isomorphic to $\theta(\mathfrak A)$. Let $\overline\theta$ be the
induced homomorphism. At each rectangular level,
\[
 \norm{\overline\theta^{(s,r)}}=\norm{\theta^{(s,r)}}.
\]
Indeed, the quotient map is completely contractive, and a quotient
matrix of norm at most one has, for every $\varepsilon>0$, a lift
of norm at most $1+\varepsilon$. Applying $\theta$ to such lifts
and letting $\varepsilon\to0$ gives the reverse inequality.
It is therefore enough to carry out the following construction for
$\overline\theta$ on $\mathfrak A/J$.

Choose a minimal rectangular extremal as in
Theorem~\ref{thm:abstract-rectangular-extremal}. Thus there are
$F\in M_{s,r}(\mathfrak A/J)$ and unit vectors
$x\in\C^r\otimes\C^d$, $y\in\C^s\otimes\C^d$ such that
\[
        \norm{F}\leq1,
        \qquad
        \overline\theta^{(s,r)}(F)x=Cy.
\]
Let $E$ and $G$ be the subspaces of $\C^d$ spanned by the Schmidt
vectors of $x$ and $y$, respectively. Theorem~\ref{thm:abstract-rectangular-extremal}
gives $E\perp G$. Define
\[
 S_\gamma
 =\sqrt\gamma\,P_E
 +\frac1{\sqrt\gamma}P_G
 +P_{(E\oplus G)^\perp}.
\]
Then $S_\gamma$ is positive and
$\cond(S_\gamma)=\gamma$. Moreover,
\[
 (I_r\otimes S_\gamma)x=\sqrt\gamma\,x,
 \qquad
 (I_s\otimes S_\gamma^{-1})y=\sqrt\gamma\,y.
\]
Therefore
\[
 \begin{aligned}
 &\left(\operatorname{Ad}_{S_\gamma^{-1}}\circ\overline\theta\right)^{(s,r)}
 (F)x\\
 &\quad=(I_s\otimes S_\gamma^{-1})
       \overline\theta^{(s,r)}(F)(I_r\otimes S_\gamma)x
 =\gamma Cy.
 \end{aligned}
\]
This proves the lower bound in
\eqref{eq:cb-similarity-amplification}. The reverse bound follows from
\[
 \norm{\operatorname{Ad}_{S_\gamma^{-1}}\circ\overline\theta}_{\mathrm{cb}}
 \leq\norm{S_\gamma^{-1}}\norm{S_\gamma}
       \norm{\overline\theta}_{\mathrm{cb}}
 =\gamma C.
\]

The same approximate-lift argument applies to
$\operatorname{Ad}_{S_\gamma^{-1}}\circ\theta$. It transfers the
computed norm $\gamma C$ from the quotient to $\mathfrak A$.
\par
\end{proof}

\section{\texorpdfstring{Complete bounds and equality in dimension three}{Complete bounds and equality in dimension three}}\label{sec:complete-three}

We first transfer the square functions from the uniform algebra to its
quotient.

\begin{lemma}\label{lem:quotient-square}
Let $A\in M_d(\C)$. Suppose that
$X=(x_1,\ldots,x_q)^{\mathsf T}\in M_{q,1}(\mathfrak A_A)$ satisfies
$\norm{X}\leq1$, and put $T_j=\theta_A(x_j)$. Then, for every
$0<\tau<1$,
\begin{equation}\label{eq:quotient-column-square}
 \sum_{k=0}^{\infty}\tau^{2k}
 \sum_{\alpha\in\mathcal W_k}T_\alpha^*T_\alpha
 \preceq\frac1{1-\tau}I.
\end{equation}
The analogous row estimate holds for every
$X\in M_{1,q}(\mathfrak A_A)$ with $\norm{X}\leq1$.
\end{lemma}

\begin{proof}
Fix $\varepsilon>0$. Lemma~\ref{lem:quotient-lifts} gives a lift
$F=(f_1,\ldots,f_q)^{\mathsf T}$ with
\[
        \max_{z\in W(A)}
        \left(\sum_{j=1}^q\abs{f_j(z)}^2\right)^{1/2}
        \leq1+\varepsilon.
\]
The operators $f_j(A)$ are the fixed operators $T_j$. Apply
Theorem~\ref{thm:tuple-square} to $f_j/(1+\varepsilon)$. For every
finite partial sum, let $\varepsilon\to0^+$. This gives
\eqref{eq:quotient-column-square}. The row assertion is proved in the
same way.
\end{proof}

\begin{proof}[Proof of Theorem~\ref{thm:intro-complete-three}]
Let $A\in M_3(\C)$ and put $C=\psi_{\mathrm{cb}}(A)$. By
Proposition~\ref{prop:quotient-calculus},
$C=\norm{\theta_A}_{\mathrm{cb}}$. If $C\leq1$, there is nothing to
prove. Suppose that $C>1$. Theorem~\ref{thm:abstract-rectangular-extremal}
gives positive integers $r,s$ and
$X\in M_{s,r}(\mathfrak A_A)$ such that
\[
        r+s\leq3,
        \qquad \norm{X}\leq1,
        \qquad \norm{\theta_A^{(s,r)}(X)}=C.
\]
The possible pairs are
\[
        (r,s)=(1,1),\qquad(1,2),\qquad(2,1).
\]
They are respectively a scalar, a two-component column, and a
two-component row. Lemma~\ref{lem:quotient-square}, at $\tau=1/2$ and
after retaining the terms of degrees zero and one, gives $C\leq2$ in
each case. This proves \eqref{eq:intro-complete-three}.

The same proof applies in dimension two. For the lower estimate, let
\[
        J=\begin{pmatrix}0&2\\0&0\end{pmatrix}.
\]
Then $W(J)=\overline\D$ and $\norm{J}=2$. Direct sums with a zero block
give the same lower estimate in dimension three. Hence
$Q_{\mathrm{cb}}(2)=Q_{\mathrm{cb}}(3)=2$.

More generally, if $\psi_{\mathrm{cb}}(A)>2$, a minimal extremal cannot
have $r=1$ or $s=1$, by Lemma~\ref{lem:quotient-square}. Thus
$r,s\geq2$ and $r+s\leq d$. When $d=4$, this leaves only
$(r,s)=(2,2)$.
\end{proof}

We now identify equality in dimension three.

\begin{proof}[Proof of Theorem~\ref{thm:intro-complete-equality}]
The implication from right to left follows from
\[
        \psi(A)\leq\psi_{\mathrm{cb}}(A)\leq2.
\]
Suppose that $\psi_{\mathrm{cb}}(A)=2$ and $\psi(A)<2$. Choose a
minimal rectangular extremal as in
Theorem~\ref{thm:abstract-rectangular-extremal}. The scalar pair
$(1,1)$ cannot attain the completely bounded norm. Hence the pair is
either $(1,2)$ or $(2,1)$.

Suppose first that the pair is $(1,2)$. Write
\[
        X=\begin{pmatrix}x_1\\x_2\end{pmatrix},
        \qquad T_j=\theta_A(x_j).
\]
There is a unit vector $\xi\in\C^3$ such that
\begin{equation}\label{eq:column-equality-vector}
        \norm{T_1\xi}^2+\norm{T_2\xi}^2=4.
\end{equation}
At $\tau=1/2$, Lemma~\ref{lem:quotient-square} gives
\[
 I+\frac14\sum_{j=1}^2T_j^*T_j
 +\frac1{16}\sum_{i,j=1}^2(T_iT_j)^*(T_iT_j)
 +\cdots\preceq2I.
\]
Taking the quadratic form at $\xi$ and using
\eqref{eq:column-equality-vector}, the first two displayed terms already
have value $2$. Hence
\begin{equation}\label{eq:second-words-vanish}
        T_iT_j\xi=0
        \qquad(1\leq i,j\leq2).
\end{equation}

Let
\[
        \eta=\frac12
        \begin{pmatrix}T_1\xi\\T_2\xi\end{pmatrix}
        \in\C^2\otimes\C^3.
\]
The minimality assertion in
Theorem~\ref{thm:abstract-rectangular-extremal} says that $\eta$ has
Schmidt rank two. Thus $T_1\xi$ and $T_2\xi$ are linearly independent.
The orthogonality relation \eqref{eq:rectangular-orthogonality} gives
\begin{equation}\label{eq:column-vector-orthogonality}
        \xi\perp T_1\xi,
        \qquad
        \xi\perp T_2\xi.
\end{equation}
Since $\dim\C^3=3$, the two vectors $T_1\xi,T_2\xi$ span
$\xi^\perp$. It follows from \eqref{eq:second-words-vanish} that each
$T_i$ vanishes on $\xi^\perp$. Therefore
\begin{equation}\label{eq:square-zero-space}
        T_iT_j=0
        \qquad(1\leq i,j\leq2).
\end{equation}

The space $L=\operatorname{span}\{T_1,T_2\}$ is two-dimensional and
$L^2=\{0\}$. It is contained in $\alg(A)$. Since $I\notin L$ and
$\dim\alg(A)\leq3$, one has
\[
        \alg(A)=\C I\oplus L.
\]
Thus $A=aI+N$ for some $N\in L$. By \eqref{eq:square-zero-space},
$N^2=0$, and consequently
\[
        \alg(A)=\operatorname{span}\{I,N\}.
\]
This contradicts $\dim(\C I\oplus L)=3$.

If the minimal pair is $(2,1)$, apply the preceding argument to $A^*$
and the adjoint column. Lemma~\ref{lem:adjoint-symmetry} preserves both
the scalar and completely bounded ratios. This completes the proof.
\end{proof}

We next prove the quantitative estimate for diagonalizable matrices.

\begin{proof}[Proof of Theorem~\ref{thm:intro-complete-diagonalizable}]
Put $\kappa=\cond(V)$ and $C=\psi_{\mathrm{cb}}(A)$. If $C=1$, the
conclusion is immediate. Otherwise, choose a minimal rectangular
extremal $X\in M_{s,r}(\mathfrak A_A)$. As before,
$(r,s)$ is one of $(1,1),(1,2),(2,1)$.

Consider first a column $X=(x_1,\ldots,x_q)^{\mathsf T}$, where
$q=1$ or $q=2$, and put $T_j=\theta_A(x_j)$. The class $x_j$ has
well-defined values on the eigenvalues of $A$. Indeed, if
$g\in\mathcal I_A$, then
\[
        0=V^{-1}g(A)V=g(\Lambda),
\]
so $g$ vanishes at every diagonal entry of $\Lambda$. Let
\[
        N_j=V^{-1}T_jV.
\]
The matrices $N_j$ are diagonal. Lemma~\ref{lem:quotient-lifts} gives,
for every $\varepsilon>0$, a column lift whose pointwise norm is at
most $1+\varepsilon$. Evaluation at the eigenvalues and passage to the
limit give
\[
        \sum_{j=1}^qN_j^*N_j\preceq I.
\]
Lemma~\ref{lem:quotient-square} gives the required column square
function for the tuple $(T_j)$. Theorem~\ref{thm:quantitative-tuple-normal}
therefore gives
\[
        \norm{(T_1,\ldots,T_q)^{\mathsf T}}
        \leq\min\{\kappa,\Phi(\kappa)\}.
\]
The row case follows by taking adjoints. Since the chosen extremal has
norm $C$, the displayed estimate proves
\eqref{eq:intro-complete-diagonalizable}. The formula
\eqref{eq:intro-Phi} gives $\Phi(\kappa)<2$ for every finite
$\kappa$.
\end{proof}

\section{The abstract Crouzeix conjecture}
\label{sec:abstract-cor}

We first consider the scalar estimate. Badea, O'Loughlin, and Virtanen
proved the scalar Clou\^atre--Ostermann--Ransford conjecture for bounded
homomorphisms on uniform algebras, without requiring the antilinear map
to be contractive \cite[Theorem~2.1]{BadeaOloughlinVirtanen2026}.
The following theorem removes the initial boundedness assumptions and
allows more general Banach algebras, under an explicit commutation
hypothesis.

A complex Banach algebra $\Bcal$ with unity satisfies von Neumann's
inequality if
\[
 \norm{p(b)}\leq\max_{|z|\leq1}|p(z)|
 \qquad(\norm{b}\leq1)
\]
for every polynomial $p$. Uniform algebras and operator algebras with unity
have this property. We use the disk algebra $A(\D)$ with the uniform
norm on $\overline\D$.

\begin{theorem}\label{v3:scalar-automatic}
Let $\Bcal$ be a complex Banach algebra with unity satisfying von Neumann's
inequality, let $\mathcal H$ be a nonzero Hilbert space, and let
$\theta:\Bcal\to B(\mathcal H)$ be a complex-linear homomorphism preserving the identity, not assumed to be bounded. Suppose that
$\Psi:\Bcal\to B(\mathcal H)$ is a complex-linear contraction preserving the identity
and
\begin{equation}\label{v3:eq-scalar-commutation}
 [\theta(b),2\Psi(b^n)^*-\theta(b^n)^*]=0
 \qquad(b\in\Bcal,\ n\geq1).
\end{equation}
Then $\theta$ is bounded and $\norm{\theta}\leq2$.

In particular, let $\alpha:\Bcal\to\Bcal$ be an antilinear map preserving the identity
and put
\[
 \Phi_{\theta,\alpha}(b)
 =\frac12\bigl(\theta(b)+\theta(\alpha(b))^*\bigr).
\]
If $\Phi_{\theta,\alpha}$ is contractive and
\[
 [\theta(b),\theta(\alpha(b^n))]=0
 \qquad(b\in\Bcal,\ n\geq1),
\]
then $\theta$ is bounded and $\norm{\theta}\leq2$.
The commutation condition holds whenever $\theta(\Bcal)$ is
commutative. No boundedness assumption on $\alpha$ is needed.
\end{theorem}

\begin{proof}
Fix $b\in\Bcal$ with $\norm{b}\leq1$ and put $T=\theta(b)$.
Since an algebra homomorphism preserving the identity preserves inverses,
\[
 \sigma(T)\subseteq\sigma_{\Bcal}(b),
 \qquad r(T)\leq1,
\]
where $r(T)$ is the spectral radius. This uses only the algebraic
properties of $\theta$. By von Neumann's inequality, polynomial
evaluation at $b$ extends to a contractive homomorphism $C_b:A(\D)\to\Bcal$ preserving the identity. Thus $L=\Psi\circ C_b$ is a contraction preserving the identity.

We recall why ordinary contractivity of $L$ gives a unitary dilation.
The formula
\[
 \widetilde L(f+\overline g)=L(f)+L(g)^*,
 \qquad f,g\in A(\D),
\]
is well defined on $A(\D)+A(\D)^*$, because the intersection of the
disk algebra and its conjugate consists of the constants. For a unit
vector $x\in\mathcal H$, the functional
$f\mapsto\ip{L(f)x}{x}$ has norm and value at $1$ equal to one.
Its norm-preserving Hahn--Banach extension to $C(\T)$ is a state.
Testing against these extensions shows that $\widetilde L$ is positive.
Positivity and the identity $\widetilde L(1)=I$ bound $\widetilde L$ on self-adjoint elements,
and hence on the whole space. Since the trigonometric polynomials are
dense in $C(\T)$, it extends to a positive map preserving the identity on $C(\T)$.
This map is completely positive because $C(\T)$ is commutative
\cite[Theorem~4]{Stinespring1955}. Stinespring's representation theorem,
with the normalization in \cite[Theorem~1 and Section~3]{Stinespring1955},
therefore gives a Hilbert space $\mathcal K$,
an isometry $V:\mathcal H\to\mathcal K$, and a unitary
$U\in B(\mathcal K)$ such that
\[
 \Psi(b^n)=V^*U^nV \qquad(n\geq1).
\]
Set
\[
 E_n=2V^*U^{*n}V-T^{*n}
     =2\Psi(b^n)^*-\theta(b^n)^*.
\]
These operators commute with $T$ by
\eqref{v3:eq-scalar-commutation}. For $0<t<1$, the operators $tT$ and
$tU$ have the corresponding perturbations $t^nE_n$. They are uniformly
bounded in $n$, since
\[
 \norm{t^nE_n}\leq2+\norm{(tT)^n}
 \quad\text{and}\quad r(tT)<1.
\]
Lemma~1 of Lorist and Schwenninger \cite{LoristSchwenninger2026}
therefore gives $t\norm{T}\leq2$. Letting $t\uparrow1$ proves
$\norm{\theta(b)}\leq2$, and taking the supremum over $b$ proves
boundedness and the asserted norm estimate.

For the final assertion, $\Phi_{\theta,\alpha}$ is complex-linear, satisfies $\Phi_{\theta,\alpha}(1)=I$, and
\[
 2\Phi_{\theta,\alpha}(b^n)^*-\theta(b^n)^*
 =\theta(\alpha(b^n)).
\]
Apply the first part with $\Psi=\Phi_{\theta,\alpha}$.
\end{proof}

The commutation hypothesis in Theorem~\ref{v3:scalar-automatic}
cannot be omitted in general, even if both $\alpha$ and
$\Phi_{\theta,\alpha}$ are completely contractive.
Theorem~\ref{v5:thm:cor-counterexample} gives such an example, and
\eqref{v5:eq:commutator} displays the failure of commutation.
\par

\begin{remark}\label{v3:no-prescribed-dilation}
The dilation in the proof of Theorem~\ref{v3:scalar-automatic}
represents the contractive map to which it is applied. Its compressed
coordinate need not agree with the homomorphism $\theta$.

For example, take $\Bcal=A(\D)$ and define
\[
 \theta(f)=f(1/4),
 \qquad
 \alpha(f)(\zeta)=\overline{f(-\overline\zeta)}.
\]
Then $\theta:\Bcal\to\C$ is a completely contractive homomorphism preserving the identity. The map $\alpha$ is antilinear and multiplicative, with $\alpha(1)=1$ and $\alpha^2=I$. For $F=[f_{ij}]\in M_n(A(\D))$,
\[
 [\alpha(f_{ij})](\zeta)=\overline{F(-\overline\zeta)},
 \qquad
 \norm{[\alpha(f_{ij})]}=\norm{F},
\]
where the bar denotes entrywise complex conjugation. Thus $\alpha$
is a complete isometry. The symmetrized map is
\[
 \Phi(f):=\Phi_{\theta,\alpha}(f)
 =\frac12\bigl(f(1/4)+f(-1/4)\bigr).
\]
It is completely contractive, since
\[
 \norm{\Phi^{(n)}(F)}
 \leq\frac12\bigl(\norm{F(1/4)}+\norm{F(-1/4)}\bigr)
 \leq\norm{F}_{M_n(A(\D))}.
\]
Let $z$ denote the coordinate function. Every unitary moment dilation
of $\Phi$, with isometry $V$ and unitary $U$, satisfies
\[
 V^*UV=\Phi(z)=0,
 \qquad
 \theta(z)=\frac14.
\]
Hence it cannot also satisfy $V^*UV=\theta(z)$. The example satisfies
the hypotheses of the complete abstract problem and
$\norm{\theta}_{\mathrm{cb}}=1$. It therefore excludes this additional
moment requirement, without contradicting the abstract Crouzeix bound.
\end{remark}

For an operator algebra with unity, a complex-linear map $P$ is called real
positive if $\Ree b\geq0$ implies $\Ree P(b)\geq0$, where
$\Ree b=(b+b^*)/2$. The next consequence extends
\cite[Theorem~3.3]{BadeaOloughlinVirtanen2026} from uniform algebras to
operator algebras with unity and removes the initial boundedness
assumptions on both maps.

\begin{corollary}\label{v3:real-positive}
Let $\Bcal$ be an operator algebra with unity, let $\mathcal H$ be a nonzero
Hilbert space, and let $\theta:\Bcal\to B(\mathcal H)$ be a
complex-linear homomorphism preserving the identity. Suppose that
$P:\Bcal\to B(\mathcal H)$ is complex-linear and real positive,
$P(1)=2I$, and
\[
 [P(b^n)-\theta(b^n),\theta(b)^*]=0
 \qquad(b\in\Bcal,\ n\geq1).
\]
Then $P$ and $\theta$ are bounded, $\norm{P}=2$, and
$\norm{\theta}\leq2$.
\end{corollary}

\begin{proof}
By \cite[Lemma~2.1]{BlecherNeal2021}, a real-positive map between
operator algebras with unity is bounded and has norm equal to the norm of
its value at the identity. Hence $\norm{P}=2$ and $\Psi=P/2$ is a
contraction preserving the identity. Taking adjoints in the commutation condition gives
\eqref{v3:eq-scalar-commutation}. Theorem~\ref{v3:scalar-automatic}
now applies.
\end{proof}

We next give the bounds for operator radii, using the exact
formula of Badea, Crouzeix, and Klaja
\cite[Theorem~3.2]{BadeaCrouzeixKlaja2018}. For $\rho\geq1$, let
$\mathcal C_\rho$ consist of the operators $S\in B(\mathcal H)$ for
which there are a Hilbert space $\mathcal K\supset\mathcal H$ and a
unitary $U\in B(\mathcal K)$ satisfying
\[
 S^n=\rho P_{\mathcal H}U^n|_{\mathcal H}\qquad(n\geq1).
\]
Define
\[
 w_\rho(S)=\inf\{t>0:t^{-1}S\in\mathcal C_\rho\},
 \qquad
 \norm{\theta}_{w_\rho}
 =\sup_{\norm{b}\leq1}w_\rho(\theta(b)).
\]
Thus $w_1(S)=\norm{S}$ and $w_2(S)=w(S)$, where
$w(S)=\sup_{\norm{x}=1}|\ip{Sx}{x}|$ is the numerical radius.
For $K\geq1$, put
\[
 \begin{aligned}
 \kappa_\rho(K)
 &=\frac{K^2+1+\sqrt{(K^2+1)^2-4\rho(2-\rho)K^2}}{2\rho K},\\
 c_\rho=\kappa_\rho(2)
 &=\frac{5+\sqrt{9+16(\rho-1)^2}}{4\rho}.
 \end{aligned}
\]

\begin{corollary}\label{v3:scalar-radii}
Under the hypotheses of Theorem~\ref{v3:scalar-automatic} or
Corollary~\ref{v3:real-positive}, let $K=\norm{\theta}$. Then
\[
 \norm{\theta}_{w_\rho}=\kappa_\rho(K)\leq c_\rho
 \qquad(\rho\geq1).
\]
Each constant $c_\rho$ is best possible. In particular,
$\sup_{\norm{b}\leq1}w(\theta(b))\leq5/4$, and $5/4$ is best possible.
\end{corollary}

\begin{proof}
The preceding results give $1\leq K\leq2$. The equality is
\cite[Theorem~3.2(a)]{BadeaCrouzeixKlaja2018}. To obtain the bound, write
\[
 x_K=\frac12(K+K^{-1}),
 \qquad
 \kappa_\rho(K)
 =\frac{x_K+\sqrt{x_K^2-\rho(2-\rho)}}{\rho}.
\]
Both expressions are increasing for $K\geq1$, so
$\kappa_\rho(K)\leq\kappa_\rho(2)$. Sharpness follows from the example
below.
\end{proof}

For a completely bounded map on an operator algebra, write
\[
 \norm{\theta}_{w_\rho,\mathrm{cb}}
 =\sup_{m\geq1}\sup_{\norm{F}_{M_m(\Bcal)}\leq1}
 w_\rho(\theta^{(m)}(F)).
\]
The following consequence uses only ordinary contractivity of the
symmetrized map, even when the domain algebra is noncommutative.

\begin{corollary}\label{v3:rank-two}
Let $\Bcal$ be an operator algebra with unity and let
$\theta:\Bcal\to B(\mathcal H)$ be a complex-linear homomorphism preserving the identity, where $\mathcal H\ne\{0\}$ and
$\dim\theta(\Bcal)\leq2$. Let $\alpha:\Bcal\to\Bcal$ be an antilinear map preserving the identity. If $\Phi_{\theta,\alpha}$ is contractive, then
\[
 \begin{aligned}
 \norm{\theta}_{\mathrm{cb}}&=\norm{\theta}\leq2,\\
 \norm{\theta}_{w_\rho,\mathrm{cb}}
 &=\kappa_\rho(\norm{\theta})\leq c_\rho
 \qquad(\rho\geq1).
 \end{aligned}
\]
These conclusions apply, in particular, to homomorphisms into
$M_2(\C)$ with commutative range. The constants $2$ and $c_\rho$ are
best possible.
\end{corollary}

\begin{proof}
An algebra of dimension at most two with unity is commutative, so
Theorem~\ref{v3:scalar-automatic} gives boundedness and
$\norm{\theta}\leq2$. Meyer proved that a bounded homomorphism of linear rank at most two preserving the identity satisfies
$\norm{\theta}_{\mathrm{cb}}=\norm{\theta}$
\cite[Corollary~3.6]{Meyer1997}. The radius formula now follows from
\cite[Theorem~3.2(b)]{BadeaCrouzeixKlaja2018}.
A commutative subalgebra of $M_2(\C)$ containing $I_2$ has dimension at most two:
if it contains a nonscalar matrix $X$, it is contained in the commutant
$\operatorname{span}\{I,X\}$ of $X$. The scalar case is immediate.
\end{proof}

To see sharpness of all the preceding constants, take
$\Bcal=A(\D)$ and define
\[
 J=\begin{pmatrix}0&2\\0&0\end{pmatrix},
 \qquad
 \theta(f)=f(0)I+f'(0)J,
 \qquad
 \alpha(f)(z)=\overline{f(0)}.
\]
Since $J^2=0$, the map $\theta$ preserves the identity and is a homomorphism with two-dimensional range, and
\[
 \Phi_{\theta,\alpha}(f)=f(J/2).
\]
Since $J/2$ is a contraction, this map is completely contractive.
The similarity $J=D(J/2)D^{-1}$ with $D=\diag(2,1)$ gives
$\norm{\theta}_{\mathrm{cb}}\leq2$, while $\norm{\theta(z)}=2$.
Thus $\norm{\theta}=\norm{\theta}_{\mathrm{cb}}=2$, and the radius
formulas give equality with $c_\rho$. The map
$P=2\Phi_{\theta,\alpha}$ also satisfies
Corollary~\ref{v3:real-positive}: the contraction
$\Phi_{\theta,\alpha}$ preserves the identity and is real positive by
\cite[Corollary~2.2]{BlecherNeal2021}, so $P$ is real positive as well.
The required commutation follows from the commutativity of
$\theta(A(\D))$.

We next prove the column and row estimates in
Theorem~\ref{thm:intro-abstract-complete-three} for quotients of uniform
algebras.

\begin{definition}\label{v3:def-Q}
A Q-algebra with unity is an operator algebra $\Bcal$ with unity that is
completely isometrically isomorphic to $\Acal/I$, where
$\Acal\subset C(X)$ is a uniformly closed algebra containing the constant function $1$, where $X$ is a compact Hausdorff space and $I$ is a proper closed ideal. The matrix norms
are the quotient matrix norms. Thus, with $q_I:\Acal\to\Acal/I$,
\[
 \norm{[b_{ij}]}_{M_{s,r}(\Acal/I)}
 =\inf\left\{\sup_{x\in X}\norm{[f_{ij}(x)]}:
       q_I(f_{ij})=b_{ij}\right\}.
\]
\end{definition}

This is the terminology of Meyer \cite[Section~1]{Meyer1997}.
It includes uniform algebras, with $I=0$, and quotients such as
$A(\D)/z^2A(\D)$. The latter contains a nonzero nilpotent element
and hence is not a uniform algebra.

For $q\geq1$, let $A(\B_q)$ denote the ball algebra on the closed unit
ball of $\C^q$.

\begin{lemma}\label{lem:abstract-ball-dilation}
Let $\Psi:A(\B_q)\to B(\mathcal H)$ be a complete contraction preserving the identity.
There are a Hilbert space $\mathcal K$, an isometry
$V:\mathcal H\to\mathcal K$, and commuting normal operators
$Q_1,\ldots,Q_q\in B(\mathcal K)$ such that
\[
        \sum_{j=1}^qQ_j^*Q_j
        =\sum_{j=1}^qQ_jQ_j^*=I
\]
and
\[
        \Psi(h)=V^*h(Q_1,\ldots,Q_q)V,
        \qquad h\in A(\B_q).
\]
\end{lemma}

\begin{proof}
Regard the ball algebra as a subalgebra of $C(\partial\B_q)$.
Arveson's theorem \cite[Theorem~1.2.9]{Arveson1969} gives a completely
positive extension $\widetilde\Psi:C(\partial\B_q)\to B(\mathcal H)$
preserving the identity. By \cite[Theorem~1 and Section~3]{Stinespring1955},
there are a $*$-representation preserving the identity
\[
        \pi:C(\partial\B_q)\longrightarrow B(\mathcal K)
\]
and an isometry $V:\mathcal H\to\mathcal K$ such that
$\Psi(h)=V^*\pi(h)V$ for $h\in A(\B_q)$.
\par
Put
$Q_j=\pi(z_j)$, where $z_j$ is the $j$th coordinate function. The
operators $Q_1,\ldots,Q_q$ commute and are normal. Since
$\sum_j|z_j|^2=1$ on $\partial\B_q$, they satisfy both spherical
identities. The compression formula holds first for
polynomials and then for every $h\in A(\B_q)$ by uniform approximation.
\end{proof}

Let $\mathcal W_n$ be the set of words of length $n$ in
$\{1,\ldots,q\}$. For
$\omega=(\omega_1,\ldots,\omega_n)\in\mathcal W_n$, put
\[
 T_\omega=T_{\omega_1}\cdots T_{\omega_n},
 \qquad
 Q_\omega=Q_{\omega_1}\cdots Q_{\omega_n}.
\]
If $j\in\{1,\ldots,q\}$, then $\omega j$ denotes the word obtained by
adjoining $j$ to $\omega$.

\begin{lemma}\label{lem:abstract-tuple-perturbation}
Let $T_1,\ldots,T_q\in B(\mathcal H)$ be commuting operators. Let
$Q_1,\ldots,Q_q\in B(\mathcal K)$ be a commuting normal tuple
satisfying
\[
        \sum_{j=1}^qQ_j^*Q_j=I,
\]
let $V:\mathcal H\to\mathcal K$ be an isometry, and put
\[
        E_\omega=2V^*Q_\omega^*V-T_\omega^*.
\]
Assume that every $E_\omega$ commutes with every $T_j$. For $n\geq1$,
define
\[
 \mathcal T_n:\mathcal H\longrightarrow\mathcal H^{\mathcal W_n},
 \qquad
 \mathcal T_nx=(T_\omega x)_{\omega\in\mathcal W_n},
\]
and
\[
 \mathcal E_n:\mathcal H^{\mathcal W_n}\longrightarrow\mathcal H,
 \qquad
 \mathcal E_n(x_\omega)_{\omega\in\mathcal W_n}
 =\sum_{\omega\in\mathcal W_n}E_\omega x_\omega.
\]
If
\[
\sup_{n\geq1}\norm{\mathcal E_n}<\infty,
\]
then
\[
        \norm{\begin{pmatrix}T_1\\ \vdots\\ T_q\end{pmatrix}}
        \leq2.
\]
\end{lemma}

\begin{proof}
Let $\mathcal A_n$ be the row operator with entries
$V^*Q_\omega^*V$, $\omega\in\mathcal W_n$. The spherical identity gives
\[
 \sum_{\omega\in\mathcal W_n}Q_\omega^*Q_\omega=I.
\]
Thus $\mathcal A_n$, being a compression of the corresponding row,
has norm at most one. Since
\[
       \mathcal E_n=2\mathcal A_n-\mathcal T_n^*,
\]
one has $\norm{\mathcal T_n}\leq2+\norm{\mathcal E_n}$.
The uniform bound for $\mathcal T_n$ therefore follows. Conversely,
$\norm{\mathcal E_n}\leq2+\norm{\mathcal T_n}$, so either one of
these bounds implies the other.
\par

Put
\[
        \kappa
        =\norm{\begin{pmatrix}T_1\\ \vdots\\ T_q\end{pmatrix}}
\]
and suppose that $\kappa>2$. If the column operator does not attain
its norm, choose unit vectors $x_m$ on which its norm tends to
$\kappa$ and pass to Hilbert-space ultrapowers with respect to a free
ultrafilter. The class of $(x_m)$ is a unit vector at which the
ultrapower of $\mathcal T_1$ attains the norm $\kappa$. The defining
relations and the uniform bounds for $\mathcal T_n$ and
$\mathcal E_n$ pass to the ultrapowers. Since
$\mathcal T_1^*\mathcal T_1\preceq\kappa^2I$, equality of its
quadratic form at this vector gives the equation below. It is
therefore enough to treat the case in which there is a unit vector
$x\in\mathcal H$ such that
\[
        \sum_{j=1}^qT_j^*T_jx=\kappa^2x.
\]

For $n\geq1$, define
\[
        m_n
        =\sum_{\omega\in\mathcal W_n}
        \Ree\ip{E_\omega T_\omega x}{x}.
\]
The sequence $(m_n)$ is bounded because
\[
        |m_n|
        \leq\norm{\mathcal E_n}\norm{\mathcal T_n}.
\]
Set
\[
        y=\sum_{j=1}^qQ_j^*VT_jx,
        \qquad
        u=y-\kappa Vx.
\]
The row $(Q_1^*,\ldots,Q_q^*)$ is a contraction, so
$\norm{y}\leq\kappa$. Moreover,
\[
        2\Ree\ip{y}{Vx}=\kappa^2+m_1.
\]
Consequently,
\begin{equation}\label{eq:abstract-u-upper}
        \norm{u}^2
        \leq2\kappa^2-\kappa m_1-\kappa^3.
\end{equation}

For every $\omega\in\mathcal W_n$, commutativity of the $T_j^*$ gives
\[
        T_{\omega j}^*T_jx=T_\omega^*T_j^*T_jx.
\]
It follows that
\[
 \sum_{j=1}^qE_{\omega j}T_jx
 =2V^*Q_\omega^*y-\kappa^2T_\omega^*x.
\]
Using $u=y-\kappa Vx$, we obtain
\[
 \left(\kappa E_\omega-
 \sum_{j=1}^qE_{\omega j}T_j\right)x
 =\kappa(\kappa-1)T_\omega^*x
  -2V^*Q_\omega^*u.
\]
Take real parts of the inner products with $T_\omega^*x$ and sum over
$\omega\in\mathcal W_n$. Since every $E_\omega$ commutes with every
$T_j$, the two terms on the left sum respectively to
$\kappa m_n$ and $m_{n+1}$. Completing the square termwise and using
that $V^*$ is a contraction gives
\begin{align}
 \kappa m_n-m_{n+1}
 &\geq
 -\frac1{\kappa(\kappa-1)}
 \sum_{\omega\in\mathcal W_n}\norm{V^*Q_\omega^*u}^2 \notag\\
 &\geq
 -\frac1{\kappa(\kappa-1)}
 \sum_{\omega\in\mathcal W_n}\norm{Q_\omega^*u}^2 \notag\\
 &=-\frac{\norm{u}^2}{\kappa(\kappa-1)}.
 \label{eq:abstract-perturbation-recurrence}
\end{align}
A commuting family of normal operators is doubly commuting. Hence
\[
        \sum_{\omega\in\mathcal W_n}Q_\omega Q_\omega^*
        =\left(\sum_{j=1}^qQ_jQ_j^*\right)^n=I.
\]
After division of
\eqref{eq:abstract-perturbation-recurrence} by $\kappa^n$ and
summation from $n=1$ to $N$, we get
\[
        m_1-\frac{m_{N+1}}{\kappa^N}
        \geq
        -\frac{\norm{u}^2}{\kappa(\kappa-1)}
        \sum_{n=1}^N\frac1{\kappa^n}.
\]
Letting $N\to\infty$ and using boundedness of $(m_n)$ gives
\[
        \kappa m_1
        \geq-\frac{\norm{u}^2}{(\kappa-1)^2}.
\]
Together with \eqref{eq:abstract-u-upper}, this yields
\[
        \norm{u}^2
        \left(1-\frac1{(\kappa-1)^2}\right)
        \leq\kappa^2(2-\kappa),
\]
which is impossible when $\kappa>2$.
\end{proof}

\begin{remark}
For $q=1$, Lemma~\ref{lem:abstract-tuple-perturbation} reduces to the
unitary-dilation case of the perturbation lemma of Lorist and
Schwenninger \cite[Lemma~1]{LoristSchwenninger2026}. The proof above
adds the summation over words needed for tuples.
\end{remark}

\begin{lemma}\label{v3:Q-ball}
Let $\Bcal$ be a Q-algebra with unity. Transposition is an isometry from
$M_{s,r}(\Bcal)$ onto $M_{r,s}(\Bcal)$ for all $r,s\geq1$.
If $b_1,\ldots,b_q\in\Bcal$ satisfy
\[
 \norm{(b_1,\ldots,b_q)^{\mathsf T}}\leq1,
\]
then the polynomial calculus extends to a complete contraction preserving the identity
\[
 C_b:A(\B_q)\longrightarrow\Bcal,
 \qquad C_b(h)=h(b_1,\ldots,b_q).
\]
For every $n\geq1$, the row and column formed by
$b_\omega=b_{\omega_1}\cdots b_{\omega_n}$,
$\omega\in\mathcal W_n$, both have norm at most one.
In particular, $\Bcal$ satisfies von Neumann's inequality.
\end{lemma}

\begin{proof}
The transpose assertion is \cite[Corollary~5.3]{Meyer1997}; its
rectangular form follows by adjoining zero rows and columns.
Identify $\Bcal$ with $\Acal/I$. For every $\varepsilon>0$, the
quotient column norm gives simultaneous lifts $f_j\in\Acal$ with
\[
 q_I(f_j)=b_j,
 \qquad
 \sup_{x\in X}\sum_{j=1}^q|f_j(x)|^2\leq(1+\varepsilon)^2.
\]
For any matrix-valued polynomial $P$ in $q$ commuting variables,
complete contractivity of $q_I$ gives
\[
 \norm{P(b_1,\ldots,b_q)}
 \leq\sup_{\norm{z}\leq1+\varepsilon}\norm{P(z)}.
\]
Letting $\varepsilon\to0^+$ and using polynomial density in
$A(\B_q)$ proves the asserted extension at every matrix level.
Finally,
\[
 \sum_{\omega\in\mathcal W_n}|z_\omega|^2
 =\left(\sum_{j=1}^q|z_j|^2\right)^n\leq1
 \qquad(z\in\overline{\B_q}).
\]
Apply the rectangular amplifications of $C_b$ to the monomial row and
column. The case $q=1$ gives von Neumann's inequality.
\end{proof}

\begin{theorem}\label{thm:abstract-column-row}
Let $\Bcal$ be a Q-algebra with unity, let
$\theta:\Bcal\to B(\mathcal H)$ be a completely bounded homomorphism preserving the identity, and let $\alpha:\Bcal\to\Bcal$ be an antilinear map preserving the identity. Suppose that $\Phi_{\theta,\alpha}$ is completely contractive.
If $b_1,\ldots,b_q\in\Bcal$ satisfy
\[
    \norm{(b_1,\ldots,b_q)^{\mathsf T}}\leq1,
\]
then
\begin{equation}\label{eq:abstract-column-row}
    \norm{\begin{pmatrix} \theta(b_1)\\ \vdots\\ \theta(b_q)
    \end{pmatrix}} \leq2,\qquad
    \norm{\begin{pmatrix}
    \theta(b_1)&\cdots&\theta(b_q)
    \end{pmatrix}}\leq2.
\end{equation}
No boundedness assumption on $\alpha$ is required.
\end{theorem}

\begin{proof}
By Lemma~\ref{v3:Q-ball},
$\Psi=\Phi_{\theta,\alpha}\circ C_b$ is a complete contraction preserving the identity.
Let $V:\mathcal H\to\mathcal K$ and $(Q_1,\ldots,Q_q)$ be the
normal spherical dilation from Lemma~\ref{lem:abstract-ball-dilation}.

Set $T_j=\theta(b_j)$ and $L=2\Phi_{\theta,\alpha}-\theta$.
Then $L$ is completely bounded,
$\norm{L}_{\mathrm{cb}}\leq2+\norm{\theta}_{\mathrm{cb}}$, and
\[
    L(b)=\theta(\alpha(b))^*.
\]
For $\omega\in\mathcal W_n$, write
$b_\omega=b_{\omega_1}\cdots b_{\omega_n}$. Since
$V^*Q_\omega V=\Phi_{\theta,\alpha}(b_\omega)$, we have
\[
    2V^*Q_\omega^*V-T_\omega^*=L(b_\omega)^*
    =\theta(\alpha(b_\omega)).
\]
These operators commute with every $T_j$ because $\Bcal$ is
commutative. Lemma~\ref{v3:Q-ball} gives
\[
    \norm{(T_\omega)_{\omega\in\mathcal W_n}^{\mathsf T}}
    \leq\norm{\theta}_{\mathrm{cb}},
    \qquad
    \norm{(L(b_\omega)^*)_{\omega\in\mathcal W_n}}
    =\norm{(L(b_\omega))_{\omega\in\mathcal W_n}^{\mathsf T}}
    \leq\norm{L}_{\mathrm{cb}}.
\]
Lemma~\ref{lem:abstract-tuple-perturbation} therefore provides the
column estimate.

For the row estimate, apply the same lemma to $S_j=T_j^*$ and
$R_j=Q_j^*$. Both tuples are commuting, and $(R_1,\ldots,R_q)$ is
again a normal spherical tuple. Since the $T_j$ commute, reversal of
the indexed products does not change the word operators. The
corresponding correction term is
\[
    2V^*R_\omega^*V-S_\omega^*=2V^*Q_\omega V-T_\omega
    =L(b_\omega).
\]
It commutes with every $S_j$. By Lemma~\ref{v3:Q-ball}, the word row
and column over $\Bcal$ both have norm at most one. Hence
\[
    \norm{(S_\omega)_{\omega\in\mathcal W_n}^{\mathsf T}}
    =\norm{(T_\omega)_{\omega\in\mathcal W_n}}
    \leq\norm{\theta}_{\mathrm{cb}},
\]
while
\[
    \norm{(L(b_\omega))_{\omega\in\mathcal W_n}}
    \leq\norm{L}_{\mathrm{cb}}.
\]
The perturbation lemma now gives the row estimate.
\end{proof}

The column and row estimates also control arbitrary rectangular
levels with a dimension-dependent constant.

\begin{lemma}\label{lem:abstract-rectangular-slice}
Under the hypotheses of
Theorem~\ref{thm:abstract-column-row}, for all positive integers
$r,s$,
\[
        \norm{\theta^{(s,r)}}
        \leq2\sqrt{\min\{r,s\}}.
\]
The same estimate holds for the homomorphism induced by $\theta$ on a quotient $\Bcal/J$, where $J$ is a closed ideal contained in $\ker\theta$ and the quotient is equipped with the operator-algebra matrix norms.
\end{lemma}

\begin{proof}
Let $F=[f_{ij}]\in M_{s,r}(\Bcal)$ satisfy $\norm{F}\leq1$, and put
$Z=\theta^{(s,r)}(F)$. Each column of $F$ has norm at most one.
Theorem~\ref{thm:abstract-column-row} therefore gives column norm at
most two for its image. For
$x=\sum_{j=1}^r e_j\otimes x_j\in\C^r\otimes\mathcal H$,
\[
 \norm{Zx}
 \leq2\sum_{j=1}^r\norm{x_j}
 \leq2\sqrt r\,\norm{x}.
\]
Thus $\norm{Z}\leq2\sqrt r$. Each row of $F$ also has norm at most
one. By Lemma~\ref{v3:Q-ball}, its transpose is a column contraction.
Applying the row estimate to the $s$ rows gives
\[
        \norm{Zx}^2\leq4s\norm{x}^2,
\]
so $\norm{Z}\leq2\sqrt s$.

For the quotient assertion, lift a rectangular matrix of quotient norm at most one to matrices over $\Bcal$ of norm at most $1+\varepsilon$, apply the preceding estimate after division by $1+\varepsilon$, and let $\varepsilon\to0^+$.
\end{proof}

\begin{theorem}\label{thm:abstract-complete-three}
Let $\Bcal$ be a Q-algebra with unity, let $\theta:\Bcal\to M_d(\C)$ be a homomorphism preserving the identity, and let $\alpha:\Bcal\to\Bcal$ be an antilinear map preserving the identity. Suppose that $\Phi_{\theta,\alpha}$ is completely contractive. If $d\leq3$, then
\begin{equation}\label{eq:abstract-complete-three}
    \norm{\theta}_{\mathrm{cb}}\leq2.
\end{equation}
More generally, for $d\geq2$,
\begin{equation}\label{eq:abstract-dimension-bound}
    \norm{\theta}_{\mathrm{cb}}
    \leq2\sqrt{\lfloor d/2\rfloor}.
\end{equation}
If $\norm{\theta}_{\mathrm{cb}}>2$, then
$\norm{\theta^{(s,r)}}=\norm{\theta}_{\mathrm{cb}}$ for some
$r,s\geq2$ with $r+s\leq d$.
\end{theorem}

\begin{proof}
Lemma~\ref{v3:Q-ball} shows that $\Bcal$ satisfies von Neumann's
inequality. Since $\Bcal$ is commutative and $\Phi_{\theta,\alpha}$
is contractive, Theorem~\ref{v3:scalar-automatic} gives
$\norm{\theta}\leq2$. Smith's lemma now shows that $\theta$ is
completely bounded. Let $J=\ker\theta$, $q:\Bcal\to\Bcal/J$ be the quotient map, and
$$
    \overline\theta:\Bcal/J\longrightarrow M_d(\C)
$$
be the induced injective homomorphism. At every rectangular level,
\begin{equation}\label{eq:abstract-quotient-level}
    \norm{\overline\theta^{(s,r)}}
    =\norm{\theta^{(s,r)}}.
\end{equation}
Indeed, one inequality follows from $\theta=\overline\theta\circ q$, while the reverse inequality follows by taking lifts whose norms approach the quotient norm. The quotient $\Bcal/J$ is finite-dimensional because it is algebraically isomorphic to the subalgebra $\theta(\Bcal)$ of $M_d(\C)$. All estimates below use approximate lifts in $\Bcal$; no invariance of $J$ under $\alpha$ is required.

Let $C=\norm{\theta}_{\mathrm{cb}}$. There is nothing to prove if $C\leq1$. If $C>1$, apply
Theorem~\ref{thm:abstract-rectangular-extremal} to $\overline\theta$. There are positive integers $r,s$ and $X\in M_{s,r}(\Bcal/J)$ such that
$$
    r+s\leq d,\qquad
    \norm{X}\leq1,\qquad
    \norm{\overline\theta^{(s,r)}(X)}=C.
$$
If $d=1$, the inequality $r+s\leq d$ is impossible, so $C\leq1$. If $d=2$ or $3$, the only possible pairs are $(r,s)=(1,1)$, $(1,2)$, and $(2,1)$ with only the first pair available when $d=2$. Using lifts of norm at most $1+\varepsilon$ and Theorem~\ref{thm:abstract-column-row}, every available level is bounded by $2(1+\varepsilon)$. Letting $\varepsilon\to0^+$ proves \eqref{eq:abstract-complete-three}.

For arbitrary $d\geq2$, Lemma~\ref{lem:abstract-rectangular-slice} and \eqref{eq:abstract-quotient-level} imply
$$
    C\leq2\sqrt{\min\{r,s\}}
    \leq2\sqrt{\lfloor d/2\rfloor},
$$
which proves \eqref{eq:abstract-dimension-bound}. Finally, if $C>2$, the detecting pair cannot have $r=1$ or $s=1$, by the scalar, column, and row estimates and the same approximate-lift argument. Hence $r,s\geq2$.
\end{proof}

The disk-algebra example following Corollary~\ref{v3:rank-two} also
proves sharpness in Theorem~\ref{thm:abstract-column-row} and in the
case $d=2$ of Theorem~\ref{thm:abstract-complete-three}. Taking the
direct sum with $f\mapsto f(0)$ gives the same sharp example for $d=3$.

Unlike Theorem~\ref{thm:tuple-square}, the abstract argument above
gives the endpoint column and row estimates but not corresponding
higher-order positive square-function series. We therefore do not claim an abstract analog of Theorem~\ref{thm:intro-complete-equality}.

\section{Noncommutative operator algebras}\label{v5:sec:noncommutative}
The scalar estimate in Theorem~\ref{v3:scalar-automatic} uses a
commutation hypothesis, and the complete Q-algebra estimates use
commutativity of the domain. We now give the counterexample in
Theorem~\ref{v5:thm:cor-counterexample} and prove
Theorem~\ref{v5:thm:digraph} for matrix algebras containing the diagonal.

For a linear or antilinear map $L$, we use the entrywise amplification
$L^{(k)}([b_{ij}])=[L(b_{ij})]$ and the norm
$\norm{L}_{\mathrm{cb}}=\sup_{k\geq1}\norm{L^{(k)}}$.
This is the convention in
\cite[Section~6]{ClouatreOstermannRansford2023}.
In particular,
\[
 [\Phi_{\theta,\alpha}(b_{ij})]
 =\frac12\bigl([\theta(b_{ij})]+[\theta(\alpha(b_{ji}))]^*\bigr),
\]
where the adjoint is taken on the whole block matrix.

We first define the algebra and the two maps, construct completely
positive maps proving the required complete contractivity, and then
establish the strict norm inequality.

Put
\begin{equation}\label{v5:eq:basis}
 P=\begin{pmatrix}1&-5/2&1\\0&0&0\\0&0&0\end{pmatrix},
 \qquad N=E_{12},\qquad I=I_3,
\end{equation}
and define
\begin{equation}\label{v5:eq:source}
 \begin{split}
 X(a,c,u)&=aP+c(I-P)+uN
 =\begin{pmatrix}a&u-\frac52(a-c)&a-c\\0&c&0\\0&0&c\end{pmatrix},\\
 \Acal&=\{X(a,c,u):a,c,u\in\C\}\subset M_3(\C).
 \end{split}
\end{equation}
The matrix norms on $\Acal$ are those inherited from $M_3(\C)$.
The matrices $I,P,N$ have real entries and form a complex basis of
$\Acal$. Define $\theta$ by complex-linear extension of
\begin{equation}\label{v5:eq:theta}
 \theta(I)=I_2,\qquad \theta(P)=E_{11},\qquad
 \theta(N)=\frac{16}{19}E_{12},
\end{equation}
where the matrix units on the right belong to $M_2(\C)$. Define $\alpha$
by antilinear extension of
\begin{equation}\label{v5:eq:alpha}
 \alpha(I)=I,\qquad
 \alpha(P)=\frac7{17}P+\frac38N,\qquad
 \alpha(N)=\frac4{21}(I-P).
\end{equation}
In particular,
\[
 \theta(X(a,c,u))=\begin{pmatrix}a&16u/19\\0&c\end{pmatrix},
 \qquad \alpha(\Acal)\subset\Acal.
\]

\begin{lemma}\label{v5:lem:algebra}
The space $\Acal$ is a three-dimensional noncommutative operator algebra
with unity, and $\Acal\cap\Acal^*=\C I$. The map $\theta$ is an injective
homomorphism onto $T_2$, preserves the identity, and satisfies
\[
 \norm\theta_{\mathrm{cb}}\le\frac{8\sqrt{29}}{19}.
\]
\end{lemma}
\begin{proof}
Direct multiplication gives
\[
 P^2=P,\qquad PN=N,\qquad NP=N^2=0,
\]
and hence
\begin{equation}\label{v5:eq:product}
 X(a,c,u)X(a',c',u')=X(aa',cc',au'+uc').
\end{equation}
These identities prove closure under multiplication and multiplicativity
of $\theta$. The coordinate formula proves injectivity, surjectivity
onto $T_2$, and preservation of the identity. If $b=X(a,c,u)$ and
$b^*\in\Acal$, the lower triangular entries of $b^*$ force $a=c$ and
$u=0$. Thus $b=cI$.

The plane $\mathcal M=\operatorname{span}\{e_1,2e_2+5e_3\}$ is
invariant for $\Acal$. In the orthonormal basis
$e_1,(2e_2+5e_3)/\sqrt{29}$, restriction to $\mathcal M$ gives the
completely contractive homomorphism
\[
 \pi(X(a,c,u))=\begin{pmatrix}a&2u/\sqrt{29}\\0&c\end{pmatrix}.
\]
With $L=\operatorname{diag}(8\sqrt{29}/19,1)$, we have
$\theta(b)=L\pi(b)L^{-1}$. Since $8\sqrt{29}/19>1$, this gives the
asserted bound at every matrix level.
\end{proof}

We next construct a completely positive extension of the symmetrized
map. Set
\begin{equation}\label{v5:eq:frame}
 v=(4,0,-1)^T,\qquad w=(1,2,4)^T,
\end{equation}
and define
\begin{equation}\label{v5:eq:Phi_extension}
 \Gamma_{\Phi}(Z)=\begin{pmatrix}
 v^*Zv/17&v^*Zw/19\\
 w^*Zv/19&w^*Zw/21
 \end{pmatrix}
 \qquad(Z\in M_3(\C)).
\end{equation}

\begin{lemma}\label{v5:lem:phi_extension}
The map $\Gamma_{\Phi}:M_3(\C)\to M_2(\C)$ is completely positive, preserves the
identity, and satisfies $\Gamma_{\Phi}|_{\Acal}=\Phi_{\theta,\alpha}$.
\end{lemma}
\begin{proof}
The coefficient matrix
\begin{equation}\label{v5:eq:small_Gram}
 G=\begin{pmatrix}1/17&1/19\\1/19&1/21\end{pmatrix}
\end{equation}
is positive definite because $19^2-17\cdot21=4>0$.
If $G=\sum_j g_jg_j^*$, its Schur multiplier is
\[
 Z\longmapsto\sum_j\operatorname{diag}(g_j)Z
                           \operatorname{diag}(g_j)^*.
\]
This is a completely positive map on $M_2(\C)$. Its composition with the
completely positive map $Z\mapsto(v,w)^*Z(v,w)$ is $\Gamma_{\Phi}$.
The identities $v^*w=0$, $\norm v^2=17$, and $\norm w^2=21$ give
$\Gamma_{\Phi}(I)=I_2$.

Evaluation on the basis yields
\begin{equation}\label{v5:eq:Phi_basis}
 \Gamma_{\Phi}(P)=\begin{pmatrix}12/17&0\\3/19&0\end{pmatrix},\qquad
 \Gamma_{\Phi}(N)=\begin{pmatrix}0&8/19\\0&2/21\end{pmatrix}.
\end{equation}
On the other hand,
\[
 \theta(\alpha(P))=\begin{pmatrix}7/17&6/19\\0&0\end{pmatrix},\qquad
 \theta(\alpha(N))=\begin{pmatrix}0&0\\0&4/21\end{pmatrix}.
\]
Together with \eqref{v5:eq:theta}, these give
\begin{equation}\label{v5:eq:Phi_restriction}
 \Gamma_{\Phi}|_{\Acal}=\Phi_{\theta,\alpha},
\end{equation}
since both sides are complex-linear and agree on $I,P,N$.
\end{proof}

To prove complete contractivity of $\alpha$, we construct a positive
Choi matrix. For a complex-linear map $\Gamma:M_3(\C)\to M_m(\C)$, write
\[
 C_\Gamma=[\Gamma(E_{ij})]_{i,j=1}^3.
\]
The blocks are indexed by the matrix units of $M_3(\C)$. If $C=[C_{ij}]$
is positive with $m\times m$ blocks, the associated map is
$\Gamma_C(Z)=\sum_{i,j}z_{ij}C_{ij}$.

\begin{lemma}\label{v5:lem:cp_completion}
Let $U,V\in M_m(\C)$. Suppose that $M\in M_m(\C)$ is positive definite and
\begin{equation}\label{v5:eq:completion_residual}
 Q=I-M-V^*M^{-1}V-(U-M)^*M^{-1}(U-M)\ge0.
\end{equation}
There is a completely positive map $\Gamma:M_3(\C)\to M_m(\C)$ preserving the
identity such that
\[
 \Gamma(E_{11}+E_{13})=U,\qquad \Gamma(E_{12})=V.
\]
\end{lemma}
\begin{proof}
Put $C=U-M$ and define
\begin{equation}\label{v5:eq:completion_Gram}
 C_\Gamma=
 \begin{pmatrix}M\\V^*\\C^*\end{pmatrix}
 M^{-1}\begin{pmatrix}M&V&C\end{pmatrix}
 +\operatorname{diag}(0,0,Q).
\end{equation}
This matrix is positive. Choi's criterion
\cite[Theorem~2]{Choi1975} gives a completely positive map.
The diagonal blocks sum to $I$, so the map preserves the identity.
The $(1,1)$, $(1,2)$, and $(1,3)$ blocks of $C_\Gamma$
are $M,V,C$, which prove the interpolation identities.

For completeness, this implication follows from a Gram factorization.
Write $C_\Gamma=\sum_\nu z_\nu z_\nu^*$,
with $z_\nu=(z_{\nu,1},z_{\nu,2},z_{\nu,3})^T$ and
$z_{\nu,i}\in\C^m$. Let $K_\nu:\C^m\to\C^3$ have $i$th row
$z_{\nu,i}^*$. Then
\[
 \Gamma(Z)=\sum_\nu K_\nu^*ZK_\nu,\qquad
 \sum_\nu K_\nu^*K_\nu=I_m.
\]
The resulting isometric column factorization also proves complete
contractivity at every matrix level.
\end{proof}

\begin{lemma}\label{v5:lem:alpha_extension}
There is a completely positive map $\Gamma_{\alpha}:M_3(\C)\to M_3(\C)$ preserving the
identity such that
\begin{equation}\label{v5:eq:alpha_restriction}
 \alpha(b)=\Gamma_{\alpha}(\overline b)\qquad(b\in\Acal),
\end{equation}
where the bar denotes entrywise conjugation. Thus $\Gamma_{\alpha}$ extends the
complex-linear map $b\mapsto\alpha(\overline b)$ on $\Acal$.
\end{lemma}
\begin{proof}
Apply Lemma~\ref{v5:lem:cp_completion} with
\begin{equation}\label{v5:eq:completion_data}
 \begin{split}
 U&=\alpha(P)+\frac52\alpha(N)
   =\frac7{17}P+\frac38N+\frac{10}{21}(I-P),\\
 V&=\alpha(N)=\frac4{21}(I-P),\qquad
 M=\frac16\begin{pmatrix}4&1&0\\1&2&0\\0&0&2\end{pmatrix}.
 \end{split}
\end{equation}
The matrix $M$ is positive definite and
\[
 M^{-1}=\begin{pmatrix}12/7&-6/7&0\\-6/7&24/7&0\\0&0&3\end{pmatrix}.
\]
For the matrix $Q$ in \eqref{v5:eq:completion_residual}, put
$D_0=\operatorname{diag}(1,2,1)$. Direct calculation gives
\begin{equation}\label{v5:eq:small_certificate}
 D_0QD_0=\frac1{1189524}
 \begin{pmatrix}
 237356&-18081&-22540\\
 -18081&454919&374228\\
 -22540&374228&508268
 \end{pmatrix}.
\end{equation}
The diagonal entries of the integer matrix exceed the sums of the
absolute off-diagonal entries in their rows by
\[
 196735,\qquad62610,\qquad111500,
\]
respectively. For a Hermitian matrix $H=[h_{ij}]$,
\[
 z^*Hz\ge\sum_i\left(h_{ii}-\sum_{j\ne i}|h_{ij}|\right)|z_i|^2.
\]
This follows from $2|z_i z_j|\le |z_i|^2+|z_j|^2$. It proves positive
definiteness in \eqref{v5:eq:small_certificate}, and hence $Q>0$.

The resulting map satisfies
\[
 \Gamma_{\alpha}(P)=U-\frac52V=\alpha(P),\qquad
 \Gamma_{\alpha}(N)=V=\alpha(N),\qquad\Gamma_{\alpha}(I)=I,
\]
because $P+\frac52N=E_{11}+E_{13}$. Since $I,P,N$ have real entries,
complex linearity of $\Gamma_{\alpha}$ and antilinearity of $\alpha$ give
\eqref{v5:eq:alpha_restriction}.
\end{proof}

\begin{proof}[Proof of Theorem~\ref{v5:thm:cor-counterexample}]
Use \eqref{v5:eq:source}--\eqref{v5:eq:alpha}.
Lemma~\ref{v5:lem:algebra} proves the algebraic assertions and the complete
upper bound. Lemmas~\ref{v5:lem:phi_extension} and
\ref{v5:lem:alpha_extension} give identity-preserving completely positive
maps. They are complete contractions
\cite[Proposition~1.2.10]{Arveson1969}. At matrix level $k$, conjugating
each input block is entrywise conjugation of the whole scalar matrix
in $M_{3k}$, and therefore preserves its norm. Consequently,
\[
 \norm{[\alpha(b_{ij})]}
 =\norm{\Gamma_{\alpha}^{(k)}([\overline{b_{ij}}])}\le\norm{[b_{ij}]},
 \qquad
 \norm{[\Phi_{\theta,\alpha}(b_{ij})]}\le\norm{[b_{ij}]}.
\]
Both maps preserve the identity, so their completely bounded norms
are equal to one.

For the lower estimate, take the integer matrix
\begin{equation}\label{v5:eq:witness}
 B=2P+6N-I=\begin{pmatrix}1&1&2\\0&-1&0\\0&0&-1\end{pmatrix},
 \qquad
 \theta(B)=\begin{pmatrix}1&96/19\\0&-1\end{pmatrix}.
\end{equation}
For $t\ge0$, set
\[
 f(t)=\left\|\begin{pmatrix}1&t\\0&-1\end{pmatrix}\right\|
     =\frac{\sqrt{4+t^2}+t}{2}.
\]
A unitary change of basis in the last two coordinates of $\C^3$ gives
\begin{equation}\label{v5:eq:norms}
 \norm B=\frac{3+\sqrt5}{2},\qquad
 \norm{\theta(B)}=\frac{48+\sqrt{2665}}{19}.
\end{equation}
For $k>1$, the inequality $f(t)>k$ is equivalent to $t>k-k^{-1}$.
With $k=3+\sqrt5$,
\[
 k-k^{-1}=\frac{9+5\sqrt5}{4}<\frac{96}{19},
\]
since $213>0$ and
\begin{equation}\label{v5:eq:rational_margin}
 213^2-125\cdot19^2=244>0.
\end{equation}
Thus $\norm{\theta(B)}>2\norm B$. The ratio in \eqref{v5:eq:norms}
proves the lower bound in \eqref{v5:eq:cor-bound}.

Finally, Theorem~\ref{thm:abstract-rectangular-extremal} applies to
$\theta:\Acal\to M_2(\C)$. Since $\norm\theta_{\mathrm{cb}}>1$, it gives an extremal
at a rectangular level with $r,s\geq1$ and $r+s\leq2$. Thus $r=s=1$,
and $\norm\theta=\norm\theta_{\mathrm{cb}}$.
\end{proof}

The algebra is a nonunitary similarity of a triangular algebra. Put
\[
 S_0=I+\frac52E_{12}-E_{13},\qquad
 \mathcal A_0=
 \left\{\begin{pmatrix}a&u&0\\0&c&0\\0&0&c\end{pmatrix}:
 a,c,u\in\C\right\}.
\]
Since $(S_0-I)^2=0$, direct multiplication gives
\[
 \Acal=S_0\mathcal A_0S_0^{-1}.
\]
The similarity preserves the multiplication but changes the matrix
norms. The algebraic identification
\[
 X(a,c,u)\longmapsto\begin{pmatrix}a&u\\0&c\end{pmatrix}
\]
is not isometric, since $\norm P=\sqrt{33}/2$ whereas the image of $P$
has norm one.

The vectors in \eqref{v5:eq:frame} satisfy $Nv=0$ and $Pw=0$.
The symmetrization identities then determine \eqref{v5:eq:alpha}.
Replacing $19$ by $t>0$ in the off-diagonal entries of
\eqref{v5:eq:Phi_extension}, and $16/19$ by $16/t$ in
\eqref{v5:eq:theta}, leaves $\alpha$ unchanged. Complete positivity and
the strict norm violation for $B$ hold simultaneously when
\[
 \sqrt{357}\leq t<\frac{384}{9+5\sqrt5}.
\]
The integer $19$ belongs to this interval by
$19^2-357=4$ and \eqref{v5:eq:rational_margin}.

The construction was obtained by semidefinite programming.
With $\theta$ in \eqref{v5:eq:theta} fixed, preservation of the identity,
the range condition for $\alpha$, and the symmetrization identities on
$I,P,N$ impose affine constraints over the real scalars on two Choi
matrices. For a Hermitian matrix $H=H_1+iH_2$ with $H_1,H_2$ real,
\[
 H\geq0\quad\Longleftrightarrow\quad
 \begin{pmatrix}H_1&-H_2\\H_2&H_1\end{pmatrix}\geq0.
\]
The positivity constraints therefore give a semidefinite feasibility
problem in the form of \cite[Section~1.1]{VandenbergheBoyd1996}.
The displayed rational matrices give an exact verification independent
of the numerical calculation.

Since $\theta$ is injective, its kernel is invariant under $\alpha$.
Thus kernel invariance alone does not imply the bound. Moreover,
\begin{equation}\label{v5:eq:commutator}
 [\theta(P),\theta(\alpha(P))]
 =\begin{pmatrix}0&6/19\\0&0\end{pmatrix}\ne0,
 \qquad B^2=I_3.
\end{equation}
This also shows where the example fails the commutation hypothesis of
Theorem~\ref{v3:scalar-automatic}. Thus the norm inequality fails on
an involution.

\begin{corollary}\label{v5:cor:minimal}
Among counterexamples to the complete COR conjecture, the smallest
possible target matrix order is two and the smallest possible complex
dimension of the domain algebra is three.
\end{corollary}
\begin{proof}
A scalar target, or a domain of complex dimension at most two, gives
$\dim\theta(\Bcal)\leq2$.
Corollary~\ref{v3:rank-two} therefore excludes a counterexample in either
case. Theorem~\ref{v5:thm:cor-counterexample} attains both minimum dimensions.
\end{proof}

We next prove the bound for digraph algebras. The following consequence
of Lemma~\ref{lem:potapov} gives the correction estimate needed in the
proof.

\begin{lemma}\label{v5:lem:positive_correction}
Let $\Bcal$ be a finite-dimensional operator algebra with unity,
and let $D\subset\Bcal$ be a $C^*$-subalgebra with the same identity.
Let $E:\Bcal\to D$ be complex-linear, and let
$\theta:\Bcal\to B(\mathcal H)$
be a homomorphism preserving the identity whose restriction to $D$ is a
$*$-homomorphism. If
\[
 \Psi=\frac12(\theta+\theta E)
\]
is completely contractive, then $\norm\theta_{\mathrm{cb}}\leq2$.
\end{lemma}
\begin{proof}
Represent $\Bcal$ completely isometrically and preserving the identity as an operator algebra
on a Hilbert space. The finite-dimensional argument below uses
compactness of the unit balls of $M_m(\Bcal)$ and Lemma~\ref{lem:potapov};
it does not require a finite-dimensional ambient representation.
\par
First suppose that $\mathcal H$ is finite-dimensional and nonzero.
Fix $m\geq1$ and put $c=\norm{\theta^{(m)}}$.
If $c\leq1$, there is nothing to prove. Otherwise, finite dimensionality
gives a contraction $F\in M_m(\Bcal)$ and a unit vector
$x\in\mathcal H^m$ such that, with
\[
 R=\theta^{(m)}(F),\qquad
 K=\theta^{(m)}(E^{(m)}(F)),
\]
one has $\norm{Rx}=c$.
Apply \eqref{v5:eq:diagonal-orthogonality} with
$B=E^{(m)}(F)\in M_m(D)$. It gives $\ip{Rx}{Kx}=0$, and hence
\[
 c^2+\norm{Kx}^2=\norm{(R+K)x}^2\leq4,
\]
since $(R+K)/2=\Psi^{(m)}(F)$ is a contraction.
Thus $\norm{\theta^{(m)}}\leq2$ at every matrix level.

Now let $\mathcal H$ be arbitrary. For a unit vector
$x=(x_1,\ldots,x_m)\in\mathcal H^m$, put
\[
 \mathcal K=\operatorname{span}
 \{\theta(b)x_j:b\in\Bcal,\ 1\leq j\leq m\}.
\]
This space is finite-dimensional, contains every $x_j$, and is invariant
under $\theta(\Bcal)$ by multiplicativity. Since
$\theta(d)^*=\theta(d^*)$ for $d\in D$, it reduces $\theta(D)$.
The restriction $\theta_{\mathcal K}(b)=\theta(b)|_{\mathcal K}$ is
therefore a homomorphism preserving the identity, and its restriction
to $D$ is a $*$-homomorphism. The map
\[
 \frac12(\theta_{\mathcal K}+\theta_{\mathcal K}E)
 =P_{\mathcal K}\Psi(\,\cdot\,)|_{\mathcal K}
\]
is completely contractive. The finite-dimensional case gives
$\norm{\theta_{\mathcal K}}_{\mathrm{cb}}\leq2$.
For every $F\in M_m(\Bcal)$, invariance now yields
\[
 \norm{\theta^{(m)}(F)x}
 =\norm{\theta_{\mathcal K}^{(m)}(F)x}
 \leq2\norm{F}.
\]
Taking the supremum over unit vectors and matrix levels proves the
assertion. The case $\mathcal H=\{0\}$ is immediate.
\end{proof}

For the finite-dimensional diagonal, we first use an argument with
central projections as in \cite[Section~5]{ClouatreOstermannRansford2023}.
We then use states to treat the noncommutative summands with arbitrary
representation multiplicities.

\begin{lemma}\label{v5:lem:positive_diagonal}
Let $D$ be a finite-dimensional $C^*$-algebra with identity. Let
$\eta:D\to B(\mathcal H)$ be a homomorphism preserving the identity,
and let $\gamma:D\to D$ be an antilinear contraction preserving the
identity. If
\[
 \Xi(d)=\frac12\bigl(\eta(d)+\eta(\gamma(d))^*\bigr)
\]
is positive, then $\eta$ is a $*$-homomorphism.
\end{lemma}
\begin{proof}
The assertion is immediate if $\mathcal H=\{0\}$, so assume otherwise.
The map $\eta$ is bounded because $D$ is finite-dimensional.
The complex-linear map $d\mapsto\gamma(d^*)$ is a contraction
preserving the identity on $D$. It is therefore positive and preserves
adjoints \cite[Proposition~1.2.8]{Arveson1969}
(see also \cite[Lemma~1.2]{BlecherLabuschagne2003}). Consequently,
$\gamma$ is positive and preserves adjoints.

First let $p$ be a central projection of $D$, and set
$U=\eta(2p-I)$ and $V=\eta(\gamma(2p-I))$.
Then $U^2=I$, $UV=VU$, and $(U+V^*)/2$ is a self-adjoint contraction.
Since $(U+I)\mathcal H\subset\ker(U-I)$ and $V$ leaves $\ker(U-I)$
invariant, self-adjointness of $(U+V^*)/2$ gives, relative to this
subspace and its orthogonal complement,
\[
 U=\begin{pmatrix}I&X\\0&-I\end{pmatrix},\qquad
 V=\begin{pmatrix}A&X\\0&-E\end{pmatrix},\qquad
 A=A^*,\quad E=E^*,\quad AX+XE=2X.
\]
The element $\gamma(2p-I)$ is a self-adjoint contraction in the
finite-dimensional algebra $D$. It satisfies a polynomial all of
whose zeros belong to $[-1,1]$. The same polynomial annihilates $V$
and $A$, and its reflected polynomial annihilates $E$.
Thus $A$ and $E$ are self-adjoint contractions with finite spectra.
Let $P_a$ and $Q_e$ be the spectral projections of $A$ and $E$
corresponding to eigenvalues $a$ and $e$. The equality $AX+XE=2X$
gives $(a+e-2)P_aXQ_e=0$. Since $a,e\le1$, the only possible
nonzero term has $a=e=1$. Thus $X$ maps the eigenspace of $E$
at $1$ into the eigenspace of $A$ at $1$ and is zero elsewhere.
On these two spaces the symmetrization has the compression
\[
 \begin{pmatrix}I&X/2\\X^*/2&-I\end{pmatrix}.
\]
Its square has diagonal blocks $I+XX^*/4$ and $I+X^*X/4$.
Contractivity forces $X=0$. Hence $\eta(p)$ is an orthogonal projection.

With normalized Haar measure on the compact unitary group of $D$, put
\[
 B=\int_{\mathcal U(D)}\eta(u)^*\eta(u)\,du.
\]
Since $\eta(u)^{-1}=\eta(u^*)$, one has
$\|\eta\|^{-2}I\le B\le\|\eta\|^2I$.
Haar invariance gives $\eta(u)^*B\eta(u)=B$ for every unitary $u\in D$.
Thus $\pi(d)=B^{1/2}\eta(d)B^{-1/2}$ is a $*$-representation.
Set $L=B^{-1/2}$ and $G=L^2$. Then
\[
 \eta(d)=L\pi(d)L^{-1},\qquad L>0.
\]
For a central projection $p\in D$, the orthogonal projection
$\eta(p)$ commutes with $\eta(u)$ and $\eta(u)^*$ for every unitary
$u\in D$. It therefore commutes with $B$ and $L$, and
$\pi(p)=\eta(p)$. Thus $L$ and $G$ are block diagonal with respect
to the represented simple summands of $D$.
Set $\beta(d)=\pi(\gamma(d^*))$, a positive complex-linear map.
The symmetrization identity becomes
\[
 2L\Xi(d)L=G\pi(d)+\beta(d)G.
\]

Fix a simple summand of $D$ on which $\pi$ is nonzero, with target space
$\mathbb C^r\otimes\mathcal L$, on which $\pi(D)$ acts as
$M_r(\C)\otimes I$. Let $j:M_r(\C)\to D$ be the inclusion of that input summand.
Compression of $\beta j$ has the form $b(\,\cdot\,)\otimes I$,
where $b$ is positive and $b(I)\le I$.
Write $L_r$ and $G_r$ for the corresponding blocks of $L$ and $G$, and
let $\Xi_r(a)$ be the compression of $\Xi(j(a))$ to this target block.
For a state $\omega$ on $B(\mathcal L)$, put
\[
 g=(\mathrm{id}\otimes\omega)(G_r),\qquad
 \Lambda(a)=(\mathrm{id}\otimes\omega)(L_r\Xi_r(a)L_r).
\]
A state is completely positive \cite[Theorem~3]{Stinespring1955}.
Compression, positive congruence, and the state slice preserve positivity,
so $\Lambda$ is positive. Since $G_r$ is
bounded below by a positive scalar multiple of the identity, $g>0$.
Also $0\le j(I)\le I$ and $\Xi(I)=I$, which give $\Lambda(I)\le g$.
The compressed identity becomes
\[
 2\Lambda(a)=ga+b(a)g.
\]
The positive map $b$ and the positive map $\Lambda$
preserve adjoints. Comparing this identity at $a$ and at $a^*$
therefore gives $[g,a-b(a)]=0$.
If $g$ has unequal eigenvalues $\lambda_i,\lambda_j$, take the
matrix unit $e_{ij}$ in an eigenbasis. Then $b(e_{ij})_{ij}=1$.
The map $T(a)=g^{-1/2}\Lambda(a)g^{-1/2}$ is positive and
$0\leq T(I)\leq I$. We claim that $\norm T\leq1$.
For a unitary $u=\sum_k\lambda_k e_k$, put $A_k=T(e_k)$ and
$W\xi=(A_k^{1/2}\xi)_k$. Then $\norm W\leq1$ and
\[
       T(u)=W^*\diag(\lambda_k I)W,
\]
so $\norm{T(u)}\leq1$. Every contraction in $M_r(\C)$ is a convex
combination of unitaries: apply a singular-value decomposition and
express its diagonal contraction as a convex combination of diagonal
matrices with entries in $\{-1,1\}$. This proves the claim.
\par However,
\[
 T(e_{ij})_{ij}
 =\frac12\left(\sqrt{\lambda_i/\lambda_j}
                   +\sqrt{\lambda_j/\lambda_i}\right)>1,
\]
which is impossible. Hence $g$ is scalar for every state $\omega$.
Writing $G_r=[G_{ij}]$ with $G_{ij}\in B(\mathcal L)$, we obtain
$\omega(G_{ij})=0$ for $i\ne j$ and
$\omega(G_{ii}-G_{jj})=0$ for all $i,j$ and all states $\omega$.
States separate operators, so $G_r=I_r\otimes H$ for a positive
invertible $H\in B(\mathcal L)$. This holds on every represented
simple summand, including arbitrary multiplicities. Therefore
$L=G^{1/2}$ commutes with $\pi(D)$ and $\eta=\pi$.
\end{proof}

\begin{proof}[Proof of Theorem~\ref{v5:thm:digraph}]
We may assume $\mathcal H\ne\{0\}$. Let $E_{ij}$ be the standard
matrix units in $M_n(\C)$. For $b=[b_{ij}]\in\Bcal$,
\[
 E_{ii}bE_{jj}=b_{ij}E_{ij}\in\Bcal.
\]
Thus $\Bcal$ is spanned by the matrix units it contains.
Define $i\preceq j$ when $E_{ij}\in\Bcal$.
This relation is reflexive, since every $E_{ii}$ belongs to $\Bcal$,
and transitive, since $E_{ij}E_{jk}=E_{ik}$.
The relation
\[
 i\sim j\quad\Longleftrightarrow\quad
 i\preceq j\ \text{and}\ j\preceq i
\]
is therefore an equivalence relation. Let $C_1,\ldots,C_r$ be its
classes, and put $p_k=\sum_{i\in C_k}E_{ii}$.
These are mutually orthogonal projections in $\Bcal$ with sum $I$.
Each $p_k\Bcal p_k$ is the full matrix algebra on the coordinate
space indexed by $C_k$.

Suppose that both $p_k\Bcal p_l$ and $p_l\Bcal p_k$ are nonzero.
There are $i,i'\in C_k$ and $j,j'\in C_l$ such that
$E_{ij},E_{j'i'}\in\Bcal$. Since
$E_{jj'},E_{i'i}\in\Bcal$, their product gives
\[
 E_{ji}=E_{jj'}E_{j'i'}E_{i'i}\in\Bcal.
\]
Hence $i\sim j$, so $k=l$. Consequently,
\[
 p_k\Bcal p_l=\{0\}\ \text{or}\ p_l\Bcal p_k=\{0\}
 \qquad(k\ne l).
\]
It follows that
\[
 D=\Bcal\cap\Bcal^*
   =\bigoplus_{k=1}^r p_k\Bcal p_k.
\]
The map $E_D(x)=\sum_k p_kxp_k$ is a completely positive projection
from $M_n(\C)$ onto $D$ and preserves the identity.

Write $\Phi=\Phi_{\theta,\alpha}$. Since $\Phi$ is a complete contraction
preserving the identity on $\Bcal\subset M_n(\C)$,
Arveson's theorem \cite[Theorem~1.2.9]{Arveson1969} gives a completely
positive extension $\widetilde\Phi:M_n(\C)\to B(\mathcal H)$ preserving
the identity. In particular, $\Phi|_D$ is completely positive.
The complex-linear map $d\mapsto\alpha(d^*)$ is a contraction
from $D$ to $M_n(\C)$ preserving the identity. It therefore preserves
adjoints by \cite[Proposition~1.2.8]{Arveson1969}. Thus
\[
 \alpha(d)^*=\alpha(d^*)\in\Bcal\qquad(d\in D),
\]
which gives $\alpha(D)\subset D$.
Lemma~\ref{v5:lem:positive_diagonal} now shows that
$\eta=\theta|_D$ is a $*$-homomorphism. Put $P_k=\eta(p_k)$.

By \cite[Theorem~1 and Section~3]{Stinespring1955}, there are a Hilbert
space $\mathcal L$, a $*$-representation $\sigma:M_n(\C)\to B(\mathcal L)$
preserving the identity, and an isometry $V:\mathcal H\to\mathcal L$
such that
\[
 \widetilde\Phi(x)=V^*\sigma(x)V\qquad(x\in M_n(\C)).
\]
The identity $\widetilde\Phi(I)=I$ gives $V^*V=I$.
Integrate over the compact unitary group $\mathcal U(D)$ with
normalized Haar measure and set
\[
 Z=\int_{\mathcal U(D)}\sigma(u)^*V\eta(u)\,du.
\]
The integrand is norm-continuous, so this is an operator-norm integral.
Each integrand has norm one, and therefore $\|Z\|\le1$.
For $v\in\mathcal U(D)$, right invariance of Haar measure gives
\[
 \sigma(v)Z\eta(v)^*
 =\int_{\mathcal U(D)}\sigma(uv^*)^*V\eta(uv^*)\,du=Z.
\]
Since the unitaries linearly span $D$, it follows that
$\sigma(d)Z=Z\eta(d)$ for every $d\in D$.
Set $b_0=Z^*Z$ and define on $M_n(\C)$ the completely positive map
\[
 \Omega_0(x)=Z^*\sigma(x)Z.
\]
Then $0\le b_0\le I$, $b_0$ commutes with $\eta(D)$, and
\[
 \Omega_0(d_1xd_2)=\eta(d_1)\Omega_0(x)\eta(d_2),\qquad
 \Omega_0(d)=\eta(d)b_0
\]
for $x\in M_n(\C)$ and $d,d_1,d_2\in D$.

Let $b\in p_k\Bcal p_l$, where $k\ne l$ and this block is nonzero.
Then $p_l\Bcal p_k=\{0\}$. For $u,v\in\mathcal U(D)$,
\[
 P_k\eta(u)^*\Phi(ubv^*)\eta(v)P_l=\frac12\theta(b).
\]
Indeed, the term containing $\theta(ubv^*)$ gives $\theta(b)/2$,
since $\eta$ preserves adjoints. The other term vanishes because
\[
 P_l\theta(\alpha(ubv^*))P_k
   =\theta\bigl(p_l\alpha(ubv^*)p_k\bigr)=0.
\]
Here we used only $\alpha(ubv^*)\in\Bcal$, not invariance of
the individual blocks under $\alpha$.
Expanding the two norm integrals gives
\[
 \Omega_0(b)=\int_{\mathcal U(D)}\int_{\mathcal U(D)}
       \eta(u)^*\Phi(ubv^*)\eta(v)\,du\,dv.
\]
The bimodule identity gives $\Omega_0(b)=P_k\Omega_0(b)P_l$.
Compressing the last integral to this block therefore yields
\[
 \Omega_0(b)=\frac12\theta(b)
 \qquad(b\in p_k\Bcal p_l,\ k\ne l).
\]

Since $\eta E_D$ is completely positive on $M_n(\C)$, so is
\[
 \Omega(x)=\Omega_0(x)
       +(I-b_0)^{1/2}\eta(E_Dx)(I-b_0)^{1/2}.
\]
For $d\in D$, commutation of $b_0$ with $\eta(D)$ gives
\[
 \Omega(d)=\eta(d)b_0+\eta(d)(I-b_0)=\eta(d).
\]
In particular, $\Omega(I)=I$. On every off-diagonal block of
$\Bcal$, one has $E_D(b)=0$ and $\Omega(b)=\theta(b)/2$.
Summing the blocks gives
\[
 \Omega|_{\Bcal}=\frac12(\theta+\theta E_D).
\]
The left side is completely contractive. Lemma~\ref{v5:lem:positive_correction},
applied to $E_D|_{\Bcal}$, now gives $\norm\theta_{\mathrm{cb}}\le2$.

For sharpness, take $\Bcal=T_2$, $L=\operatorname{diag}(2,1)$,
$\theta(b)=LbL^{-1}$, and
$\alpha(b)=\overline{\operatorname{diag}(b)}$.
The map $\alpha$ is an antilinear complete contraction preserving the
identity, and the symmetrization is the inclusion of $T_2$ in $M_2(\C)$.
Since $\norm{\theta(E_{12})}=2$, the constant is best possible.
\end{proof}
\par

\begin{corollary}\label{cor:triangular-corners}
Let $\Bcal\subset B(\mathcal H_0)$ be a finite-dimensional operator
algebra with unity. Suppose that it contains mutually orthogonal
projections $p_1,\ldots,p_r$ with sum $I$, that each $p_k\Bcal p_k$
is self-adjoint, and that, for $k\ne l$, at least one of
$p_k\Bcal p_l$ and $p_l\Bcal p_k$ is zero. Let
$\theta:\Bcal\to B(\mathcal H)$ be a homomorphism preserving the
identity and let $\alpha:\Bcal\to\Bcal$ be an antilinear contraction
preserving the identity. If
\[
 \Phi(b)=\frac12\bigl(\theta(b)+\theta(\alpha(b))^*\bigr)
\]
is completely contractive, then $\norm\theta_{\mathrm{cb}}\leq2$.
The constant is sharp.
\end{corollary}
\begin{proof}
Put $\mathcal C=C^*(\Bcal)$ and
$D=\bigoplus_kp_k\Bcal p_k$. The absence of reverse off-diagonal
corners gives $D=\Bcal\cap\Bcal^*$. The pinching map
$E(b)=\sum_kp_kbp_k$ is a complete contraction from $\Bcal$ onto
$D$ preserving the identity. A completely positive extension
$\widetilde\Phi:\mathcal C\to B(\mathcal H)$ preserving the identity
exists by Arveson's theorem. The contraction $d\mapsto\alpha(d^*)$
preserves adjoints, so $\alpha(D)\subset D$.
Lemma~\ref{v5:lem:positive_diagonal} gives a $*$-homomorphism
$\eta=\theta|_D$.

Write $\widetilde\Phi(x)=V^*\sigma(x)V$ in Stinespring form and put
\[
 Z=\int_{\mathcal U(D)}\sigma(u)^*V\eta(u)\,du,
 \qquad b_0=Z^*Z,\qquad\Omega_0(x)=Z^*\sigma(x)Z.
\]
Then $\norm Z\leq1$, $\sigma(d)Z=Z\eta(d)$, and
$0\leq b_0\leq I$ belongs to $\eta(D)'$. Consequently,
\[
 \Omega_0(d_1xd_2)=\eta(d_1)\Omega_0(x)\eta(d_2),
 \qquad\Omega_0(d)=\eta(d)b_0.
\]
For an allowed off-diagonal corner $b\in p_k\Bcal p_l$, the
opposite corner is zero. With $P_k=\eta(p_k)$, this implies
\[
 P_k\eta(u)^*\Phi(ubv^*)\eta(v)P_l=\tfrac12\theta(b),
 \qquad u,v\in\mathcal U(D),
\]
because $p_l\alpha(ubv^*)p_k=0$. Expanding the two integrals in
$Z^*\sigma(b)Z$ and using its $(P_k,P_l)$ corner therefore gives
$\Omega_0(b)=\theta(b)/2$.

The map $\eta E:\Bcal\to B(\mathcal H)$ is a complete contraction
preserving the identity. Extend it to a completely positive map
$J:\mathcal C\to B(\mathcal H)$ preserving the identity, and set
\[
 \Omega(x)=\Omega_0(x)+(I-b_0)^{1/2}J(x)(I-b_0)^{1/2}.
\]
This map is completely positive and $\Omega(I)=I$. On $D$ it
equals $\eta$, while on every off-diagonal corner of $\Bcal$ it
equals $\theta/2$. Thus
$\Omega|_{\Bcal}=(\theta+\theta E)/2$ is completely contractive.
Lemma~\ref{v5:lem:positive_correction} gives the conclusion.
The example $\Bcal=T_2$ in Theorem~\ref{v5:thm:digraph} gives sharpness.
\end{proof}
For example, if $W\subset B(\mathcal H_2,\mathcal H_1)$ is a
finite-dimensional operator space, this applies to
\[
 \left\{\begin{pmatrix}aI_{\mathcal H_1}&w\\0&cI_{\mathcal H_2}\end{pmatrix}
       :a,c\in\C,\ w\in W\right\}.
\]
\par

\section{Similarity transfer and scaled \texorpdfstring{$q$}{q}-numerical ranges}\label{sec:q-range}

Let $D\subset\C$ be a bounded convex domain and let
$A\in M_d(\C)$ satisfy $\sigma(A)\subset D$. Define the optimal
complete spectral constant of $D$ for $A$ by
\begin{equation}\label{eq:MDcb}
 M_{D,\mathrm{cb}}(A)
 =\sup_{m\geq1}\sup\left\{\norm{F(A)}:
 F\in M_m(\Acal(\overline D)),\
 \norm{F}_{\overline D}\leq1\right\}.
\end{equation}
The same value is obtained with rectangular matrix-valued functions.

The following proposition applies Theorem~\ref{thm:cb-similarity-amplification} to the functional calculus on a fixed domain.

\begin{proposition}\label{cor:functional-calculus-amplification}
Suppose that $M_{D,\mathrm{cb}}(A)>1$. For every $\gamma\geq1$, there
is a positive invertible matrix $S_\gamma$ such that
\begin{equation}\label{eq:domain-amplification}
 \cond(S_\gamma)=\gamma,
 \qquad
 M_{D,\mathrm{cb}}(S_\gamma^{-1}AS_\gamma)
 =\gamma M_{D,\mathrm{cb}}(A).
\end{equation}
\end{proposition}
\par

\begin{proof}
Let
\[
 \mathfrak A_{A,D}
 =\Acal(\overline D)/\ker\{f\mapsto f(A)\},
\]
with its quotient matrix norms, and let
$\theta_{A,D}:\mathfrak A_{A,D}\to M_d(\C)$ be the induced
homomorphism. The proof of Proposition~\ref{prop:quotient-calculus}
applies without change and gives
\[
        \norm{\theta_{A,D}}_{\mathrm{cb}}
        =M_{D,\mathrm{cb}}(A).
\]
The quotient is finite-dimensional because its image is $\alg(A)$.
Apply Theorem~\ref{thm:cb-similarity-amplification}. Since
\[
 f(S^{-1}AS)=S^{-1}f(A)S,
\]
the amplified homomorphism is the functional calculus of
$S^{-1}AS$ on the same domain $D$. This proves
\eqref{eq:domain-amplification}.
\end{proof}

The following theorem is the completely bounded form of the transfer
principle of Aouichaoui and O'Loughlin
\cite[Theorem~3.3]{AouichaouiOloughlin2026}.

\begin{theorem}\label{thm:complete-transfer}
For each $B\in M_d(\C)$, let $K(B)\subset\C$ be a compact complete
$C$-spectral set for $B$, where $C\geq1$ is independent of $B$. Let
$A\in M_d(\C)$ and let $D$ be a bounded convex domain containing
$\sigma(A)$. Suppose that $\gamma>1$ and
\begin{equation}\label{eq:complete-transfer-hypothesis}
 K(S^{-1}AS)\Subset D
 \qquad\text{whenever }S\text{ is invertible and }\cond(S)<\gamma.
\end{equation}
Then
\begin{equation}\label{eq:complete-transfer-conclusion}
        M_{D,\mathrm{cb}}(A)
        \leq\max\left\{1,\frac C\gamma\right\}.
\end{equation}
\end{theorem}

\begin{proof}
The conclusion is immediate if $M_{D,\mathrm{cb}}(A)=1$. Otherwise,
fix $1<\lambda<\gamma$. Proposition~\ref{cor:functional-calculus-amplification}
gives a positive invertible $S_\lambda$ such that
\[
 \cond(S_\lambda)=\lambda,
 \qquad
 M_{D,\mathrm{cb}}(S_\lambda^{-1}AS_\lambda)
 =\lambda M_{D,\mathrm{cb}}(A).
\]
By \eqref{eq:complete-transfer-hypothesis},
$K(S_\lambda^{-1}AS_\lambda)\Subset D$. Every matrix-valued function
in $\Acal(\overline D)$ is holomorphic on a neighborhood of this
compact set. By Runge's theorem, applied entrywise, the complete
$C$-spectral estimate extends from rational functions to such
functions. Hence
\[
 M_{D,\mathrm{cb}}(S_\lambda^{-1}AS_\lambda)\leq C.
\]
Thus $\lambda M_{D,\mathrm{cb}}(A)\leq C$. Letting
$\lambda\to\gamma^-$ proves \eqref{eq:complete-transfer-conclusion}.
\end{proof}

In the following discussion of $q$-numerical ranges, assume that $d\geq2$.
For $0<|q|\leq1$, define
\[
 W_q(A)=\{\ip{Ax}{y}:\norm{x}=\norm{y}=1,\ \ip{x}{y}=q\},
 \qquad
 \Omega_q(A)=q^{-1}W_q(A).
\]
The set $\Omega_q(A)$ depends only on $|q|$. It is compact and convex,
and it contains the spectrum of $A$; see
\cite{Tsing1984,LiMehtaRodman1994,AouichaouiOloughlin2026}. The
scaled range and its spectral constants were studied in
\cite{OloughlinRani2026}. For $0<r\leq1$, put
\begin{equation}\label{eq:chi-definition}
        \chi(r)=\frac{1+\sqrt{1-r^2}}r.
\end{equation}
Aouichaoui and O'Loughlin proved the exact similarity formula
\begin{equation}\label{eq:q-similarity-formula}
 \Omega_r(A)
 =\bigcup_{\cond(S)\leq\chi(r)}W(S^{-1}AS).
\end{equation}
If $A$ is not scalar, $0<r<1$, and $\cond(S)<\chi(r)$, then
\begin{equation}\label{eq:q-strict-inclusion}
        W(S^{-1}AS)\Subset\Omega_r(A)^\circ
\end{equation}
\cite[Theorem~2.2 and Corollary~2.3]{AouichaouiOloughlin2026}.

For $d\geq2$, let $C^{\mathrm{cb}}_{q,d}$ be the least $C\geq1$ such
that
\begin{equation}\label{eq:Cqd-definition}
 \norm{F(A)}
 \leq C\max_{z\in\Omega_q(A)}\norm{F(z)}
\end{equation}
for every $A\in M_d(\C)$ and every square matrix-valued polynomial
$F$. Since $\Omega_q(A)$ is convex, the same constant is obtained for
matrix-valued functions holomorphic on a neighborhood of
$\Omega_q(A)$.

\begin{theorem}\label{thm:q-complete-constants}
For every $d\geq2$ and $0<|q|\leq1$,
\begin{equation}\label{eq:q-complete-two-sided}
 \max\left\{1,\frac{2}{\chi(|q|)}\right\}
 \leq C^{\mathrm{cb}}_{q,d}
 \leq
 \max\left\{1,\frac{Q_{\mathrm{cb}}(d)}{\chi(|q|)}\right\}.
\end{equation}
In particular,
\begin{equation}\label{eq:q-complete-universal}
 C^{\mathrm{cb}}_{q,d}
 \leq
 \max\left\{1,
 \frac{(1+\sqrt2)|q|}{1+\sqrt{1-|q|^2}}\right\}.
\end{equation}
Consequently, $\Omega_q(A)$ is a complete spectral set for every
matrix $A$ of every order at least two when $|q|\leq1/\sqrt2$. For $d=2,3$,
\begin{equation}\label{eq:q-complete-sharp}
 C^{\mathrm{cb}}_{q,d}
 =\max\left\{1,
 \frac{2|q|}{1+\sqrt{1-|q|^2}}\right\}.
\end{equation}
Thus, in these two dimensions, $\Omega_q(A)$ is a complete spectral
set for every $A$ if and only if $|q|\leq4/5$.
\end{theorem}

\begin{proof}
If $|q|=1$, then $\Omega_q(A)=W(A)$ up to the harmless reduction to
$|q|$, and the upper bound in \eqref{eq:q-complete-two-sided} is the
definition of $Q_{\mathrm{cb}}(d)$. Suppose that $0<|q|<1$. Scalar
matrices have complete spectral constant one on $\Omega_q(A)$. If $A$
is not scalar, apply Theorem~\ref{thm:complete-transfer} with
\[
 K(B)=W(B),
 \qquad C=Q_{\mathrm{cb}}(d),
 \qquad D=\Omega_{|q|}(A)^\circ,
 \qquad \gamma=\chi(|q|).
\]
The hypothesis follows from \eqref{eq:q-strict-inclusion}. This proves
the upper bound.

For the lower bound, let
$N_d=2e_1e_2^*\in M_d(\C)$. The rank-one $q$-numerical-radius formula
\cite[Theorem~2.1]{DencicEtAl2025} gives
\[
        \max_{z\in\Omega_q(N_d)}|z|=\chi(|q|).
\]
Taking the scalar polynomial $p(z)=z$ gives the lower bound
$2/\chi(|q|)$, while constant functions give the lower bound one.
This proves \eqref{eq:q-complete-two-sided}.

The complete Crouzeix--Palencia inequality
\cite{CrouzeixPalencia2017,CrouzeixKressner2025} gives
$Q_{\mathrm{cb}}(d)\leq1+\sqrt2$, which proves
\eqref{eq:q-complete-universal}. Since
$\chi(1/\sqrt2)=1+\sqrt2$, the dimension-free complete spectral-set
assertion follows. Theorem~\ref{thm:intro-complete-three} gives
$Q_{\mathrm{cb}}(2)=Q_{\mathrm{cb}}(3)=2$, so the two bounds coincide
for $d=2,3$. Finally, $\chi(r)\geq2$ exactly when $r\leq4/5$.
\end{proof}

\begin{remark}
Theorems~\ref{thm:cb-similarity-amplification} and
\ref{thm:complete-transfer} hold in every finite dimension. They give
complete estimates on domains larger than the numerical range, but
they do not control the rectangular levels left by
Theorem~\ref{thm:abstract-rectangular-extremal}. In particular, at
$q=1$ Theorem~\ref{thm:q-complete-constants} returns
$Q_{\mathrm{cb}}(d)$ itself. For $d=4$, the unresolved level is the
$2\times2$ matrix level.
\end{remark}

\section{Equality and quantitative estimates}\label{sec:equality}
At $\tau=1/2$, both square-function estimates are needed to identify
the reducing subspace associated with the singular value $2$. The
results in this section are first stated for bounded operators. For an
operator $T$ satisfying these estimates, put
\[
        E_T=\ker(4I-T^*T),
        \qquad
        F_T=\frac12TE_T.
\]

\begin{theorem}\label{thm:equality-block}
Let $T\in\Bcal(H)$ satisfy
\[
 \sum_{k=0}^{\infty}4^{-k}T^{*k}T^k\preceq2I,
 \qquad
 \sum_{k=0}^{\infty}4^{-k}T^kT^{*k}\preceq2I.
\]
Then $E_T\perp F_T$, the subspace $E_T\oplus F_T$ reduces $T$, and,
after identifying $E_T$ with $F_T$ by $x\mapsto Tx/2$,
\[
 T\simeq_u
 \begin{pmatrix}0&0\\2I_{E_T}&0\end{pmatrix}\oplus R.
\]
If $\norm{T}=2$ is attained, then $E_T\neq\{0\}$. If $H$ is
finite-dimensional and $\norm{T}=2$, then $\norm{R}<2$, and $T$ is not
diagonalizable.
\end{theorem}

\begin{proof}
Put $E=E_T$. The map $U=T|_E/2$ is an isometry from $E$ onto
$F=F_T$. Let $x\in E$. The first square-function estimate gives
\[
\begin{aligned}
2\norm{x}^2
&\geq\sum_{k=0}^{\infty}4^{-k}\norm{T^kx}^2\\
&=2\norm{x}^2+
  \sum_{k=2}^{\infty}4^{-k}\norm{T^kx}^2.
\end{aligned}
\]
Thus $T^2x=0$. If $y=Ux$, then
\[
        Ty=0,
        \qquad
        T^*y=2x.
\]
The second square-function estimate applied to $y$ gives
$(T^*)^2y=0$, and therefore $T^*x=0$.

If $x,z\in E$, then
\[
        \ip{x}{Uz}=\frac12\ip{T^*x}{z}=0.
\]
Hence $E\perp F$, and
\[
        TE=2F,
        \qquad TF=0,
        \qquad T^*F=2E,
        \qquad T^*E=0.
\]
Thus $E\oplus F$ reduces $T$ and gives the stated block.

If $\norm{T}=2$ is attained at a unit vector $x$, then
$T^*Tx=4x$, so $x\in E$. Suppose now that $H$ is finite-dimensional and
$\norm{T}=2$. Then $E\neq\{0\}$. If the restriction $R$ to
$(E\oplus F)^\perp$ had norm $2$, a unit vector at which its norm is
attained would belong to $E\cap(E\oplus F)^\perp$, a contradiction.
Hence $\norm{R}<2$. The nonzero nilpotent summand shows that $T$ is not
diagonalizable.
\end{proof}

\begin{remark}\label{rem:similar-normal-proof2}
As a consequence, we obtain another proof that $\psi(A)<2$ if $A=V\Lambda V^{-1}$ is diagonalizable. Indeed, $\psi(A)$ is attained, say at a scalar-valued function $f$ with $\|f\|_{W(A)}=1$. If $\psi(A)=2$, then $T=f(A)$ satisfies $\|T\|=2$, and Theorem~\ref{thm:tuple-square} with $\tau=1/2$ allows us to apply the preceding theorem to conclude that $T$ is not diagonalizable. On the other hand,
$T=f(A)=Vf(\Lambda)V^{-1}$, so $T$ is diagonalizable, which is a contradiction.
\end{remark}

When the norm is not attained, $E_T$ may be trivial. The next theorem
replaces it by two-dimensional subspaces on which $T$ approaches the
nilpotent equality block.

\begin{theorem}\label{thm:approximate-equality}
Let $T\in\Bcal(H)$ satisfy the two square-function estimates in
Theorem~\ref{thm:equality-block}, and assume that $\norm{T}=2$. Then
there are orthogonal decompositions
\[
        H=\mathcal M_j\oplus\mathcal M_j^\perp,
        \qquad
        \dim\mathcal M_j=2,
\]
and operators $D_j\in\Bcal(\mathcal M_j^\perp)$ such that, in a
suitable orthonormal basis of $\mathcal M_j$,
\[
 \norm{T-\left(
 \begin{pmatrix}0&0\\2&0\end{pmatrix}
 \oplus D_j\right)}\longrightarrow0.
\]
\end{theorem}

\begin{proof}
Since $4=\norm{T^*T}$ belongs to the spectrum of the positive operator
$T^*T$, there are unit vectors $x_j$ such that
\[
        \norm{(T^*T-4I)x_j}\longrightarrow0.
\]
After discarding finitely many indices, put
\[
        s_j=\norm{Tx_j}>0,
        \qquad
        y_j=\frac{Tx_j}{s_j}.
\]
Then $s_j\to2$ and
\[
        \norm{T^*y_j-2x_j}\longrightarrow0.
\]
The first square-function estimate gives
$\norm{T^2x_j}\to0$, and hence $\norm{Ty_j}\to0$. Since
$\norm{T^*y_j}\to2$, the second square-function estimate gives
$\norm{(T^*)^2y_j}\to0$. Indeed,
\[
 2T^*x_j=(T^*)^2y_j-T^*(T^*y_j-2x_j)\longrightarrow0,
\]
and
\[
 \abs{\ip{y_j}{x_j}}
 =\frac{\abs{\ip{Tx_j}{x_j}}}{s_j}
 =\frac{\abs{\ip{x_j}{T^*x_j}}}{s_j}
 \longrightarrow0.
\]

Write
\[
 \beta_j=\ip{y_j}{x_j},
 \qquad
 \rho_j=(1-\abs{\beta_j}^2)^{1/2},
 \qquad
 e_j=\frac{y_j-\beta_jx_j}{\rho_j},
\]
after discarding finitely many further indices. Put
$\mathcal M_j=\operatorname{span}\{x_j,e_j\}$. Then
$\{x_j,e_j\}$ is orthonormal and
\[
 Tx_j-2e_j\longrightarrow0,
 \qquad
 Te_j\longrightarrow0,
 \qquad
 T^*x_j\longrightarrow0,
 \qquad
 T^*e_j-2x_j\longrightarrow0.
\]
Let $P_j$ be the orthogonal projection onto $\mathcal M_j$, let
$N_jx_j=2e_j$ and $N_je_j=0$, and set
\[
        D_j=(I-P_j)T|_{\mathcal M_j^\perp}.
\]
The four preceding limits give
\[
 \norm{P_jTP_j-N_j}\longrightarrow0,
 \qquad
 \norm{(I-P_j)TP_j}\longrightarrow0,
 \qquad
 \norm{P_jT(I-P_j)}\longrightarrow0.
\]
These are the nonzero blocks of $T-(N_j\oplus D_j)$, and the result
follows.
\end{proof}

We next solve the optimization problem used in the quantitative
estimates. For $0<s\leq2$, let $\tau_s$ be the unique
number in $[1/2,3/4)$ satisfying
\begin{equation}\label{eq:tau-s}
        2s^2\tau_s(1-\tau_s)^2+4\tau_s-3=0,
\end{equation}
and put
\[
        \eta(s)
        =\left(
        \frac{2\tau_s-1}
        {2\tau_s^3(1-\tau_s)^2}
        \right)^{1/2}.
\]
Thus $\eta(2)=0$. The functions $L_\tau$ and $L$ were defined in
\eqref{eq:L-tau-definition}.

\begin{lemma}\label{lem:tau-optimization}
For $0<s\leq2$,
\begin{equation}\label{eq:eta-minimum}
 \eta(s)
 =\min_{0<\tau<1}
 \frac1{\tau^2}
 \left(
 \frac{\tau}{1-\tau}-s^2\tau^2
 \right)^{1/2}.
\end{equation}
Moreover,
\begin{equation}\label{eq:eta-simple-bound}
        \eta(s)\leq2\sqrt{4-s^2},
\end{equation}
with strict inequality when $s<2$. For $0<s<2$,
\begin{equation}\label{eq:L-maximum}
        L(s)=\max_{0<\tau<1}L_\tau(s)
        >L_{1/2}(s)=\frac{s^4}{4(4-s^2)}.
\end{equation}
The functions $s\mapsto\tau_s$ and $s\mapsto\eta(s)$ are continuous on
$(0,2]$, and $L$ is continuous on $(0,2)$.
\end{lemma}

\begin{proof}
Put
\[
        q_s(\tau)=2s^2\tau(1-\tau)^2+4\tau-3.
\]
Its derivative satisfies
\[
        q_s'(\tau)
        =2s^2(1-\tau)(1-3\tau)+4
        \geq4-\frac{2s^2}{3}>0.
\]
Also, $q_s(0)=-3$, $q_s(1)=1$,
$q_s(1/2)=s^2/4-1\leq0$, and $q_s(3/4)=3s^2/32>0$.
Thus \eqref{eq:tau-s} has the asserted unique solution.

The square of the expression minimized in \eqref{eq:eta-minimum} is
\[
 d_s(\tau)
 =\frac1{\tau^3(1-\tau)}-\frac{s^2}{\tau^2}.
\]
A direct calculation gives
\[
        d_s'(\tau)
        =\frac{q_s(\tau)}{\tau^4(1-\tau)^2}.
\]
Since $d_s(\tau)$ tends to infinity at both endpoints, its unique
minimum occurs at $\tau_s$. Substitution of \eqref{eq:tau-s} gives
\[
        d_s(\tau_s)
        =\frac{2\tau_s-1}
        {2\tau_s^3(1-\tau_s)^2}.
\]
This proves \eqref{eq:eta-minimum}. Choosing $\tau=1/2$ gives
\eqref{eq:eta-simple-bound}. The inequality is strict for $s<2$ because
then $\tau_s>1/2$.

For $0<s<2$, the denominator in $L_\tau(s)$ is positive. Differentiation
gives
\[
 \frac{\partial}{\partial\tau}L_\tau(s)
 =-\frac{s^4\tau^2q_s(\tau)}
 {\{1-s^2\tau(1-\tau)\}^2}.
\]
Hence the unique maximum occurs at $\tau_s$. Since $\tau_s>1/2$,
\eqref{eq:L-maximum} follows. Continuity follows from the implicit
function theorem applied to \eqref{eq:tau-s}.
\end{proof}

Put
\[
        \vartheta(s)=\min\{s^2,\eta(s)\},
        \qquad 0<s\leq2.
\]

\begin{theorem}\label{thm:tuple-second-word}
Let $T_1,\ldots,T_q\in\Bcal(H)$ satisfy, for every $0<\tau<1$,
\begin{equation}\label{eq:tuple-second-assumption}
 \sum_{k=0}^{\infty}\tau^{2k}
 \sum_{\alpha\in\mathcal W_k}T_\alpha^*T_\alpha
 \preceq\frac1{1-\tau}I.
\end{equation}
Let
\[
        s=\norm{(T_1,\ldots,T_q)^{\mathsf T}}>0,
\]
and suppose that this norm is attained at a unit vector $x$. Then
$s\leq2$ and
\begin{equation}\label{eq:tuple-second-word-bound}
 \left(\sum_{i,j=1}^q\norm{T_iT_jx}^2\right)^{1/2}
 \leq\vartheta(s).
\end{equation}
The function $\vartheta$ cannot be replaced by a smaller function of
$s$, since the scalar case $q=1$ is sharp.
\end{theorem}

\begin{proof}
Keeping the terms of degrees zero and one in
\eqref{eq:tuple-second-assumption}, and taking $\tau=1/2$, gives
$s\leq2$. Applying the quadratic form at $x$ and retaining the terms
through degree two gives
\[
 \tau^4\sum_{i,j=1}^q\norm{T_iT_jx}^2
 \leq\frac{\tau}{1-\tau}-\tau^2s^2.
\]
Lemma~\ref{lem:tau-optimization} gives the bound by $\eta(s)$. On the
other hand, the definition of the column norm gives
\[
 \sum_{i,j=1}^q\norm{T_iT_jx}^2
 \leq s^2\sum_{j=1}^q\norm{T_jx}^2=s^4.
\]
Taking the smaller of these two bounds proves
\eqref{eq:tuple-second-word-bound}. The sharp scalar examples are given
in Theorem~\ref{thm:sharp-T2}.
\end{proof}

The same optimization gives a lower bound on every eigenvector
condition number that realizes a diagonalization in dimension three.

\begin{corollary}\label{cor:complete-condition-number}
Let $A\in M_3(\C)$ be diagonalizable, write
$A=V\Lambda V^{-1}$, and put
\[
        C=\psi_{\mathrm{cb}}(A),
        \qquad \kappa=\cond(V).
\]
Then
\begin{equation}\label{eq:complete-condition-lower}
        \kappa^4\geq L(C)
        >\frac{C^4}{4(4-C^2)}.
\end{equation}
In particular, if $A_j=V_j\Lambda_jV_j^{-1}\in M_3(\C)$ are
diagonalizable and $\psi_{\mathrm{cb}}(A_j)\to2$, then
$\cond(V_j)\to\infty$ for every choice of eigenvector matrices $V_j$.
\end{corollary}

\begin{proof}
If $C=1$, apply Theorem~\ref{thm:quantitative-tuple-normal} to
$I=VIV^{-1}$. We may therefore assume that $C>1$.
Choose a minimal scalar, column, or row extremal as in the proof of
Theorem~\ref{thm:intro-complete-diagonalizable}. Apply
\eqref{eq:tuple-condition-lower} to that scalar, column, or adjoint row
tuple. This gives \eqref{eq:complete-condition-lower}. The final
assertion follows.
\end{proof}

\begin{theorem}\label{thm:sharp-T2}
Let $T\in\Bcal(H)$ and assume that, for every $0<\tau<1$,
\[
 \sum_{k=0}^{\infty}\tau^{2k}T^{*k}T^k
 \preceq\frac1{1-\tau}I,
 \qquad
 \sum_{k=0}^{\infty}\tau^{2k}T^kT^{*k}
 \preceq\frac1{1-\tau}I.
\]
Assume that $s=\norm{T}>0$ is attained, and choose unit singular vectors
$x,y$ such that
\[
        Tx=sy,
        \qquad
        T^*y=sx.
\]
Then $s\leq2$ and
\[
        \norm{T^2x},\ \norm{(T^*)^2y}\leq\vartheta(s).
\]
The estimate for $\norm{T^2x}$ is sharp for every $0<s\leq2$, already
among three-dimensional nilpotent weighted shifts. Taking adjoints gives
sharpness for $\norm{(T^*)^2y}$. More precisely,
\[
 \vartheta(s)=
 \begin{cases}
 s^2,&0<s\leq3/2,\\
 \eta(s),&3/2\leq s\leq2.
 \end{cases}
\]
\end{theorem}

\begin{proof}
Theorem~\ref{thm:tuple-second-word}, applied with $q=1$, gives
$s\leq2$ and $\norm{T^2x}\leq\vartheta(s)$. Applying the same theorem
to $T^*$ and $y$ gives
$\norm{(T^*)^2y}\leq\vartheta(s)$.

For $s>0$ and $a\geq0$, let
\[
 W_{s,a}=
 \begin{pmatrix}
 0&0&0\\
 s&0&0\\
 0&a&0
 \end{pmatrix}.
\]
Let $e_1,e_2,e_3$ be the standard basis vectors of $\C^3$. Then
$W_{s,a}e_1=se_2$, $W_{s,a}e_2=ae_3$, and
\[
\begin{aligned}
 \sum_{k=0}^{\infty}\tau^{2k}W_{s,a}^{*k}W_{s,a}^k
 &=\diag\bigl(1+s^2\tau^2+s^2a^2\tau^4,
                1+a^2\tau^2,1\bigr),\\
 \sum_{k=0}^{\infty}\tau^{2k}W_{s,a}^kW_{s,a}^{*k}
 &=\diag\bigl(1,1+s^2\tau^2,
                1+a^2\tau^2+s^2a^2\tau^4\bigr).
\end{aligned}
\]

Suppose first that $0<s\leq3/2$ and take $a=s$. Put $b=s^2$ and
\[
        R_b(\tau)=1-b\tau+b\tau^2-b^2\tau^3+b^2\tau^4.
\]
Since
\[
 \frac1{1-\tau}-(1+b\tau^2+b^2\tau^4)
 =\frac{\tau R_b(\tau)}{1-\tau}
\]
and
\[
 \frac{\partial R_b(\tau)}{\partial b}
 =-\tau(1-\tau)(1+2b\tau^2)\leq0,
\]
one has
\[
 R_b(\tau)\geq R_{9/4}(\tau)
 =\frac{(3\tau-2)^2(9\tau^2+3\tau+4)}{16}\geq0.
\]
Thus $W_{s,s}$ satisfies both square-function estimates. Its norm is
$s$ and $\norm{W_{s,s}^2e_1}=s^2$.

Now suppose that $3/2\leq s\leq2$. Evaluating
\eqref{eq:eta-minimum} at $\tau=2/3$ gives
\[
        \eta(s)^2\leq\frac{81}{8}-\frac94s^2\leq s^4.
\]
Set $a=\eta(s)/s$. Then $0\leq a\leq s$, and the definition of
$\eta(s)$ gives
\[
        1+s^2\tau^2+\eta(s)^2\tau^4
        \leq\frac1{1-\tau}
        \qquad(0<\tau<1).
\]
Each remaining diagonal entry of either square function is bounded by
$1+s^2\tau^2+\eta(s)^2\tau^4$ or by $1+s^2\tau^2$. Moreover,
$1+s^2\tau^2\leq(1-\tau)^{-1}$ because $s\leq2$. Hence
$W_{s,a}$ satisfies both estimates, $\norm{W_{s,a}}=s$, and
\[
        \norm{W_{s,a}^2e_1}=sa=\eta(s).
\]
This proves sharpness and the displayed formula for $\vartheta$.
\end{proof}

We next obtain an operator-norm estimate. More precisely, we choose a
two-dimensional subspace $\mathcal M$ and compare $T$ with
$\left(\begin{smallmatrix}0&0\\2&0\end{smallmatrix}\right)\oplus D$
relative to the decomposition
$H=\mathcal M\oplus\mathcal M^\perp$.

\begin{theorem}\label{thm:operator-stability}
Let $T\in\Bcal(H)$ satisfy
\[
 \sum_{k=0}^{\infty}4^{-k}T^{*k}T^k\preceq2I,
 \qquad
 \sum_{k=0}^{\infty}4^{-k}T^kT^{*k}\preceq2I.
\]
Put
\[
        s=\norm{T},
        \qquad
        \delta=1-\frac{s^2}{4}.
\]
If $0\leq\delta\leq1/16$, and either $\delta>0$ or $\norm T$ is attained, then there are an orthogonal decomposition
$H=\mathcal M\oplus\mathcal M^\perp$, with $\dim\mathcal M=2$, and an
operator $D\in\Bcal(\mathcal M^\perp)$ such that
\[
 \norm{T-\left(
 \begin{pmatrix}0&0\\2&0\end{pmatrix}
 \oplus D\right)}
 \leq8\sqrt\delta.
\]
The displayed $2\times2$ matrix is written in a suitable orthonormal
basis of $\mathcal M$.
\end{theorem}

\begin{proof}
First suppose that $\norm T$ is attained.
Choose unit singular vectors $x,y$ such that
\[
        Tx=sy,
        \qquad
        T^*y=sx.
\]
Since $s\geq\sqrt{15}/2>\sqrt2$, the vectors $x$ and $y$ are linearly
independent. Indeed, if $y$ were a constant of modulus one times $x$,
then the first square-function series on $x$ would be a geometric series
with ratio $s^2/4>1/2$, and its sum would be larger than $2$. Put
\[
 \beta=\ip{y}{x},
 \qquad
 \rho=(1-\abs{\beta}^2)^{1/2},
 \qquad
 e=\frac{y-\beta x}{\rho},
 \qquad
 \mathcal M=\operatorname{span}\{x,e\}.
\]
Then $\{x,e\}$ is an orthonormal basis of $\mathcal M$ and
$y=\beta x+\rho e$. Let $P$ be the orthogonal projection onto
$\mathcal M$. Define $N$ on $\mathcal M$ by
\[
        Nx=2e,
        \qquad
        Ne=0,
\]
and put $D=(I-P)T|_{\mathcal M^\perp}$.

First,
\[
        (PTP-N)x=s\beta x+(s\rho-2)e,
\]
while
\[
        \norm{(PTP-N)e}
        \leq\frac{\norm{Ty}+s\abs{\beta}}{\rho}.
\]
Using the Hilbert--Schmidt norm of $(PTP-N)|_{\mathcal M}$ gives
\begin{equation}\label{eq:stability-compression}
 \norm{PTP-N}
 \leq
 \left(
 s^2\abs{\beta}^2+\abs{s\rho-2}^2
 +\frac{(\norm{Ty}+s\abs{\beta})^2}{\rho^2}
 \right)^{1/2}.
\end{equation}
Moreover,
\begin{equation}\label{eq:stability-lower-block}
        \norm{(I-P)TP}
        \leq\frac{\norm{Ty}}{\rho}.
\end{equation}
Indeed, $(I-P)TPx=0$, and the estimate on $e$ follows from
$e=(y-\beta x)/\rho$. If $w=(I-P)T^*x$, then
\[
 (I-P)T^*Px=w,
 \qquad
 (I-P)T^*Pe=-\frac{\beta}{\rho}w,
\]
because $T^*y=sx$. Hence
\begin{equation}\label{eq:stability-upper-block}
        \norm{PT(I-P)}
        \leq\frac{\norm{T^*x}}{\rho}.
\end{equation}
The norm of $PT(I-P)+(I-P)TP$ is at most the larger of the two bounds in
\eqref{eq:stability-lower-block} and \eqref{eq:stability-upper-block}.

The two square-function estimates give
\[
        \norm{T^2x},\ \norm{(T^*)^2y}\leq4\sqrt\delta.
\]
Since $T^2x=sTy$, $(T^*)^2y=sT^*x$, and
$s^2=4(1-\delta)$, one obtains
\[
 \norm{Ty},\ \norm{T^*x}\leq\frac{11}{5}\sqrt\delta,
 \qquad
 \abs{\beta}\leq\frac{\sqrt\delta}{1-\delta}\leq\frac4{15}.
\]
Thus $\rho\geq15/16$, $s\abs{\beta}\leq(11/5)\sqrt\delta$, and
\[
        2-s\rho
        \leq2\delta+2\abs{\beta}^2
        \leq\frac{11}{10}\sqrt\delta.
\]
It follows from \eqref{eq:stability-compression} that
\[
        \norm{PTP-N}\leq\frac{28}{5}\sqrt\delta,
\]
and $PT(I-P)+(I-P)TP$ has norm at most $(12/5)\sqrt\delta$.
The triangle inequality gives the stated bound.
\par
Now suppose that $0<\delta\leq1/16$, without assuming norm
attainment. Choose unit vectors $x_j$ with
$s_j=\norm{Tx_j}\to s$, and put $y_j=Tx_j/s_j$.
Since $s^2I-T^*T\succeq0$,
$\norm{(T^*T-s^2I)x_j}\to0$, whence
\[
 r_j=T^*y_j-sx_j\longrightarrow0.
\]
Applying the column estimate through degree two to $x_j$ gives
$\limsup_j\norm{T^2x_j}\leq4\sqrt\delta$. Applying the row
estimate to $y_j$ and using $T^*y_j-sx_j\to0$ gives the analogous
bound for $(T^*)^2y_j$. It follows that
\[
 \limsup_j\norm{Ty_j},\quad\limsup_j\norm{T^*x_j}
       \leq\frac{4\sqrt\delta}{s},
 \qquad
 \limsup_j|\ip{y_j}{x_j}|\leq\frac{\sqrt\delta}{1-\delta}.
\]
Define $\beta_j,\rho_j,e_j,\mathcal M_j,P_j,N_j,D_j$ as above.
These are well-defined for large $j$, since
$\limsup_j|\beta_j|\leq4/15<1$. The compression and lower-block
calculations apply with $s_j$ in place of $s$. For the upper block,
\[
 (I-P_j)T^*e_j
 =\frac{-\beta_j(I-P_j)T^*x_j+(I-P_j)r_j}{\rho_j},
\]
so its norm is at most $\norm{T^*x_j}/\rho_j+o(1)$.
Using $4/s\leq11/5$, $\liminf_j\rho_j\geq15/16$, and the same
compression estimate, we obtain
\[
 \limsup_j\norm{T-(N_j\oplus D_j)}
 \leq\frac{\sqrt{631741}+352}{150}\sqrt\delta
 <8\sqrt\delta.
\]
The strict inequality follows from $631741<848^2$.
Choosing a sufficiently large $j$ proves the assertion.
\par
\end{proof}

\begin{corollary}\label{cor:stability-Crouzeix}
Let $A\in M_n(\C)$ and let $f\in\Acal(W(A))$ satisfy
$\norm{f}_{W(A)}=1$. If
\[
        \norm{f(A)}\geq2-\varepsilon,
        \qquad
        0\leq\varepsilon\leq2-\frac{\sqrt{15}}2,
\]
then there are an orthogonal decomposition
$\C^n=\mathcal M\oplus\mathcal M^\perp$, with
$\dim\mathcal M=2$, and a matrix $D$ on $\mathcal M^\perp$ such that
\[
 \norm{f(A)-\left(
 \begin{pmatrix}0&0\\2&0\end{pmatrix}
 \oplus D\right)}
 \leq8\sqrt\varepsilon.
\]
The displayed block is written in a suitable orthonormal basis of
$\mathcal M$.
\end{corollary}

\begin{proof}
Put $T=f(A)$. The bound on $\varepsilon$ gives
$\norm{T}\geq\sqrt{15}/2$, so Theorem~\ref{thm:operator-stability}
applies. Since
\[
        1-\frac{\norm{T}^2}{4}
        \leq2-\norm{T}\leq\varepsilon,
\]
the result follows.
\end{proof}

\section{Representing measures for extremal functions}\label{sec:measures}
We now apply the square-function estimates to extremal functions on a
fixed bounded Jordan domain. Let $\Omega$ be such a domain, let
$A\in M_n(\C)$, and assume
\[
        W(A)\subset\overline\Omega,
        \qquad
        \sigma(A)\subset\Omega.
\]
We write $\Acal(\overline\Omega)$ for the functions that are continuous
on $\overline\Omega$ and holomorphic in $\Omega$. For
$h\in H^\infty(\Omega)$, the matrix $h(A)$ is defined by
the holomorphic functional calculus. Put
\[
 \psi_\Omega(A)
 =\sup\{\norm{h(A)}:h\in H^\infty(\Omega),\
                    \norm{h}_{\infty,\Omega}\leq1\}.
\]
We call $f$ a normalized extremal function for $(A,\Omega)$ if
\[
        \norm{f}_{\infty,\Omega}=1,
        \qquad
        \norm{f(A)}=\psi_\Omega(A).
\]

\begin{proposition}\label{prop:domain-extremals}
Under the preceding assumptions, $\psi_\Omega(A)\leq2$. The supremum is
attained, and every normalized extremal function has the form
\[
        f=B\circ\varphi,
\]
where $\varphi:\Omega\to\D$ is a conformal bijection and $B$ is a
finite Blaschke product, that is, a unimodular constant
times a possibly empty finite product of automorphisms of $\D$,
of degree at most $d_A-1$, where $d_A=\deg m_A$. In particular,
\begin{equation}\label{eq:boundary-inner}
        \abs{f}=1
        \qquad\text{on }\partial\Omega.
\end{equation}
If $\Omega=W(A)^\circ$, then
\[
        \psi_\Omega(A)=\psi(A).
\]
\end{proposition}

\begin{proof}
Let $\varphi:\Omega\to\D$ be a conformal bijection. Since
$\partial\Omega$ is a Jordan curve, $\varphi$ and $\varphi^{-1}$ extend
continuously to the closures \cite[Chapter~2]{Duren1983}. For $h\in H^\infty(\Omega)$ and
$0<r<1$, put
\[
        h_r(z)=h\bigl(\varphi^{-1}(r\varphi(z))\bigr).
\]
Then $h_r\in\Acal(\overline\Omega)$,
$\norm{h_r}_{\overline\Omega}\leq\norm{h}_{\infty,\Omega}$, and
$h_r(A)\to h(A)$ as $r\to1^-$. By Mergelyan's theorem, $h_r$ can be
approximated uniformly on $\overline\Omega$ by polynomials. Since
$W(A)\subset\overline\Omega$, Theorem~\ref{thm:tuple-square} gives
\[
        \norm{h_r(A)}
        \leq2\norm{h}_{\infty,\Omega}.
\]
Letting $r\to1^-$ proves $\psi_\Omega(A)\leq2$.

To apply Crouzeix's extremal theorem on a disk, put $T=\varphi(A)$.
For $h\in H^\infty(\Omega)$ and $g=h\circ\varphi^{-1}$, the
composition rule gives $h(A)=g(T)$, and their supremum norms agree.
Thus $\psi_\Omega(A)=\psi_\D(T)$, and normalized extremals correspond
in both directions. The disk satisfies the hypotheses of the cited
theorem. The same transfer applies to every invariant restriction of $A$.

Crouzeix's finite-dimensional extremal theorem gives attainment
\cite[Theorem~2.1]{Crouzeix2004}. Let $f$ be any normalized extremal
and choose a unit vector $x$ with $\norm{f(A)x}=\psi_\Omega(A)$.
The cyclic subspace $E=\operatorname{span}\{A^jx:j\geq0\}$ is
invariant under $A$ and has dimension at most $d_A$. Holomorphic
functional calculus respects this restriction. Hence
$\psi_\Omega(A|_E)\leq\psi_\Omega(A)$, while the vector $x$ gives
the reverse inequality for $f$. Thus $f$ is also a normalized
extremal for $A|_E$. Applying the same theorem to $A|_E$ gives
the stated Blaschke form with degree at most $\dim E-1\leq d_A-1$.
\par If
$\Omega=W(A)^\circ$, polynomials give
$\psi_\Omega(A)\geq\psi(A)$. Repeating the radial approximation above
and using the definition of $\psi(A)$ in place of the constant $2$
gives the reverse inequality.
\end{proof}

Put
\[
        K=\psi_\Omega(A),
\]
and assume $K>1$. Choose a normalized extremal function $f$ as in
Proposition~\ref{prop:domain-extremals} and a unit vector $x$ such that
\[
        \norm{f(A)x}=K.
\]
Write
\[
        F=f(A),
        \qquad
        y=\frac1KFx.
\]
Then
\[
        Fx=Ky,
        \qquad
        F^*y=Kx.
\]
Apply Lemma~\ref{lem:potapov} at the scalar level to the
homomorphism $h\mapsto h(A)$ and to the extremal function $f$. This
gives
\begin{equation}\label{eq:extremal-orthogonality}
        \ip{Fx}{x}=0,
        \qquad
        x\perp y.
\end{equation}
This is the usual extremal orthogonality relation. See also
\cite[Section~4.1]{BickelEtAl2020}.
The proof of \cite[Theorem~4.5]{BickelEtAl2020} also
applies to the present Jordan domain, without smoothness of the boundary.
We include it because the uniqueness statement is used below.

\begin{lemma}\label{lem:representing-measure}
There is a unique probability measure $\mu$ on $\partial\Omega$ such
that
\begin{equation}\label{eq:representing-measure}
        \ip{h(A)x}{x}
        =\int_{\partial\Omega}h\dd\mu
        \qquad(h\in\Acal(\overline\Omega)).
\end{equation}
Moreover,
\begin{equation}\label{eq:mean-zero}
        \int_{\partial\Omega}f\dd\mu=0.
\end{equation}
\end{lemma}

\begin{proof}
Define
\[
        \ell(h)=\ip{h(A)x}{x}
        \qquad(h\in\Acal(\overline\Omega)).
\]
Let $h\in\Acal(\overline\Omega)$ and suppose that $\Ree h\geq0$ on
$\Omega$. For $t\geq0$,
\[
        \norm{e^{-th}f}_{\overline\Omega}\leq1.
\]
Extremality therefore gives
\[
        \norm{e^{-t h(A)}Fx}\leq K.
\]
The function
$t\mapsto\norm{e^{-t h(A)}Fx}^2$ has a maximum at the endpoint $t=0$
of $[0,\infty)$. Since $F$ commutes with $h(A)$ and $F^*Fx=K^2x$, its
right derivative at $0$ is nonpositive. This gives
\[
        \Ree\ell(h)\geq0.
\]

Let $\mathcal R$ be the real parts on $\partial\Omega$ of functions in
$\Acal(\overline\Omega)$, and define
\[
        \mathcal L(\Ree h)=\Ree\ell(h).
\]
This is well defined. If $\Ree h=0$ on $\partial\Omega$, then $h$ is a
purely imaginary constant. If $\Ree h\geq0$ on $\partial\Omega$, the minimum principle
gives
$\Ree h\geq0$ in $\Omega$. The preceding argument and
$\mathcal L(1)=1$ therefore show that $\mathcal L$ is a positive
functional of norm one on $\mathcal R$. Let
$\varphi:\Omega\to\D$ be a conformal bijection, extended
continuously to the closures. The real and imaginary parts of
$\varphi^k$, restricted to $\partial\Omega$, belong to $\mathcal R$
and, when $\varphi(\zeta)=e^{it}$, equal $\cos(kt)$ and
$\sin(kt)$. Since
the real trigonometric polynomials are uniformly dense in
$C(\T,\mathbb R)$, the space $\mathcal R$ is uniformly dense in
$C(\partial\Omega,\mathbb R)$. Hence $\mathcal L$ extends uniquely
to a positive functional of norm one on that space. The Riesz representation theorem
gives a unique probability measure $\mu$ such that
\[
        \Ree\ell(h)=\int_{\partial\Omega}\Ree h\dd\mu.
\]
Applying this identity to $-ih$ also gives the imaginary parts and
proves \eqref{eq:representing-measure}. Finally,
\eqref{eq:extremal-orthogonality} gives \eqref{eq:mean-zero}.
\end{proof}

The following two identities are the only matrix coefficients needed
below.

\begin{proposition}\label{prop:measure-identities}
For every $h\in\Acal(\overline\Omega)$,
\begin{align}
        \ip{h(A)y}{y}
        &=\int_{\partial\Omega}h\dd\mu,
        \label{eq:same-state}\\
        K\ip{h(A)y}{x}
        &=\int_{\partial\Omega}hf\dd\mu.
        \label{eq:cross-state}
\end{align}
\end{proposition}

\begin{proof}
Since $F$ commutes with $h(A)$,
\[
\begin{aligned}
        K^2\ip{h(A)y}{y}
        &=\ip{h(A)Fx}{Fx}\\
        &=\ip{h(A)x}{F^*Fx}\\
        &=K^2\ip{h(A)x}{x}.
\end{aligned}
\]
This proves \eqref{eq:same-state}. Using $F^*y=Kx$ and then
\eqref{eq:same-state},
\[
\begin{aligned}
        K\ip{h(A)y}{x}
        &=\ip{h(A)y}{F^*y}\\
        &=\ip{Fh(A)y}{y}\\
        &=\int_{\partial\Omega}hf\dd\mu.
\end{aligned}
\]
\end{proof}

Put
\[
        m=\norm{f}_{W(A)},
        \qquad
        s=\frac Km.
\]
Since $m\leq1$ and the scalar Crouzeix inequality applied to $f$ on
$W(A)$ gives $K\leq2m$, one has $m>0$ and $1<s\leq2$.
For a Borel map $g$ and a measure $\lambda$, the
pushforward measure $g_*\lambda$ is defined by
$g_*\lambda(E)=\lambda(g^{-1}(E))$ for every Borel set $E$. We write
$m_{\T}$ for normalized arclength measure on $\T$.

\begin{theorem}\label{thm:measure-stability}
Set
\[
        \varepsilon_{K,m}=m^2\vartheta(s),
        \qquad
        s=\frac Km.
\]
Then, for every $h\in\Acal(\overline\Omega)$,
\begin{equation}\label{eq:approximate-annihilation}
 \left|\int_{\partial\Omega}fh\dd\mu\right|
 \leq\varepsilon_{K,m}\norm{h}_{\overline\Omega}.
\end{equation}
Moreover,
\begin{equation}\label{eq:measure-coarse}
 0\leq\varepsilon_{K,m}
 \leq2m\sqrt{4m^2-K^2}
 \leq2\sqrt{4-K^2}.
\end{equation}
If $\nu=f_*\mu$, then $\nu$ is a probability measure on $\T$ and
\begin{equation}\label{eq:fourier-bound}
        \left|\int_{\T}\zeta^k\dd\nu(\zeta)\right|
        \leq\varepsilon_{K,m}
        \qquad(k\in\mathbb Z\setminus\{0\}).
\end{equation}
Let $(A_j,\Omega_j,f_j,x_j)$ be a sequence satisfying the
assumptions above, and define
$K_j,m_j,s_j,\mu_j$, and $\nu_j=(f_j)_*\mu_j$ as above. If
$s_j=K_j/m_j\to2$, then $\nu_j$ converges weakly to $m_{\T}$.
\end{theorem}

\begin{proof}
The normalized matrix $T=F/m$ satisfies the two square-function
estimates for every $0<\tau<1$, has norm $s$, and attains its norm at
$x$. Theorem~\ref{thm:sharp-T2} gives
\[
        \norm{T^2x}\leq\vartheta(s).
\]
Since $Fy=F^2x/K$,
\begin{equation}\label{eq:near-extremal-Fy}
        \norm{Fy}\leq\frac{m^2\vartheta(s)}K.
\end{equation}
By \eqref{eq:cross-state} and commutativity,
\[
\begin{aligned}
 \int_{\partial\Omega}fh\dd\mu
 &=K\ip{h(A)y}{x}\\
 &=\ip{Fh(A)y}{y}\\
 &=\ip{h(A)Fy}{y}.
\end{aligned}
\]
The definition of $K$ gives
$\norm{h(A)}\leq K\norm{h}_{\overline\Omega}$. Together with
\eqref{eq:near-extremal-Fy}, this proves
\eqref{eq:approximate-annihilation}.

Since $\vartheta(s)\leq\eta(s)$,
Lemma~\ref{lem:tau-optimization} gives
\[
        \varepsilon_{K,m}
        \leq2m^2\sqrt{4-\frac{K^2}{m^2}}
        =2m\sqrt{4m^2-K^2}.
\]
Since $m\leq1$ and $K\leq2m$, both factors on the right in
\[
 4-K^2-m^2(4m^2-K^2)
 =(1-m^2)\bigl(4(1+m^2)-K^2\bigr)
\]
are nonnegative. This proves the second inequality in
\eqref{eq:measure-coarse}.

By \eqref{eq:boundary-inner}, $\nu$ is supported on $\T$. Taking
$h=f^{k-1}$ in \eqref{eq:approximate-annihilation} gives
\eqref{eq:fourier-bound} for $k\geq1$. The estimates for negative $k$
follow by conjugation. For a sequence in the final assertion,
$\vartheta(s_j)\to0$ and $m_j\leq1$, so
$\varepsilon_{K_j,m_j}\to0$. Thus every nonzero Fourier coefficient of
$\nu_j$ tends to zero. Weak convergence follows from the uniform
density of the trigonometric polynomials in $C(\T)$.
\end{proof}

\begin{corollary}\label{cor:measure-equality}
If
\[
        K=2m,
        \qquad
        m=\norm{f}_{W(A)},
\]
then
\[
        \int_{\partial\Omega}fh\dd\mu=0
        \qquad(h\in\Acal(\overline\Omega)),
\]
and
\[
        f_*\mu=m_{\T}.
\]
\end{corollary}

\begin{proof}
The assumption is $s=2$, so $\vartheta(s)=0$. The first assertion
follows from Theorem~\ref{thm:measure-stability}. Taking $h=f^{k-1}$
gives vanishing of all positive Fourier coefficients of $f_*\mu$. The
negative coefficients are their complex conjugates, so
$f_*\mu=m_{\T}$.
\end{proof}

\section{Polynomial equality and equality in dimension three}\label{sec:dimension-three}
We first characterise equality for scalar polynomials in every finite
dimension. We then combine the scalar rigidity theorem with the
complete bound to classify equality in dimension three.

We write $A\simeq B$ for unitary similarity. For a matrix $A$, put
$w(A)=\max_{z\in W(A)}|z|$. For $m\geq1$, let
$C_m\in M_{m+1}(\C)$ have nonzero entries only on the first
superdiagonal. For $m=1$, its only nonzero entry is $2$. For $m\geq2$,
its first and last superdiagonal entries are $\sqrt2$, and all
intervening entries are $1$. These are the Crabb--Choi--Crouzeix
matrices. They satisfy
\[
        W(C_m)=\overline\D,
        \qquad \norm{C_m^m}=2
\]
\cite[Section~4]{HnetynkovaTichy2018}.

\begin{proposition}\label{v3:prop:polynomial-disk}
Let $A\in M_d(\C)$ satisfy $W(A)\subset\overline\D$, and let $p$ be
a nonzero polynomial. Then
\[
        \norm{p(A)}=2\max_{|z|\leq1}|p(z)|
\]
if and only if there are $c\neq0$ and an integer $1\leq m<d$ such that
$p(z)=cz^m$ and $A$ is unitarily similar to $C_m\oplus B$, where
$w(B)\leq1$. The block $B$ is omitted when $m+1=d$.
\end{proposition}

\begin{proof}
Divide $p$ by its supremum norm on $\overline\D$. Thus
$\norm{p}_{\overline\D}=1$ and $\norm{p(A)}=2$. Suppose that
$Av=\lambda v$, $\norm{v}=1$, and $|\lambda|=1$. Since
\[
 I-\Ree(\overline\lambda A)\geq0,
 \qquad
 \ip{(I-\Ree(\overline\lambda A))v}{v}=0,
\]
positivity gives $\Ree(\overline\lambda A)v=v$. It follows that
$A^*v=\overline\lambda v$, so $\C v$ reduces $A$. Removing these
reducing subspaces successively gives
\[
        A\simeq A_0\oplus D,
        \qquad A_0\in M_{d_0}(\C),
        \qquad \sigma(A_0)\subset\D,
\]
where $D$ is diagonal and $\sigma(D)\subset\T$. Either block may be
absent before equality is used. Since $\norm{p(D)}\leq1$, equality
implies that $A_0$ is present and $\norm{p(A_0)}=2$.

Apply Proposition~\ref{prop:domain-extremals} to $(A_0,\D)$. Its
bound gives $\psi_\D(A_0)\leq2$, and $p$ attains equality. Thus $p$
is a normalized extremal function. By the same proposition, which
uses Crouzeix's finite-dimensional extremal theorem
\cite[Theorem~2.1]{Crouzeix2004}, $p$ is a finite Blaschke product
of degree at most $d_0-1$. A polynomial finite Blaschke product
is a unimodular monomial. Indeed, any nonzero zero $\alpha\in\D$
in its Blaschke factorisation would produce a pole at
$1/\overline\alpha$, which cannot cancel against a numerator zero
in $\D$. Consequently $p(z)=\eta z^m$, where $|\eta|=1$. A constant
polynomial cannot have norm $2$ at $A_0$, so $1\leq m<d$.

Now $\norm{A^m}=2$. The numerical-radius power inequality gives
\[
        2=\norm{A^m}\leq2w(A^m)\leq2w(A)^m\leq2,
\]
so $w(A)=1$. The equality theorem for powers, proved from Crabb's
theorem in \cite[Theorems~12 and~14]{HnetynkovaTichy2018}, gives a
unitary decomposition $A\simeq C_m\oplus B$, with $w(B)\leq1$.
This proves necessity and allows us to remove the normalization of
$p$.

Conversely, if $A\simeq C_m\oplus B$ and $w(B)\leq1$, then
$W(A)=\overline\D$. The numerical-radius power inequality gives
$\norm{B^m}\leq2$, while $\norm{C_m^m}=2$. Hence
$\norm{p(A)}=2|c|=2\norm{p}_{\overline\D}$ for $p(z)=cz^m$.
\end{proof}

Jin's rigidity theorem states that $\psi(A)=2$ forces $W(A)$ to be a
nondegenerate disk in every finite dimension
\cite[Theorem~1]{JinRigidity2026}. Together with the preceding
proposition, it gives the following classification. The matrix
decomposition for monomials is the theorem of Hn{\v e}tynkov{\'a} and
Tich{\'y} cited above. The additional conclusion is that, after an
affine change of variable, every polynomial attaining equality must
be a monomial.

\begin{theorem}\label{v3:thm:polynomial-classification}
Let $A\in M_d(\C)$ and let $p$ be a polynomial with $p(A)\neq0$.
Then
\[
        \norm{p(A)}=2\max_{z\in W(A)}|p(z)|
\]
if and only if there are $a\in\C$, $r>0$, $c\neq0$, and an integer
$1\leq m<d$ such that
\[
        p(z)=c(z-a)^m,
        \qquad
        A\simeq aI+r(C_m\oplus B),
        \qquad w(B)\leq1.
\]
The block $B$ is omitted when $m+1=d$.
\end{theorem}

\begin{proof}
Equality and the scalar bound give $\psi(A)=2$. Jin's rigidity
theorem \cite[Theorem~1]{JinRigidity2026} gives
$W(A)=a+r\overline\D$, with $r>0$. Put $T=(A-aI)/r$ and
$q(z)=p(a+rz)$. Proposition~\ref{v3:prop:polynomial-disk} gives
$q(z)=c_0z^m$ and $T\simeq C_m\oplus B$, with $1\leq m<d$ and
$w(B)\leq1$. Therefore $p(z)=c_0r^{-m}(z-a)^m$, as required.
The converse follows from the converse in
Proposition~\ref{v3:prop:polynomial-disk} by the same affine change
of variable.
\end{proof}

The next consequence answers the degree question in
\cite[Section~8.4]{GreenbaumOverton2018}.

\begin{corollary}\label{v3:cor:polynomial-degree}
Let $A\in M_d(\C)$ and let $p$ be a polynomial with
$\deg p\geq d$ and $p(A)\neq0$. Then
\[
        \norm{p(A)}<2\max_{z\in W(A)}|p(z)|.
\]
\end{corollary}

\begin{proof}
The scalar bound gives the non-strict inequality. Equality would
give $\deg p<d$ by Theorem~\ref{v3:thm:polynomial-classification}.
\end{proof}

In dimension three, Theorem~\ref{thm:intro-complete-equality}
transfers Jin's rigidity theorem to the complete constant.

\begin{corollary}\label{v3:cor:nondisk-strict}
Let $A\in M_d(\C)$, where $d\leq3$. If $W(A)$ is not a
nondegenerate disk, then $\psi_{\mathrm{cb}}(A)<2$.
\end{corollary}

\begin{proof}
Theorem~\ref{thm:intro-complete-three} gives
$\psi_{\mathrm{cb}}(A)\leq2$. Equality would give $\psi(A)=2$ by
Theorem~\ref{thm:intro-complete-equality}, and then $W(A)$ would be
a nondegenerate disk by \cite[Theorem~1]{JinRigidity2026}.
\end{proof}

Crouzeix classified the matrices of order three with unit-disk
numerical range and scalar constant two
\cite[Lemma~10.2]{Crouzeix2016}. The preceding results remove the
disk assumption and give the complete equality classification.

\begin{corollary}\label{v3:cor:three-classification}
For $A\in M_3(\C)$, the following statements are equivalent.
\begin{enumerate}[(i)]
\item $\psi(A)=2$.
\item $\psi_{\mathrm{cb}}(A)=2$.
\item There are $a\in\C$ and $b\in\C\setminus\{0\}$ such that $A$
is unitarily similar to $aI+bB_\xi$ or $aI+bM_s$, where
\[
 B_\xi=\begin{pmatrix}0&0&2\\0&\xi&0\\0&0&0\end{pmatrix},
 \qquad |\xi|\leq1,
\]
and
\[
 M_s=\begin{pmatrix}
 0&\sqrt{2(1-s^2)}&2s\\
 0&-s&\sqrt{2(1-s^2)}\\
 0&0&0
 \end{pmatrix},
 \qquad 0\leq s\leq1.
\]
\end{enumerate}
\end{corollary}

\begin{proof}
Statements (i) and (ii) are equivalent by
Theorem~\ref{thm:intro-complete-equality}. If (i) holds, Jin's
rigidity theorem gives $W(A)=a+r\overline\D$, with $r>0$.
Apply \cite[Lemma~10.2]{Crouzeix2016} to $(A-aI)/r$. Its two
families are precisely $B_\xi$ and $e^{i\vartheta}M_s$, with
$s=\sin\phi$, $0\leq\phi\leq\pi/2$. Absorbing the unimodular
factor into $b$ gives (iii). Conversely, the same lemma gives
$\psi(B_\xi)=\psi(M_s)=2$. Unitary similarity and a nonconstant
affine change preserve $\psi$, so (iii) implies (i).
\end{proof}

We finish with algebraic consequences of equality that apply in
arbitrary finite dimension. Write $r_A(\lambda)$ for the size of the
largest Jordan block of $A$ at $\lambda$, and put
$D(z_0,r)=\{z\in\C:|z-z_0|<r\}$.

\begin{lemma}\label{lem:algebra-generation}
Let $A$ be a matrix and let $f$ be holomorphic on a neighborhood of
$\sigma(A)$. Then
\[
        \alg(f(A))=\alg(A)
\]
if and only if
\begin{enumerate}[(i)]
\item $f$ is injective on $\sigma(A)$.
\item $f'(\lambda)\neq0$ whenever $r_A(\lambda)\geq2$.
\end{enumerate}
\end{lemma}

\begin{proof}
The inclusion $\alg(f(A))\subset\alg(A)$ is automatic. Suppose first
that $A=q(f(A))$ for a polynomial $q$. If
$\lambda,\nu\in\sigma(A)$, then
\[
        q(f(\lambda))=\lambda,
        \qquad
        q(f(\nu))=\nu.
\]
Thus $f$ is injective on $\sigma(A)$. If $r_A(\lambda)\geq2$,
comparison of first derivatives on a Jordan block at $\lambda$ gives
$q'(f(\lambda))f'(\lambda)=1$.

Conversely, assume (i) and (ii). For each $\lambda\in\sigma(A)$, the
function $f$ has a local holomorphic inverse at $f(\lambda)$ when
$r_A(\lambda)\geq2$. When $r_A(\lambda)=1$, only the value of an
inverse is needed. Since the points $f(\lambda)$ are distinct, Hermite
interpolation gives a polynomial $q$ whose derivatives through order
$r_A(\lambda)-1$ agree at $f(\lambda)$ with those of the corresponding
local inverse. It follows that $q(f(A))=A$.
\end{proof}

\begin{theorem}\label{thm:algebra-disk}
Let $A\in M_n(\C)$ and let $\Omega$ be a bounded Jordan domain. Assume
\[
        W(A)\subset\overline\Omega,
        \qquad
        \sigma(A)\subset\Omega.
\]
Let $f\in\Acal(\overline\Omega)$ satisfy
\[
        \norm{f}_{\overline\Omega}=1,
        \qquad
        \norm{f(A)}=2.
\]
If $\alg(f(A))=\alg(A)$, then there are $z_0\in\Omega$ and $c\neq0$
such that
\[
        \Omega=z_0+c\D,
        \qquad
        f(z_0+c\zeta)=\zeta
        \quad(\zeta\in\D).
\]
\end{theorem}

\begin{proof}
Put $F=f(A)$. Theorem~\ref{thm:tuple-square} gives
\[
        2=\norm{F}\leq2\norm{f}_{W(A)}\leq2.
\]
Thus $F$ satisfies the two square-function estimates at $\tau=1/2$.
By Theorem~\ref{thm:equality-block}, there is a two-dimensional reducing
subspace $V_0$ such that, in a suitable orthonormal basis,
\[
        F|_{V_0}=N,
        \qquad
        N=\begin{pmatrix}0&0\\2&0\end{pmatrix}.
\]
Choose a polynomial $q$ such that $A=q(F)$. Then $V_0$ also reduces
$A$, and
\[
        A|_{V_0}=q(N)=z_0I+cN,
        \qquad z_0=q(0),
        \qquad c=q'(0).
\]
The identity $f(A|_{V_0})=N$ gives
$f(z_0)=0$ and $cf'(z_0)=1$. Hence $c\neq0$. The number $z_0$ is the
eigenvalue of $A|_{V_0}$, so $z_0\in\sigma(A)\subset\Omega$.

Since $W(N)=\overline\D$,
\[
        \{z:|z-z_0|\leq|c|\}=W(A|_{V_0})
        \subset W(A)\subset\overline\Omega.
\]
Since the interior of the closure of a Jordan domain is the domain
itself, $D(z_0,|c|)\subset\Omega$. The function
$g(\zeta)=f(z_0+c\zeta)$ satisfies $|g|\leq1$ on $\D$. Since
$g'(0)=1$, it is nonconstant, and the maximum principle shows that it
maps $\D$ into $\D$. Moreover, $g(0)=0$ and $g'(0)=1$. Equality in
Schwarz's lemma gives $g(\zeta)=\zeta$. The identity theorem gives
$f(z)=(z-z_0)/c$ on $\Omega$. The bound $\norm{f}_{\overline\Omega}=1$
then gives $\Omega\subset z_0+c\D$. The reverse inclusion was proved
above, so $\Omega=z_0+c\D$.
\end{proof}

\begin{corollary}\label{cor:algebra-obstruction}
Under the assumptions of Theorem~\ref{thm:algebra-disk}, suppose that
$\Omega$ is not a disk. Then either
\begin{enumerate}[(i)]
\item $f(\lambda)=f(\nu)$ for two distinct eigenvalues
$\lambda,\nu\in\sigma(A)$, or
\item $f'(\lambda)=0$ for some $\lambda\in\sigma(A)$ with
$r_A(\lambda)\geq2$.
\end{enumerate}
\end{corollary}

\begin{proof}
Otherwise Lemma~\ref{lem:algebra-generation} would give
$\alg(f(A))=\alg(A)$, and Theorem~\ref{thm:algebra-disk} would imply
that $\Omega$ is a disk.
\end{proof}

\begin{proposition}\label{prop:jordan-restriction}
Let $A$ be a matrix, and let $f\in\Acal(W(A))$ be holomorphic on a
neighborhood of $\sigma(A)$. Assume
\[
        M=\max_{z\in W(A)}\abs{f(z)}>0,
        \qquad
        \norm{f(A)}=2M.
\]
Then there are an eigenvalue $\lambda$ of $A$, a Jordan block of size
$r$ at $\lambda$, and an integer $m\geq1$ such that $f$ has a zero of
order $m$ at $\lambda$ and
\begin{equation}\label{eq:jordan-order}
        m<r<3m.
\end{equation}
\end{proposition}

\begin{proof}
After replacing $f$ by $f/M$, Theorem~\ref{thm:equality-block} shows
that the Jordan form of $f(A)$ contains a block $J_2(0)$. Applying $f$
to the Jordan blocks of $A$, there must be a block $J_r(\lambda)$ for
which the Jordan form of $f(J_r(\lambda))$ contains $J_2(0)$. Hence
$f(\lambda)=0$. If the zero has order $m$, write
\[
        f(J_r(\lambda))=N^mg(N),
        \qquad
        N=J_r(\lambda)-\lambda I,
\]
where $g(N)$ is invertible and commutes with $N$. For every $k\geq1$,
\[
        (N^mg(N))^k=N^{mk}g(N)^k.
\]
The matrix $g(N)^k$ is invertible and commutes with $N$. Hence
$(N^mg(N))^k$ and $N^{mk}$ have kernels of the same dimension for every
$k$, so $f(J_r(\lambda))$ and $N^m$ have the same Jordan block sizes.

Write $r=qm+s$, where $q\geq0$ and $0\leq s<m$. The Jordan blocks
of $N^m$ have sizes $q+1$, repeated $s$ times, and $q$, repeated
$m-s$ times. Blocks of size zero are omitted. A block of size two
occurs
exactly when either $q=1$ and $s>0$, or $q=2$ and $s<m$. These
conditions are equivalent to \eqref{eq:jordan-order}.
\end{proof}

\section{A family with a line segment on the boundary}\label{sec:family}

We apply the preceding results to
\begin{equation}\label{eq:family}
 A_{\mu,r}=\begin{pmatrix}
 i\mu&0&ir\\
 0&-i\mu&ir\\
 ir&ir&-1
 \end{pmatrix},
 \qquad \mu>0,
 \qquad r\geq0.
\end{equation}
For a unit vector $x=(x_1,x_2,x_3)^{\mathsf T}$,
\[
        \operatorname{Re}(x^*A_{\mu,r}x)=-\abs{x_3}^2.
\]
Hence the line $\operatorname{Re}z=0$ supports the numerical range,
and the vectors with $x_3=0$ give
\begin{equation}\label{eq:boundary-segment}
        W(A_{\mu,r})\cap i\mathbb R=[-i\mu,i\mu].
\end{equation}
Thus the displayed line segment is contained in the boundary of the
numerical range. If $r=0$, then
\[
        A_{\mu,0}=\diag(i\mu,-i\mu,-1)
\]
is normal, and therefore
\begin{equation}\label{eq:family-normal-case}
        \psi(A_{\mu,0})=\psi_{\mathrm{cb}}(A_{\mu,0})=1.
\end{equation}
\begin{corollary}\label{v3:family-strict}
For every $\mu>0$ and $r\geq0$,
\[
        \psi_{\mathrm{cb}}(A_{\mu,r})<2.
\]
\end{corollary}

\begin{proof}
The nontrivial boundary segment in \eqref{eq:boundary-segment}
excludes a nondegenerate disk. Apply
Corollary~\ref{v3:cor:nondisk-strict}.
\end{proof}

The remainder of the section concerns $r>0$.

\begin{lemma}\label{lem:family-irreducible}
If $r>0$, then $A_{\mu,r}$ is unitarily irreducible.
\end{lemma}

\begin{proof}
Let $X$ commute with $A_{\mu,r}$ and $A_{\mu,r}^*$. It commutes with
$\operatorname{Re}A_{\mu,r}=\diag(0,0,-1)$, so
$X=B\oplus[c]$ with $B\in M_2(\C)$. It also commutes with
\[
 \frac{A_{\mu,r}-A_{\mu,r}^*}{2i}
 =\begin{pmatrix}
 \mu&0&r\\
 0&-\mu&r\\
 r&r&0
 \end{pmatrix}.
\]
The upper-left block shows that $B$ is diagonal. The upper-right block
gives $B(r,r)^{\mathsf T}=c(r,r)^{\mathsf T}$. Since $r>0$, one has
$B=cI_2$. Thus the joint commutant of $A_{\mu,r}$ and
$A_{\mu,r}^*$ consists of scalar matrices. An orthogonal projection
onto a reducing subspace belongs to this joint commutant.
\end{proof}

The family is unitarily equivalent, up to sign, to a real family used
by Crouzeix \cite[Section~9, p.~436]{Crouzeix2016}.
For $a>0$ and $b\geq0$, put
\[
 C_{a,b}=\begin{pmatrix}
 0&a&b\\
 -a&0&b\\
 -b&-b&1
 \end{pmatrix}.
\]

\begin{proposition}\label{prop:crouzeix-normal-form}
The matrix $A_{a,b}$ is unitarily equivalent to $-C_{a,b}$. Hence
\[
        \psi(A_{a,b})=\psi(C_{a,b}),
        \qquad
        \psi_{\mathrm{cb}}(A_{a,b})=\psi_{\mathrm{cb}}(C_{a,b}).
\]
\end{proposition}

\begin{proof}
Let
\[
 u_1=\frac{e^{i\pi/4}}{\sqrt2}
        \begin{pmatrix}i\\1\\0\end{pmatrix},
 \quad
 u_2=\frac{e^{3i\pi/4}}{\sqrt2}
        \begin{pmatrix}-i\\1\\0\end{pmatrix},
 \quad
 u_3=\begin{pmatrix}0\\0\\1\end{pmatrix},
\]
and put $U=(u_1,u_2,u_3)$. Direct multiplication gives
$U^*(-C_{a,b})U=A_{a,b}$. The two ratios are invariant under unitary
equivalence and under multiplication of the matrix by a nonzero
scalar.
\end{proof}

Put
\[
        h=\frac r\mu,
        \qquad q=\mu^2,
        \qquad y=h^2.
\]
With respect to the orthonormal basis
\[
        \frac{e_1-e_2}{\sqrt2},
        \qquad
        \frac{e_1+e_2}{\sqrt2},
        \qquad e_3,
\]
the matrix $A_{\mu,r}$ is unitarily equivalent to
\begin{equation}\label{eq:tilde-A}
 \widetilde A_{\mu,h}=\begin{pmatrix}
 0&i\mu&0\\
 i\mu&0&i\sqrt2h\mu\\
 0&i\sqrt2h\mu&-1
 \end{pmatrix}.
\end{equation}
Its characteristic polynomial is
\begin{equation}\label{eq:family-charpoly}
 p_{q,y}(\lambda)
 =\lambda^3+\lambda^2+(1+2y)q\lambda+q,
\end{equation}
and its discriminant is
\begin{equation}\label{eq:family-discriminant}
\begin{aligned}
 \Delta(q,y)
 =-4q\{&8q^2y^3+(12q^2-q)y^2\\
       &+(6q^2-10q)y+(q+1)^2\}.
\end{aligned}
\end{equation}
If $h>0$, neither $0$ nor $-1$ is a root of
$p_{q,y}$. For every root $\lambda$,
\begin{equation}\label{eq:family-eigenvector}
 x(\lambda)=
 \left(\frac{i\mu}{\lambda},\ 1,\
       \frac{i\sqrt2h\mu}{1+\lambda}\right)^{\mathsf T}
\end{equation}
is an eigenvector of $\widetilde A_{\mu,h}$. Every eigenspace is
one-dimensional. Consequently, $A_{\mu,r}$ is diagonalizable exactly
when $\Delta(q,y)\ne0$.

When the roots $\lambda_1,\lambda_2,\lambda_3$ are simple, let
$v_j=x(\lambda_j)/\norm{x(\lambda_j)}$ and put
$V=(v_1,v_2,v_3)$. Then
\[
        \widetilde A_{\mu,h}=V\diag(\lambda_1,\lambda_2,\lambda_3)V^{-1}.
\]
Theorem~\ref{thm:intro-complete-diagonalizable} gives the following
complete estimate.

\begin{corollary}\label{cor:family-diagonalizable}
Assume $r>0$ and $\Delta(q,y)\ne0$. With $V$ as above and
$\kappa=\cond(V)$,
\begin{equation}\label{eq:family-complete-bound}
        \psi_{\mathrm{cb}}(A_{\mu,r})
        \leq\min\{\kappa,\Phi(\kappa)\}<2.
\end{equation}
\end{corollary}

\begin{lemma}\label{lem:boundary-normal-eigenvalue}
Let $A\in M_n(\C)$, let
$\lambda\in\sigma(A)\cap\partial W(A)$, and let $u\ne0$ satisfy
$Au=\lambda u$. Then $A^*u=\overline\lambda u$. In particular, if
$A$ is unitarily irreducible and $n\geq2$, then
$\sigma(A)\subset W(A)^\circ$.
\end{lemma}

\begin{proof}
After translation and rotation, suppose that $\lambda=0$ and that the
numerical range lies in the closed right half-plane. Then
$\operatorname{Re}A\succeq0$. If $Au=0$, positivity and
$\ip{(\operatorname{Re}A)u}{u}=0$ give
$(\operatorname{Re}A)u=0$, and hence $A^*u=0$. Thus the span of $u$
reduces $A$.
\end{proof}

For $r>0$, Lemmas~\ref{lem:family-irreducible}
and~\ref{lem:boundary-normal-eigenvalue} give
\[
        \sigma(A_{\mu,r})\subset W(A_{\mu,r})^\circ.
\]
\begin{proposition}\label{prop:family-defective-curve}
Assume $r>0$ and $\Delta(q,y)=0$. There is a unique
$t\in(0,1/2)$ such that
\begin{equation}\label{eq:defective-curve}
        q=t^2(1-2t),
        \qquad
        y=\frac{(1-t)^2}{t(1-2t)}.
\end{equation}
Equivalently,
\begin{equation}\label{eq:defective-parameters}
        \mu=t\sqrt{1-2t},
        \qquad
        r=(1-t)\sqrt t.
\end{equation}
The repeated eigenvalue is $-t$ and the remaining eigenvalue is
$2t-1$. Moreover,
\begin{equation}\label{eq:defective-location}
        q\leq\frac1{27},
        \qquad
        y\geq4,
\end{equation}
with equality in either inequality exactly when $t=1/3$. If
$t\ne1/3$, the Jordan form is $J_2(-t)\oplus[2t-1]$. If $t=1/3$, the
Jordan form is $J_3(-1/3)$.
\end{proposition}

\begin{proof}
Let $\lambda$ be a repeated root of $p_{q,y}$. Since the polynomial has
real coefficients and degree three, $\lambda$ is real. The equations
$p_{q,y}(\lambda)=p_{q,y}'(\lambda)=0$ give
\[
        q=\lambda^2(2\lambda+1),
        \qquad
        y=-\frac{(\lambda+1)^2}{\lambda(2\lambda+1)}.
\]
Since $q>0$ and $y>0$, one has $-1/2<\lambda<0$. Put $t=-\lambda$.
Then
\[
        p_{q,y}(z)=(z+t)^2(z-2t+1),
\]
and \eqref{eq:defective-curve} follows. The formulas determine $t$
uniquely. Furthermore,
\[
 \frac1{27}-q
 =\frac{(1-3t)^2(1+6t)}{27},
 \qquad
 y-4
 =\frac{(1-3t)^2}{t(1-2t)},
\]
which proves \eqref{eq:defective-location}. The eigenvector calculation
in \eqref{eq:family-eigenvector} shows that every eigenspace is
one-dimensional. This gives the two Jordan forms.
\end{proof}

At the triple value, the matrix is affinely equivalent to the
$3\times3$ Kac--Murdock--Szeg\H{o} matrix considered in
\cite{CrouzeixGreenbaumLi2024}.

\begin{proposition}\label{prop:triple-KMS}
Let
\[
        \mu_* =\frac1{3\sqrt3},
        \qquad
        r_* =\frac2{3\sqrt3},
        \qquad
 K_3=\begin{pmatrix}
 0&1&1\\
 0&0&1\\
 0&0&0
 \end{pmatrix}.
\]
Then
\[
        A_{\mu_*,r_*}+\frac13I
        \simeq_u-\frac23K_3.
\]
Consequently,
\[
        \psi(A_{\mu_*,r_*})=\psi(K_3),
        \qquad
        \psi_{\mathrm{cb}}(A_{\mu_*,r_*})
        =\psi_{\mathrm{cb}}(K_3).
\]
\end{proposition}

\begin{proof}
Set
\[
U=\frac1{\sqrt3}\begin{pmatrix}
-\frac12-\frac{\sqrt3}{2}i&1&-\frac12+\frac{\sqrt3}{2}i\\
 \frac12-\frac{\sqrt3}{2}i&-1& \frac12+\frac{\sqrt3}{2}i\\
1&1&1
\end{pmatrix}.
\]
The columns of $U$ are orthonormal, and direct multiplication gives
\[
 U^*\left(A_{\mu_*,r_*}+\frac13I\right)U=-\frac23K_3.
\]
The two ratios are invariant under unitary equivalence and under
$A\mapsto\alpha A+\beta I$, where $\alpha\ne0$.
\end{proof}

\begin{lemma}\label{lem:KMS-rational-domain}
Let
\[
 f_1(z)=\frac{\sum_{j=1}^7 c_jz^j}
 {1+\sum_{j=1}^7d_jz^j},
\]
where
\[
\begin{aligned}
10^5(c_1,\ldots,c_7)
  &=(73400,49736,7268,-521,13,61,-251),\\
10^5(d_1,\ldots,d_7)
  &=(32564,-3291,1000,-400,84,-242,28).
\end{aligned}
\]
Then $f_1$ is holomorphic on a neighborhood of $\overline\D$ and
\begin{equation}\label{eq:rational-inclusion}
        f_1(\overline\D)\subset W(K_3).
\end{equation}
\end{lemma}

\begin{proof}
Write
\[
 d_0(z)=1+\sum_{j=1}^7d_jz^j.
\]
For $\abs z\leq1$, the coefficients give
\[
 \abs{d_0(z)}\geq1-\sum_{j=1}^7\abs{d_j}>0.
\]
Thus $d_0$ has no zero in $\overline\D$, and $f_1$ is holomorphic
on a neighborhood of $\overline\D$.
\par

We verify \eqref{eq:rational-inclusion} in exact rational arithmetic.
Parametrize the upper semicircle by
\[
        \zeta(u)=\frac{1+iu}{1-iu},
        \qquad 0\leq u<\infty,
\]
and put $w(u)=f_1(\zeta(u))$. The common denominator of
$\operatorname{Re}w(u)$ and $\operatorname{Im}w(u)$ is
\[
\begin{aligned}
D(u)={}&3892636881u^{14}+53940184519u^{12}
 +168968707589u^{10}\\
&+416744879107u^8+436966909139u^6
 +280138656021u^4\\
&+111953512087u^2+16833246049,
\end{aligned}
\]
which is positive on $[0,\infty)$. Let $C,S_1,S_2,S_3\in\mathbb Z[u]$
be the primitive numerators with positive leading coefficients,
after cancelling common factors in the rational expressions
obtained by substituting $w=w(u)$ into
\[
\begin{aligned}
&27\abs{w}^4-18\abs{w}^2-8\operatorname{Re}w-1,\\
&\operatorname{Re}w+\frac12,
\qquad -\operatorname{Re}w,
\qquad \frac1{12}-(\operatorname{Im}w)^2,
\end{aligned}
\]
respectively. These expressions are
$C/D^2$, $S_1/(2D)$, $2S_2/D$, and $S_3/(12D^2)$.
Their signs therefore agree with those of the corresponding polynomials.
We use Sturm sequences in exact rational arithmetic
\cite[Chapter~2]{BasuPollackRoy2006}. The following table gives the
resulting variations.
The entry $+\infty$ means that the signs of the leading coefficients
are used.
\begin{center}
\begin{tabular}{c c c c c c}
\toprule
polynomial & degree & interval & left variation & right variation & sign \\
\midrule
$C$   & $28$ & $[0,5/2]$ & $14$ & $14$ & $C(0)<0$ \\
$S_1$ & $14$ & $[12/5,+\infty)$& $7$  & $7$  & $S_1(12/5)>0$ \\
$S_2$ & $14$ & $[12/5,+\infty)$& $6$  & $6$  & $S_2(12/5)>0$ \\
$S_3$ & $28$ & $[12/5,+\infty)$& $12$ & $12$ & $S_3(12/5)>0$ \\
\bottomrule
\end{tabular}
\end{center}
Thus $w(u)$ belongs to
\[
 \mathcal C=\{z:27\abs z^4-18\abs z^2-8\operatorname{Re}z-1\leq0\}
\]
for $0\leq u\leq5/2$, and it belongs to
\[
 \mathcal R=\left\{z:-\frac12\leq\operatorname{Re}z\leq0,
          \ \abs{\operatorname{Im}z}\leq\frac{\sqrt3}{6}\right\}
\]
for $u\geq12/5$. The intervals overlap.
Letting $u\to\infty$ and using the continuity
of $f_1$ and the closedness of $\mathcal R$, we obtain
$f_1(-1)\in\mathcal R$. Both $\mathcal C$ and
$\mathcal R$ are contained in $W(K_3)$
\cite[Section~2]{CrouzeixGreenbaumLi2024}. The coefficients of $f_1$
are real, and $W(K_3)$ is symmetric with respect to the real axis.
Hence $f_1(\T)\subset W(K_3)$. Applying the maximum principle to every
supporting half-plane of the convex set $W(K_3)$ proves
\eqref{eq:rational-inclusion}.

\end{proof}

The next proof gives an independent certificate of strictness at the
triple point, using Lemma~\ref{lem:KMS-rational-domain}. It also gives
a completely bounded counterpart to Crouzeix's numerical results
\cite{Crouzeix2013}.

\begin{theorem}\label{thm:triple-strict}
At the triple point,
\begin{equation}\label{eq:triple-strict}
 \psi_{\mathrm{cb}}(A_{\mu_*,r_*})<2.
\end{equation}
Moreover, there is an $\varepsilon>0$ such that
\[
 |\mu-\mu_*|+|r-r_*|<\varepsilon
 \quad\Longrightarrow\quad
 \psi_{\mathrm{cb}}(A_{\mu,r})<2.
\]
\end{theorem}

\begin{proof}
To prove~\eqref{eq:triple-strict}, by Proposition~\ref{prop:triple-KMS}, it suffices to show that $\psi_{\mathrm{cb}}(K_3)<2$. Let $\Omega=W(K_3)^\circ$. By Proposition~\ref{prop:triple-KMS} and \eqref{eq:boundary-segment}, $W(K_3)$ has a line segment on its boundary and hence is not a disk. Trivially, $K_3$ is unitarily irreducible, its Jordan form is $J_3(0)$, and $\sigma(K_3)=\{0\}$. Moreover, Lemma~\ref{lem:boundary-normal-eigenvalue}
implies that $0\in\Omega$.

Suppose that $\psi_{\mathrm{cb}}(K_3)=2$. By
Theorem~\ref{thm:intro-complete-equality}, $\psi(K_3)=2$.
Proposition~\ref{prop:domain-extremals} gives a normalized extremal
function $f$ of Blaschke degree at most two. Since the Jordan form of
$K_3$ is $J_3(0)$, Proposition~\ref{prop:jordan-restriction} forces
$f$ to have a double zero at zero. Choose a conformal map
$\varphi:\Omega\to\D$ with $\varphi(0)=0$. It follows that
$f=\eta\varphi^2$, where $|\eta|=1$. As $K_3^3=0$,
\[
 f(K_3)=\eta\varphi'(0)^2K_3^2.
\]
Using $\norm{f(K_3)}=2$ and $\norm{K_3^2}=1$, we obtain
\[
 \operatorname{crad}_{\Omega}(0)^2
 :=\frac1{|\varphi'(0)|^2}=\frac12.
\]

On the other hand, Lemma~\ref{lem:KMS-rational-domain} gives a nonconstant holomorphic function $f_1$ such that
$$
    f_1(\overline{\D})\subset W(K_3),\qquad
    f_1(0)=0,\qquad
    f_1'(0)=\frac{367}{500}.
$$
By the open mapping theorem, $f_1(\D)\subset\Omega$. Let
$\varphi:\Omega\to\D$ be conformal with $\varphi(0)=0$. Schwarz's lemma applied to $\varphi\circ f_1$ implies that
$$
|\varphi'(0)|\,\frac{367}{500}\leq1,
$$
and hence
$$
\operatorname{crad}_{\Omega}(0)
=\frac1{|\varphi'(0)|}
\geq\frac{367}{500}
>\frac1{\sqrt2},
$$
which is a contradiction.

Every eigenvalue of $A_{\mu_*,r_*}$ lies in the interior of its
numerical range. The complete Crouzeix ratio is continuous at such
matrices \cite[Theorem~2.1]{Crouzeix2016}. Since
$(\mu,r)\mapsto A_{\mu,r}$ is continuous, the strict inequality
persists in a neighborhood of $(\mu_*,r_*)$.
\end{proof}

For $0<t<1/2$, put
\begin{equation}\label{eq:A-t}
 A_t=A_{\mu(t),r(t)},
 \qquad
 \mu(t)=t\sqrt{1-2t},
 \qquad
 r(t)=(1-t)\sqrt t.
\end{equation}

\begin{proposition}\label{prop:endpoint}
One has
\begin{equation}\label{eq:endpoint-limits}
        \lim_{t\to1/2^-}\psi(A_t)=2,
        \qquad
        \lim_{t\to1/2^-}\psi_{\mathrm{cb}}(A_t)=2.
\end{equation}
In particular,
\[
        \sup_{0<t<1/2}\psi_{\mathrm{cb}}(A_t)=2.
\]
\end{proposition}

\begin{proof}
As $t\to1/2^-$, the matrices $A_t$ converge to
\[
 A_*=\begin{pmatrix}
 0&0&i/(2\sqrt2)\\
 0&0&i/(2\sqrt2)\\
 i/(2\sqrt2)&i/(2\sqrt2)&-1
 \end{pmatrix}.
\]
With respect to the orthonormal basis
\[
 \frac1{\sqrt2}(1,-1,0)^{\mathsf T},
 \qquad
 \frac1{\sqrt2}(1,1,0)^{\mathsf T},
 \qquad
 (0,0,1)^{\mathsf T},
\]
the matrix $A_*$ is unitarily equivalent to
\[
 [0]\oplus B,
 \qquad
 B=\begin{pmatrix}0&i/2\\ i/2&-1\end{pmatrix}.
\]
The matrix
\[
        N=B+\frac12I
        =\frac12\begin{pmatrix}1&i\\ i&-1\end{pmatrix}
\]
satisfies $N^2=0$ and $\norm{N}=1$. Hence $W(B)$ is the closed disk
with center $-1/2$ and radius $1/2$. This disk contains zero, so it is
also $W(A_*)$. For $p(z)=2z+1$,
\[
        \max_{z\in W(A_*)}\abs{p(z)}=1,
        \qquad
        \norm{p(A_*)}=2.
\]
The numerical range satisfies
$d_{\mathrm H}(W(X),W(Y))\leq\norm{X-Y}$. Therefore
\[
        \norm{p(A_t)}\longrightarrow2,
        \qquad
        \max_{z\in W(A_t)}\abs{p(z)}\longrightarrow1.
\]
It follows that $\liminf_{t\to1/2^-}\psi(A_t)\geq2$. The scalar
Crouzeix inequality gives the reverse bound. Finally,
\[
        \psi(A_t)\leq\psi_{\mathrm{cb}}(A_t)\leq2
\]
by Theorem~\ref{thm:intro-complete-three}, which proves the second
limit.
\end{proof}

\begin{corollary}\label{cor:endpoint-stability}
For $t$ sufficiently close to $1/2$, let $g_t$ be a normalized scalar
extremal function for $(A_t,W(A_t)^\circ)$, and choose a unit vector
$x_t$ such that
\[
        \norm{g_t(A_t)x_t}=\psi(A_t).
\]
Let $\mu_t$ be the representing measure associated with $g_t$ and
$x_t$, and put $\nu_t=(g_t)_*\mu_t$. There are a unitary matrix $U_t$
and a scalar $d_t$ such that
\[
 \norm{U_t^*g_t(A_t)U_t-
 \left(
 \begin{pmatrix}0&0\\2&0\end{pmatrix}\oplus[d_t]
 \right)}
 \leq8\sqrt{2-\psi(A_t)}.
\]
Moreover,
\[
 \left|\int_{\T}\zeta^k\dd\nu_t(\zeta)\right|
 \leq\vartheta(\psi(A_t))
 \qquad(k\in\mathbb Z\setminus\{0\}).
\]
Consequently, the matrix estimate tends to zero and
$\nu_t$ converges weakly to $m_{\T}$ as $t\to1/2^-$.
\end{corollary}

\begin{proof}
For $0<t<1/2$, the spectrum of $A_t$ lies in the interior of its
numerical range. Proposition~\ref{prop:domain-extremals} gives $g_t$ in
the form $B_t\circ\varphi_t$, where $B_t$ is a finite Blaschke product
and $\varphi_t:W(A_t)^\circ\to\D$ is conformal. Hence
$\norm{g_t}_{W(A_t)}=1$. Proposition~\ref{prop:endpoint} and
Corollary~\ref{cor:stability-Crouzeix} give the matrix estimate.
Theorem~\ref{thm:measure-stability}, with $K=s=\psi(A_t)$ and $m=1$,
gives the Fourier estimate and weak convergence.
\end{proof}

Corollary~\ref{v3:family-strict} gives
$\psi_{\mathrm{cb}}(A_t)<2$ for every $0<t<1/2$, while
Proposition~\ref{prop:endpoint} shows that these constants tend to two
as $t\to1/2^-$. Thus every point of the defective curve satisfies a
strict bound, but no constant smaller than two works for the whole
curve.

\section*{Acknowledgments}
We are deeply grateful to our friend Antti Per\"{a}l\"{a} for his invaluable discussions and insightful comments on this paper.

\medskip

The first-named author was partially funded by the Swedish Research
Council under grant agreement no.~2025-05053 and by the Swedish Energy
Agency under project no.~P2025-04323. The second-named author was also
partially funded by the National Science Centre, Poland, under the
Weave-UNISONO programme, grant no.~UMO-2025/07/Y/ST1/00146. The
third-named author was supported in part by the Engineering and
Physical Sciences Research Council grant EP/Y008375/1 and the Research Council of Finland grant 375631.

\end{document}